\documentclass{article}
\usepackage{amsfonts}
\usepackage[utf8]{inputenc}
\usepackage{amsmath}
\usepackage{amssymb}
\usepackage{graphicx}
\usepackage{authblk}
\usepackage{url}
\usepackage{algorithm2e}
\usepackage{appendix}
\usepackage{rotating}

\usepackage[english]{babel}
\usepackage[autostyle]{csquotes}

\vspace{-40pt}
\title{Cyclicality, Periodicity and the Topology of Time Series}
\vspace{-30pt}
\author[1]{Pawe{\l} D{\l}otko\thanks{\textbf{Corresponding Author}. Full Address: Mathematics Department, College of Science, Swansea University, Bay Campus, Swansea, SA1 8EN, United Kingdom. Tel: +44 (0)1792 606325. Email:p.t.dlotko@swansea.ac.uk. }}
\affil[1]{Mathematics Department, Swansea University, United Kingdom}
\vspace{-40pt}

\vspace{-30pt}
\author[2]{Wanling Qiu\thanks{Full Address: Accounting and Finance Subject Group, School of Management, University of Liverpool, 20 Chatham Street, Liverpool, L69 7ZH, United Kingdom. Email:wanling.qiu@liverpool.ac.uk}}
\affil[2]{School of Management, University of Liverpool, United Kingdom}
\vspace{-30pt}
\author[3]{Simon Rudkin \thanks{Full Address: Economics Department, School of Management, Swansea University, Bay Campus, Swansea, SA1 8EN, United Kingdom. Email:s.t.rudkin@swansea.ac.uk}}
\affil[3]{Economics Department, Swansea University, United Kingdom}
\vspace{-20pt}
\begin{document}
\maketitle
\begin{abstract}
Periodic and semi periodic patterns are very common in nature. In this paper we introduce a topological toolbox aiming in detecting and quantifying periodicity.
%as well as some more complicated patterns in time series.
The presented technique is of a general nature and may be employed wherever there is suspected cyclic behaviour in a time series with no trend.
%Discussing the topological signatures of periodicity we provide a workflow to reconstruct the general dynamics which generate the time series.
The approach is tested on a number of real-world examples enabling us to consistently demonstrate an ability to recognise periodic behaviour where conventional techniques fail to do so. Quicker to react to changes in time series behaviour, and with a high robustness to noise, the toolbox offers a powerful way to deeper understanding of time series dynamics.

%Forecasting major changes in stock market behaviour has inherent investor value, but reversal identification has relied on structural impositions such as technical analysis or log-periodicity. Indeed the process of identifying cyclical behaviour in the market has been an inexact science, challenged not least by the presence of economic, annual, monthly and weekly cycles within the same data. Exploiting the topology of time series data it is shown that multiple consecutive cyclical behaviours may be identified within a series without assumption of functional form or complex the estimation thereof. Taking the example of the S\&P 500 the ability to identify these concurrent cyclicalities in financial data, and hence potential for forecasting major behavioural change, is demonstrated.
\end{abstract}

\section{Introduction}
Topological Data Analysis (TDA) is a branch of data science aiming to apply methods of algebraic topology, in particular persistent homology~\cite{edelsbrunner2010computational} to analyze data. As a topological tool it is robust to additive noise and can be applied to extract information from small datasets. Present applications of TDA group around cross sectional analyses of multidimensioanl data~\cite{banman2018mind,Ferri,bobrowski2018topology} for example, and to a much lesser extent the properties of time series. Recent works on time series include \cite{perea_harer},\cite{perea2017sliding} and \cite{perea2019topological} theoretical driven work and the consideration of financial landscapes in \cite{gidea2018topological}. Methods developed in this strand are the foundation for our work but currently leave many questions unanswered. A need to understand more on the topology of time series, both continuing in the spirit of the previous pioneering literature and through the embracing of new approaches, such as that presented here, is identified.

Our specific motivation for the study of periodic functions comes from the prevelance of highly volatile time series
%, such as the daily temperature, sun activity and a number of time series
in many real world applications. \cite{smith1992estimating} and \cite{turchin1992complex} offer useful introductions to the features of noisy time series in finance and biology respectively. In each case the often high levels of variation present major problems for existing dynamic behaviour detection approaches. In turn this creates a call for a method with robustness to noise. We answer the call through the use topological and geometrical techniques to identify patterns, periodicity and cyclicality in noisy time series.

Contributions to the literature are made from three critical perspectives. Firstly we idnetify a number of cases where the ability of TDA to isolate periodic and cyclical behaviour is shown to be much more accurate than other methods suggested in the literature. Secondly, TDA does not require many repeats of the cycle in order to recognize its existence. Both in the experimental and theoretical studies we show that just two cycles is sufficient, unlike spectral analysis approaches that can often miss the true periodicity by a long way. Thirdly, all of this is achieved with a robustness to noise that few others can match. Across these contributions we develop a go to source for time series analysis. Two caution points are noted, however. Because adding trends imposes a model upon the data there is no detrending of the series prior to analysis. A seperate detrending could be undertaken first using existing packages to obtain a series that can then be analysed through our approach. Secondly, although all periodicity identifiers are indebted to their paramaterisation we do note there must be some thought given to the implementation of our methodology. Notwithstanding these small points of caution, this paper makes significant strides to truly unlocking our understanding of time series.

In a theoretical exposition we compare the ability of TDA to learn from a large set of noisy time series built from two known periods. We first outline the approach, depicting periodic functions as cycles, and combinations of functions as variations thereupon. From the preliminary exposition comes the evidence that cyclicality of the sliding window embedding is a necessary, but not sufficient, condition for the periodicity of the input function. Robustness of conclusions against the selection of sliding window embedding length provides a timely confirmation of the ability of the approach to identify critical features within the time series.

\begin{figure}
    \begin{center}
    \caption{Work Flow for Periodicity Identification}
    \includegraphics[scale=0.5]{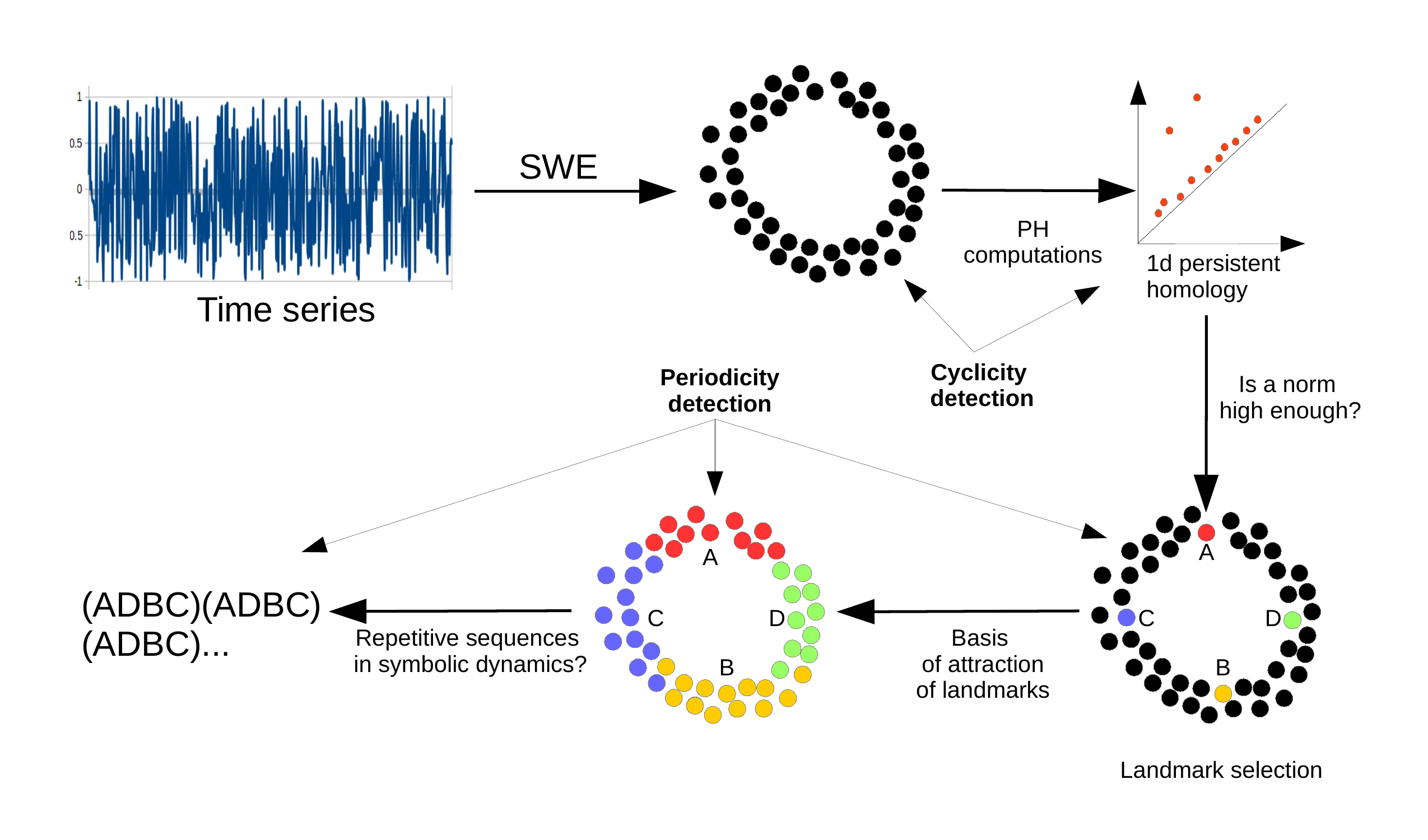}

    %\caption{The summary of the approach to time series analysis discussed in this paper; We start from the input time series (in the top right). By performing its sliding window embedding, we obtain a \emph{dynamic point cloud}. Given it, we use \emph{persistent homology} to detect long one dimensional cycles that gives information about \emph{cyclicity of the point cloud} (at top right). Given that, in order for test for periodicity, we select a uniformly spread subset of \emph{landmark points} (at bottom right) and cluster the whole point cloud by their proximity to the landmark points (bottom middle). Given the natural \emph{time ordering} of the points, we track the clusters to which following points belong. Once the clusters are tagged with letters of alphabet, as in the bottom row of the picture, we obtain a sequence of letters as in the bottom left. Discovering repetitive patterns therein is an evidence for periodicity of the initial time series.}
    \label{fig:paper_flow}
    \end{center}
\raggedright
\footnotesize{Notes: Illustration of process of the methodology for cyclicality and periodicity detection within time series. We start from a time series and construct its sliding window embedding (SWE). Subsequently a \emph{test of cyclicity} is performed. We compute the first persistent homology group of the SWE. If the norm of the diagram is high we conclude that the point cloud is \emph{cyclic}. Subsequently a \emph{test of periodicity} is performed. We select a number of equally spaced \emph{landmark} points. Each point within the SWE is assigned either to the closest (hard assignment), or to a number of closest (soft assignment) landmark point. Finally we consider points from the SWE in the order in which they appear. Each point output has the label of the landmark point assigned to it. The obtained sequence will be refereed to as \emph{symbolic dynamics}. Periodic symbolic dynamics are evidence of periodicity of the initial time series.}
\end{figure}

Our approach is neatly summarised in Figure \ref{fig:paper_flow}. Starting from a time series and moving to a sliding window embedding (SWE) is a common means to capture long, or multivariate, time series. In this paper we are looking at the dynamics of one variable and hence it is the former that is of interest. This paper focuses on univariate time series as an illustration, but the theory would extend to multivariate time series as well. Through this process, and the subsequent employment of persistent homology on the embedded point cloud, we are able to identify much more about the overall behaviour of the data generating process with a strong robustness to the noise inherently present.
%Where we refer to robustness to noise we draw on a fundamental understanding of TDA that data points are realisations from a surface and noise would see the realised points in different locations upon that surface.Should movements in the points not create changes in the homology then we can conclude that the approach is invariant to such noise.
We return to the noise question within the formal set out of the persistent homology approach.

Resulting from the SWE is a point cloud that may, if cyclicality is present, appear as a doughnut. If we were to number the points we would expect them to go round the doughnut and then return to the start point before going round again. To intuitively identify such behaviour we will spread a small number of points points $l_1,\ldots,l_k$ uniformly in the SWE of the time series. By assigning the following points  from SWE to the closest among $\{l_1,\ldots,l_k\}$, we obtain an approximation of dynamics of the points in SWE. In particular, a periodic point cloud will give us a repetitive pattern of symbols. From this simplistic interpretation of the embedding comes a powerful opportunity to identify not only cyclicality, but also the stability through time thereof as well as more complicated dynamical patterns.

%can label the start point of the time series as A. We then denote the point furthest from the start point as B and then continue to label points further from the already labelled ones, C, D etc... to create an impression of the shape of the doughnut. This results in the point B being opposite A, and C opposite D, in Figure \ref{fig:paper_flow} \textbf{ANSWERD I HOPE: I am not sure if it is clear that we are talking here about the landmark points construction. Could you please make sure?}. Should the series move from being close to D then close to C without going close to B then a change in behaviour has occurred. From this simplistic interpretation of the embedding comes a powerful opportunity to identify not only cyclicality, but also the stability through time thereof as well as more complicated dynamical patterns.

%Obtaining the period from the resulting SWE requires aggregation of the times of return. A graphical approach is demonstrated as a visualisation tool, with algorithms able to obtain all of this information in the background

Obtaining the periodicity from this is then simply a matter of looking at the times of return through the new sequence that has been generated from the SWE. Here direct comparison with extant methodologies may be made, our use of complex exemplars highlights how the approach adopted in this paper succeeds where others become distracted. Fundamentally by incorporating all data, and not seeking patterns as the driver for the quantification of periodicity, our proposed method is able to recover all constituent parts from the point cloud. As is later shown the applicability of this robustness in complex time series is strong; there are many composite functions expected behind the observed series.

We draw examples from temperature data, Sun activity and business sales statistics. As one of many geographic time series, temperature is well studied for its period; on a crude level we would expect a periodicity of around 1 year to reflect the changing seasons \cite{edsall2000tools}, and with a long enough dataset longer periodicities of global warming and cooling \cite{diodato2014climate}. Sun activity has a longer period, most estimates being around 130 months \cite{oliver1998emergence,lopes2015looking}. Finally businesses notoriously go through strong spells and weaker times . Business is also subject to the overriding seasonality of demand meaning an inherent periodicity. Such external influences may vary over the business cycle \cite{matas2004does}, demanding time series analytics can cope with varying periodicity robustly. In our applications we demonstrate the TDA approach can recover many of these features in ways that alternatives do not.

 The rest of the paper now proceeds as follows. Firstly we introduce the TDA approach, placing it within the context of prior time series work. Turning to the
 metrics that may be obtained from noisy time series Section \ref{sec:noisy_time_series} discusses the role of persistent homology in capturing dynamics. We then detail our approach to periodicity detection in Section \ref{sec:chop_and_loop}. Section \ref{sec:compare} looks at the evaluation of our artificial examples with code from the \textit{forecast} package for R \cite{hyndman2019package}. Applications to real world data are then offered in Section \ref{sec:apply}, again with comparison to the \cite{hyndman2019package} code. Finally, Section \ref{sec:conclude} concludes on the benefits evidenced and the potential for impact from our toolbox.

\section{Mathematical Summary}
\label{sec:math_approach}
A function $f : \mathbb{R} \rightarrow \mathbb{R}$, $f$ is \emph{periodic with a period $t$} if for every $x \in \mathbb{R}$, $f(x) = f(x+t)$. This classical definition of periodicity suffers from at least two serious drawbacks. Firstly, it requires exact knowledge of the period $t$. Secondly, given a noisy version of $f$, one cannot expect $f(x)$ to be \emph{precisely} aligned with $f(x+t)$, as there will inevitably be some shift due to the noise. Moreover, the shift will vary for different $x \in \mathbb{R}$ making the direct usage of the classic periodicity definition difficult.

To balance this effect in this paper we will explore an alternative characterization of periodicity that uses SWE (a.k.a. time delay embedding). Let us assume that the function $f$ is given as a time series $(x_0,f_0),\ldots,(x_k,f_k)$ so that $f_i  = f(x_i)$. %For simplicity we will assume that $x_i$'s are equidistant points.
Let us fix two positive integers $d$ and $n$. An SWE of $f$ at $x_i$ is a vector $[f_i , f_{i+d} , \ldots , f_{i+nd}]$. Note that this vector can be interpreted as a point in $\mathbb{R}^{n+1}$. Intuitively we may think about it as a encapsulation of a local shape of the time series $f$.
By applying this construction for all indices for which that make sense, the following sequence of $N$ points:
$[f_0 , f_d , \ldots$,$ f_{nd}], [f_1 , f_{1+d} , \ldots , f_{1+nd}]$, $\ldots, [f_{k-nd} , f_{k-(n-1)d} , \ldots , f_{k}]$ is obtained. This collection of points is a sliding window (time delay) embedding (SWE) of the function $f$ with the parameters $n$ and $d$.
%for a given N provided $Nd+n \leq k$.
For every point in the SWE, the subscript of its last component will be refereed to as the \emph{time index of the point}.

Under the assumption that $f$ is continuous and periodic, with a period $t$
%, assuming that $Nd > t$,
there exists an intermediate $i \in \{ 1, \ldots ,k \}$ such that $t - (x_i-x_0)$ is minimal. In this case, $[f_0 , f_d , \ldots , f_{nd}]$ and $[f_{i} , f_{i+d} , \ldots , f_{i+nd}]$ as well as any other point obtained for a multiplier of $i$ will be close to each other. If in addition under the assumption that there are no intermediate points $j$ between $0$ and $i$ such that $[f_0 , f_d , \ldots , f_{nd})]$ is close to $[f_{j} , f_{j+d} , \ldots , f_{j+nd}]$, then the obtained sequence of points will be periodically wrapping around a cosed one dimensional curve in $\mathbb{R}^{n+1}$. As an example of this construction let us consider a function $f(x) = sin(x)$, sampled every $0.1$ with $d = 5$, $n = 2$ and $N=200$. Such a function forms the basis of Figure\ref{fig:sin}.

\begin{table}[h!]
\caption{Illustrating Sliding Window Embedding}
\begin{tabular}{ c  c  }
\includegraphics[width=0.45\textwidth, height=60mm]{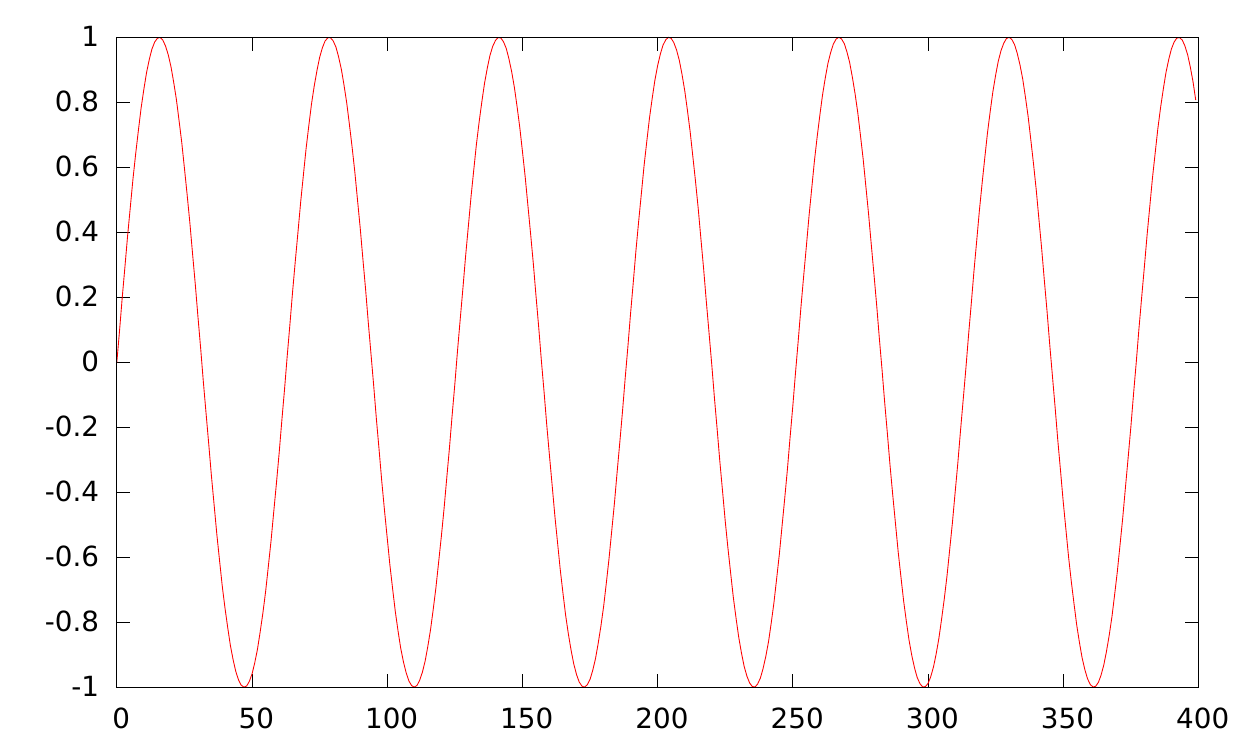}
&
\includegraphics[width=0.45\textwidth, height=60mm]{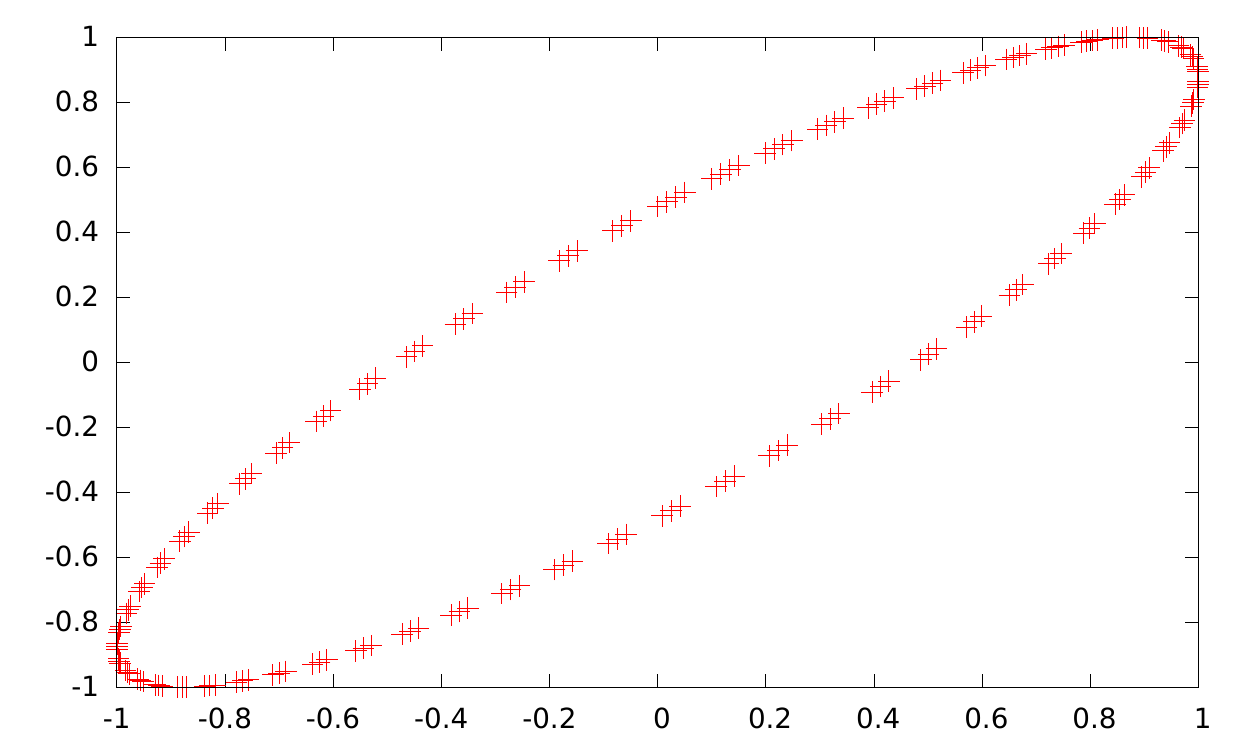}\\
Panel (a): Sin function & Panel (b): Sliding window embedding
\\
\end{tabular}
\label{fig:sin}
\raggedright
\footnotesize{Panel (a): a time series representing a sinus function sampled on a uniform grid of of a size $0.1$ (therefore the period of the time series correspond to $10\ 2\ \pi \approx 63$). Panel (b): a sliding window embedding of the time series with the sliding window size (embedding dimension) $2$ and the parameter $d = 5$. Note that periodic time series have cyclic sliding window embedding.}
\end{table}

\begin{figure}[h!]
\begin{center}
\caption{Sliding Window Embedding with Composite Functions \label{fig:cyclic_not_periodic}}
\begin{tabular}{ p{6cm}  p{6cm}  }
\includegraphics[width=0.45\textwidth, height=60mm]{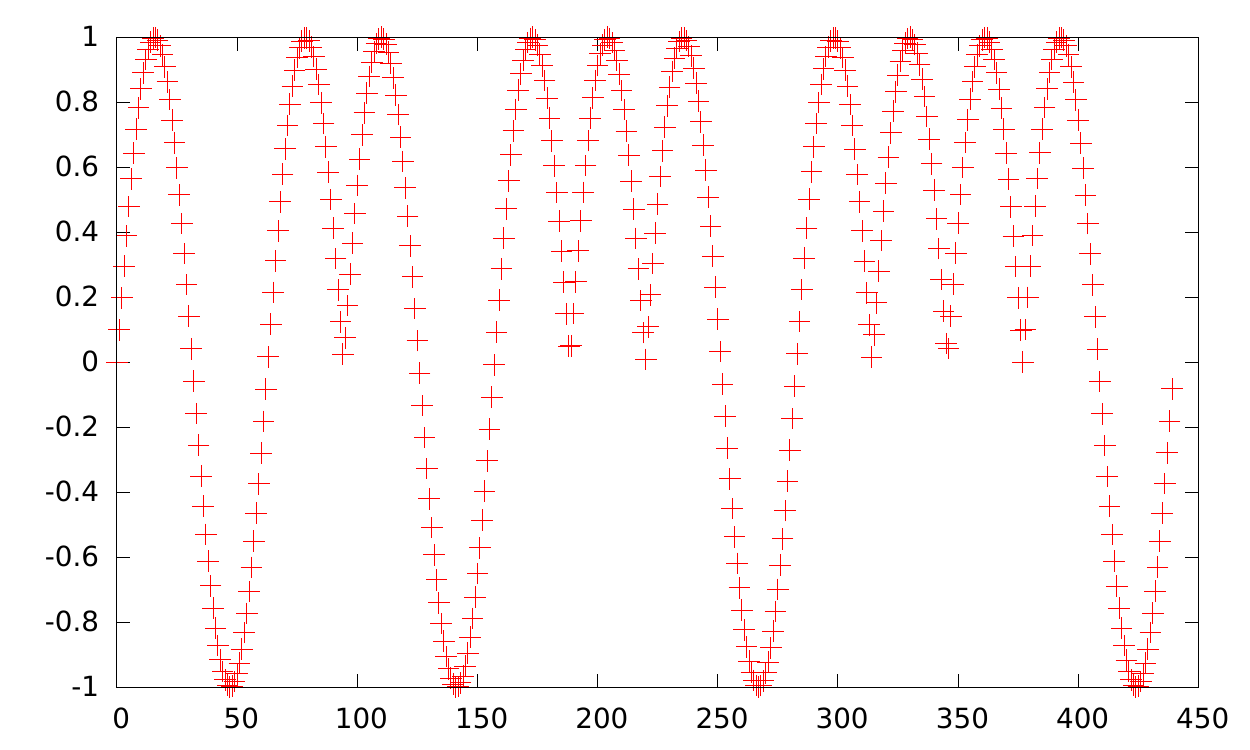}
&
\includegraphics[width=0.45\textwidth, height=60mm]{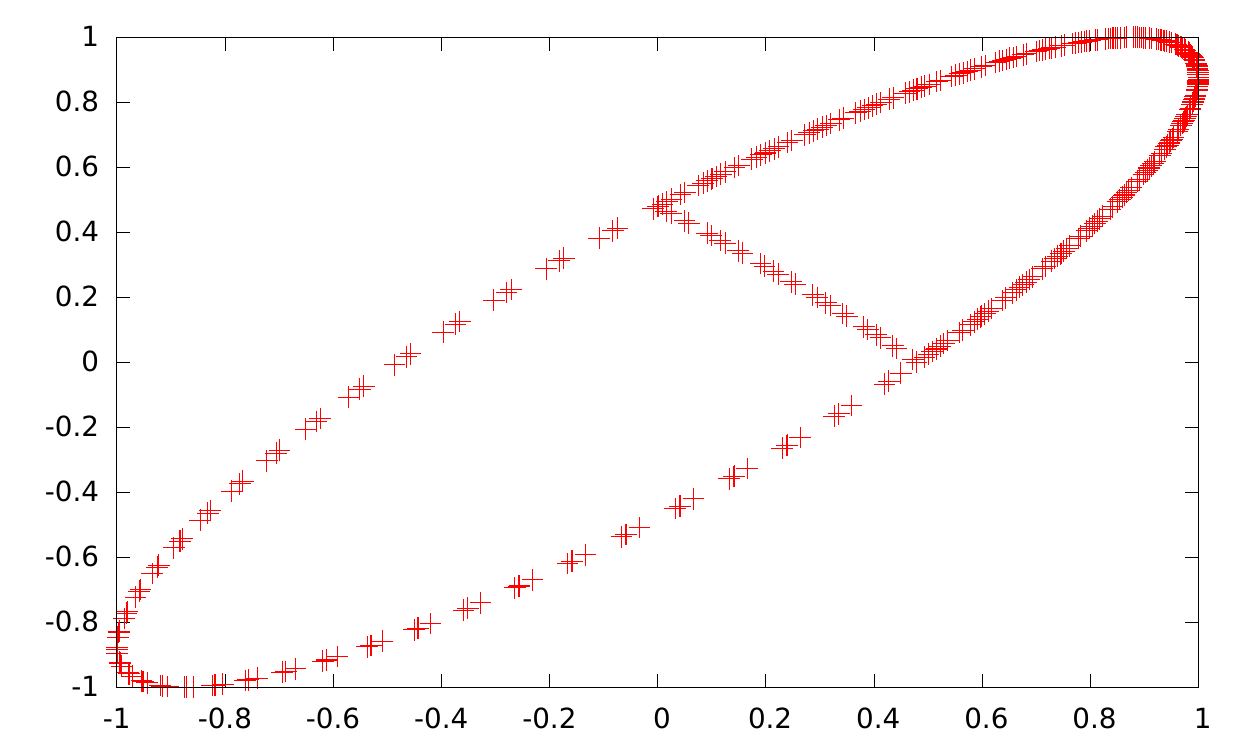}\\
Panel (a): Sin--based non periodic time series & Panel (b): Sliding window embedding\\
\end{tabular}
\end{center}
\raggedright
\footnotesize{Panel (a): The time series obtained by taking a positive part of a plot of a sinus, followed by negative part of the plot, followed by two positive pars of the plot, followed by a negative part of the plot, followed by three positive parts of the plot, followed by negative part of the plot etc. Note that this is not a periodic time series.
Panel (b): The sliding window embedding to two dimensional space with the parameter $d=5$, of the time series on the left. Note the existence of clear cycles in the sliding window embedding. This example shows that the cyclicity of the sliding window embedding is a necessary but not a sufficient condition for the periodicity of the time series.
}
\end{figure}

Figure \ref{fig:sin} illustrates that the periodicity of the function in panel (a) is clearly reflected in the cyclicality of its' SWE in panel (b). However, the reverse statement is not true. Let us consider a function that is obtained from a standard sinus by reflecting some of its parts by the x axis. For simplicity we will call a part of the plot when sinus is positive a \emph{positive part}, and a part of the plot where it is negative, a \emph{negative part}. We propose to construct a function $f$ by taking the positive part of sinus, followed by its negative part, followed by two constitutive copies of its positive parts, followed by a copy of its negative part, followed by three copies of its positive part positive parts, followed by a copy of its negative part negative part, followed by four copies of its positive parts etc. The initial part of the plot is presented in panel (a) of Figure~\ref{fig:cyclic_not_periodic}. Panel (b) depicts the obtained SWE which results into two cycles. The intermediate part, the line cutting through the cycle, together with the upper right part, is generated by the transition from a positive part of the sinus to another positive part\footnote{An example of this transition occurs at 100.}. The lower left part of the cycle in Figure \ref{fig:cyclic_not_periodic} panel (b) corresponds to the negative parts of sinus.

Despite clear cyclic behavior in Figure \ref{fig:cyclic_not_periodic}, the function $f$ is not periodic. This can be also observed on the SWE of $f$, by tracking the points forming in the SWE in the order of their appearance. While in the case of Figure~\ref{fig:sin}, they always make the whole turns, for Figure~\ref{fig:cyclic_not_periodic} they will  make a full turn only for transition from a positive to a negative domain. In case of a transition from a positive domain to another positive domain, the points will be moving through the intermediate bar and iterate in the upper half of a cycle.

Given the two examples above, we can conclude that there are two essential fingerprints of a periodic function that can be recognized from its SWE:
\begin{enumerate}
\item Points are aligned along a closed curve $C$ (a deformed cycle).
\item Given any point $x \in C$, the points from the sliding window embedding get back close to $x$ with more or less regular intervals.
\end{enumerate}
Here it should be noted that point 1 without point 2 would not represent periodic behaviour as noted in the example in the Figure~\ref{fig:cyclic_not_periodic}.

\begin{figure}[h!]
\begin{center}
  \caption{Illustration of Persistent Homology}
  \label{fig:persistence_explanation}
  \includegraphics[width=\linewidth]{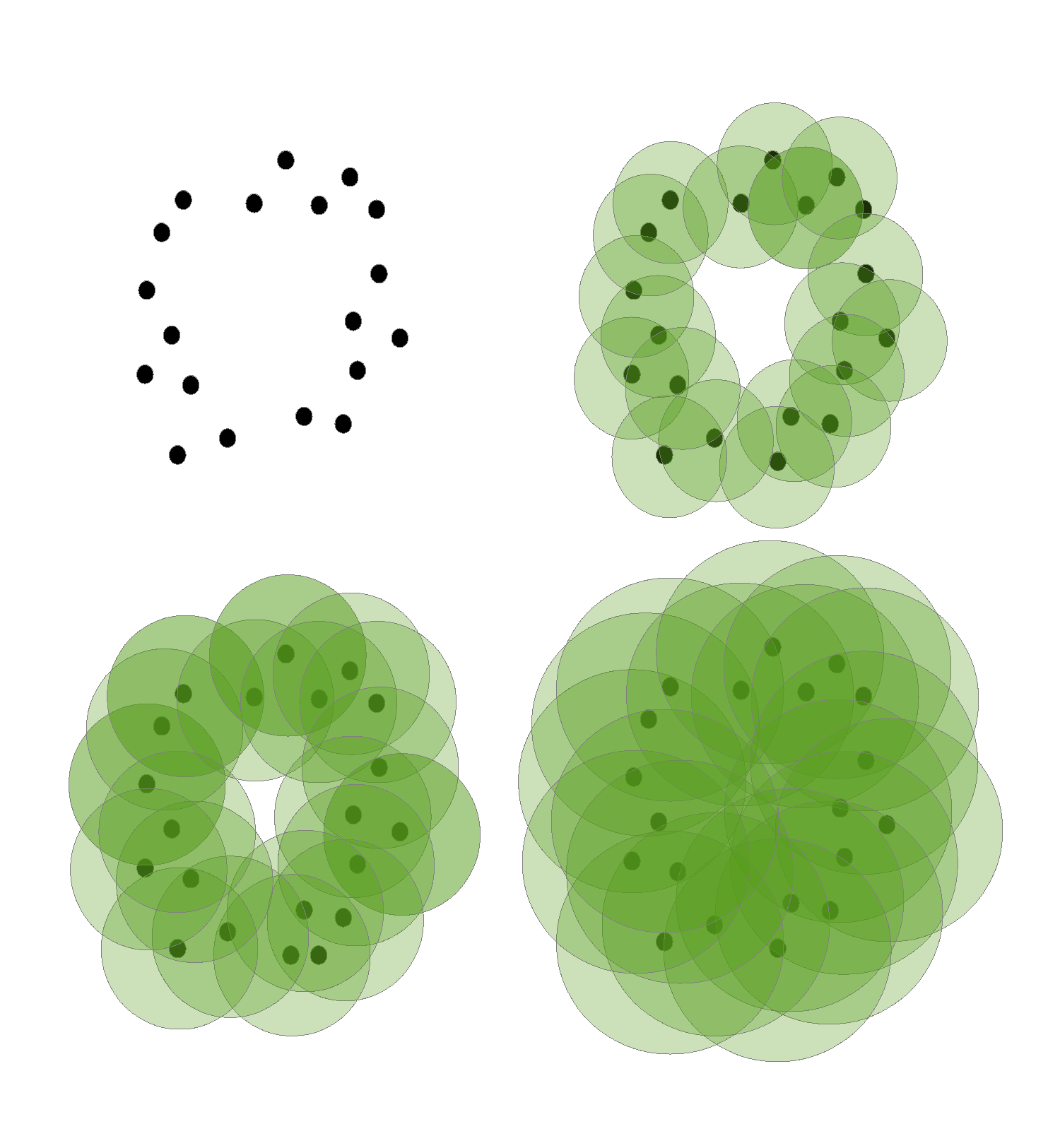}
 \end{center}
 \raggedright
 \footnotesize{Illustration of persistent homology. Consider the collection of points in the upper left. It is natural to state that it is sampled from a circle. Persistent homology is a tool to quantify this statement; Let us place a ball of a radius $r$ centered at every point in the point cloud and let $r$ grow from $0$ to $\infty$. The balls corresponding to increasing sequence of radii are present in the top right and the bottom line. For a large range of radii the sum of balls can be continuously deformed to a cycle. The persistent homology will capture this information by providing the range of radii for which the cycle can be observed.}
\end{figure}

Point (1) implies that the cyclicity of the sliding window embedding is a necessary, but not sufficient condition for the periodicity of the input function. In this paper we will use the machinery of persistent homology (PH) as per \cite{edelsbrunner2010computational}, to detect cyclic collections of points in the SWE of the time series. The idea of the construction is as follows: starting from a collection of points $\{p_1,\ldots,p_k\} \subset \mathbb{R}^{n+1}$ (i.e. each point has $n+1$ coordinates), we place a ball centered at each of the $p_i$ which has radius $r$, and we let $r$ vary between $0$ and $\infty$. Figure~\ref{fig:persistence_explanation} introduces the idea.

Given the obtained nested sequence of spaces, persistent homology will quantify that the obtained space is equivalent to a cycle for a large collection of constitutive radii. More precisely, PH will separate out a collection of \emph{persistence intervals}. Each interval is characterized by two points, $[b,d]$. The point $b$ indicates the radius on which the cycle is firstly observed, while $d$ indicates the radius in which the cycle is glued in to other cycles. For SWEs of (semi) periodic functions, we expect to see a single dominant persistence interval. Such a dominant interval would be accompanied by a number of short intervals that correspond to some short cycles resulting from the finite sampling.

\begin{figure}
\begin{center}
    \caption{Noisy Sin Function Example \label{fig:noisy_sin_fun} }
    \includegraphics[width=\linewidth]{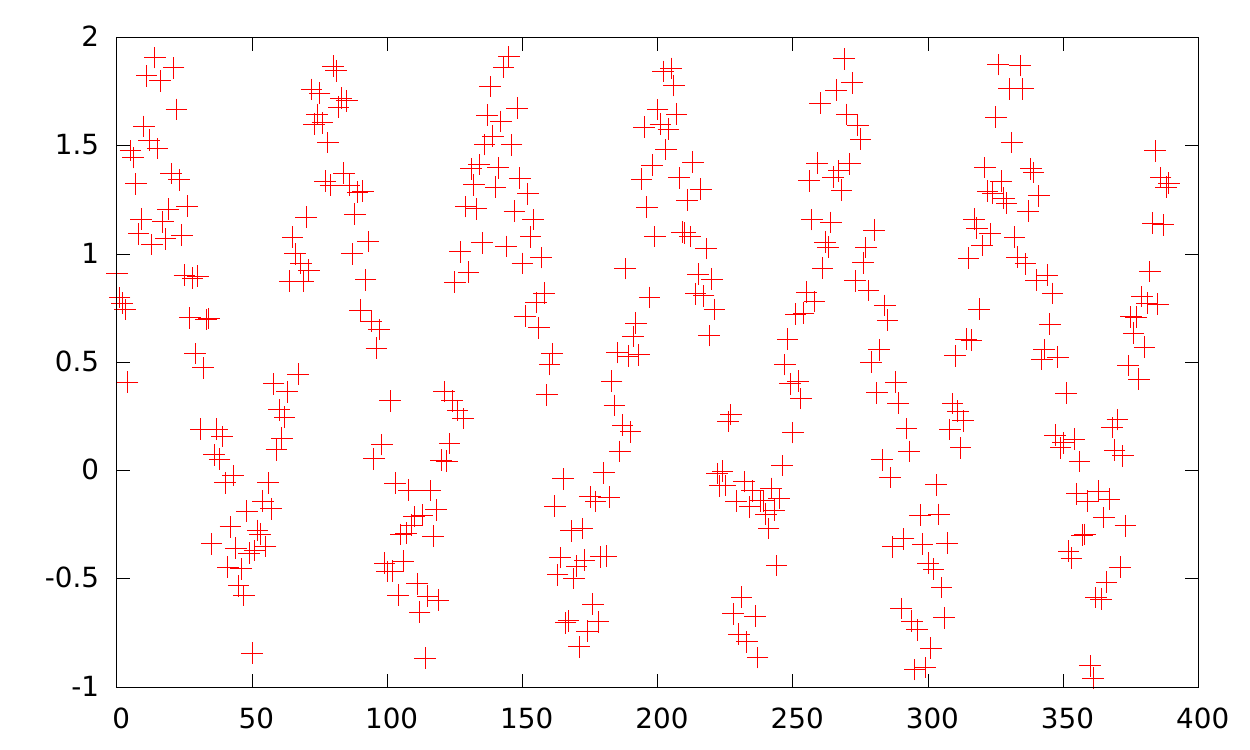}
\end{center}
\raggedright
\footnotesize{Notes: A sinusoidal time series with a uniform noise sampled from $[0,\frac{1}{2}]$. Points are plotted in lieu of a line to illustrate the point based nature of the information fed into the TDA approach developed in this paper.}
\end{figure}

There are a number of software libraries available to compute persistent homology. This paper will use a modified version of Gudhi library \cite{gudhi}. When constructing examples we will consider SWE into spaces of various dimensions, as the cyclicity of the collection of points may become more prominent as the dimension of the SWE grows. To illustrate this phenomena we have used the noisy version of sinus function presented in Figure~\ref{fig:noisy_sin_fun} and have constructed the SWE to dimension $2,5,10,15,20,25,30$ with $d = 10$. Table \ref{tab:domint} reports the lengths of three dominant persistence intervals.

\begin{table}
\begin{center}
    \caption{Sliding Window Embedding Dimensions and Dominant Persistence Intervals}
    \label{tab:domint}
    \begin{tabular}{|c|c|c|c|c|c|c|c|}
  \hline
  dim & 2 & 5 & 10 & 15 & 20 & 25 &30\\
  \hline
  1st & 0.3093 & 1.1748 & 1.7920 & 2.7821 & 3.3496 & 3.9754 & 4.2958 \\
  \hline
  2nd & 0.2356 & 0.2068 & 0.2767 & 0.2791 & 0.2406 & 0.2695 & 0.2587 \\
  \hline
  3rd & 0.1402 & 0.2003 & 0.2644 &  0.2459 & 0.2225 & 0.2535 & 0.2568\\
  \hline
\end{tabular}
\end{center}
\raggedright
\footnotesize{Notes: Estimation of dominant persistence is performed using the Gudhi library \cite{gudhi}}
\end{table}

As one can observe the ration of the longest to the second and third longest persistence interval is growing as the dimension grows. Clearly this phenomena does not extend for too long, and for very high embedding dimensions we see a reversed trend. These patterns demonstrate the importance of trying a number of parameters for the SWE to get most stable results. No algorithm exists for defining the optimal parameter choices. Rather like variable selection in multiple regression, it is left to the individual researcher to define what ultimately constitutes the optimal embedding dimension.

The relation of a dominant interval in dimension one and (semi) periodicity is known in the literature, see for example~\cite{perea_harer}, \cite{hamza} and others. However, such must be extracted from the noise of the series Section~\ref{sec:noisy_time_series} addresses this issue by discussing topological fingerprints of noisy time series.

\section{Topological Fingerprints of Noisy Time Series}
\label{sec:noisy_time_series}
Time series from finance, biology or meteorology are typically very noisy. Indeed under the efficient market hypothesis all deviations from the expected price may be regarded as noise with an expected value of 0. In this section we will study the effect of noise on the topological signatures of periodicity. Rather than providing a theoretical study that would require putting strong assumptions on the model of noise, we will present an example of a sinusoidal time series blurred with various level of a uniform noise. In the examples noise will be drawn randomly from a uniform distribution $\left[0,N_X\right]$, where the maximal noise level $N_X \in {0,1,2,3,4}$.  For each of those time series, we will compute its one dimensional persistence diagram and concentrate on the dominant intervals.

\begin{figure}
    \begin{center}
    \caption{Noisy Sin Functions}
    \label{fig:noisy_sin}
    \begin{tabular}{ c  c  }
\includegraphics[width=0.45\textwidth, height=30mm]{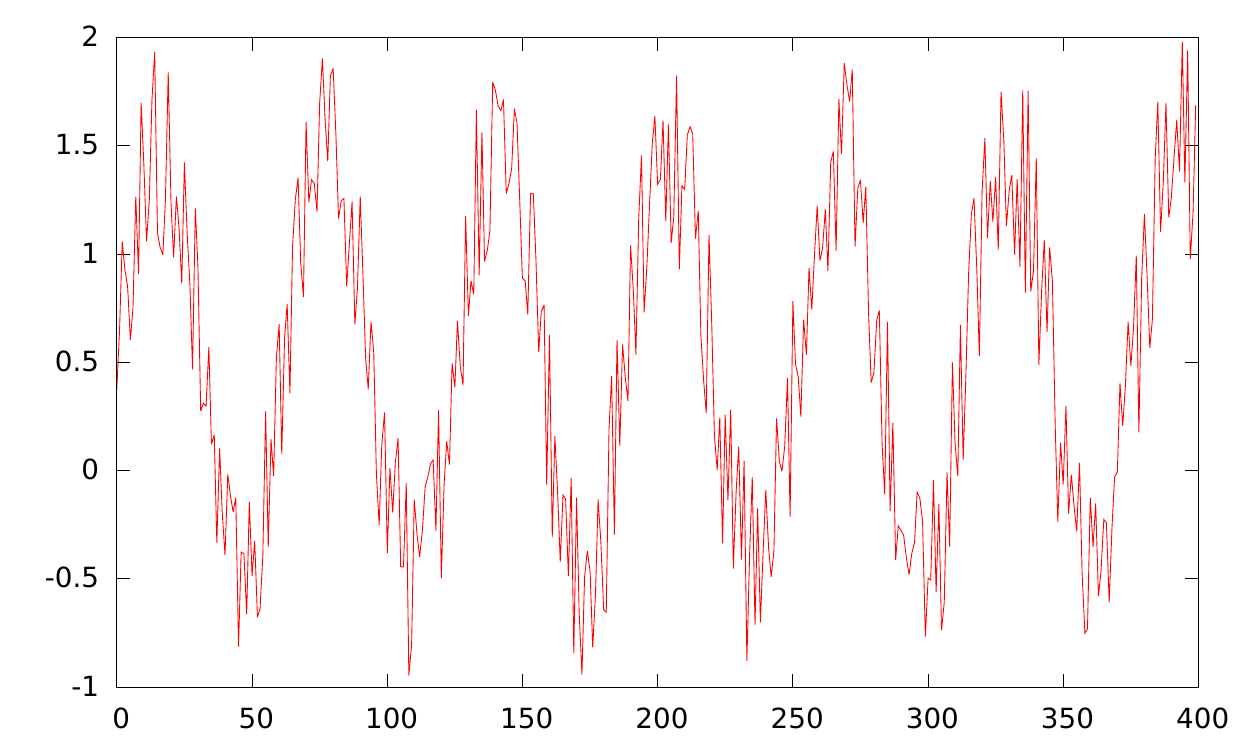}
&
\includegraphics[width=0.45\textwidth, height=30mm]{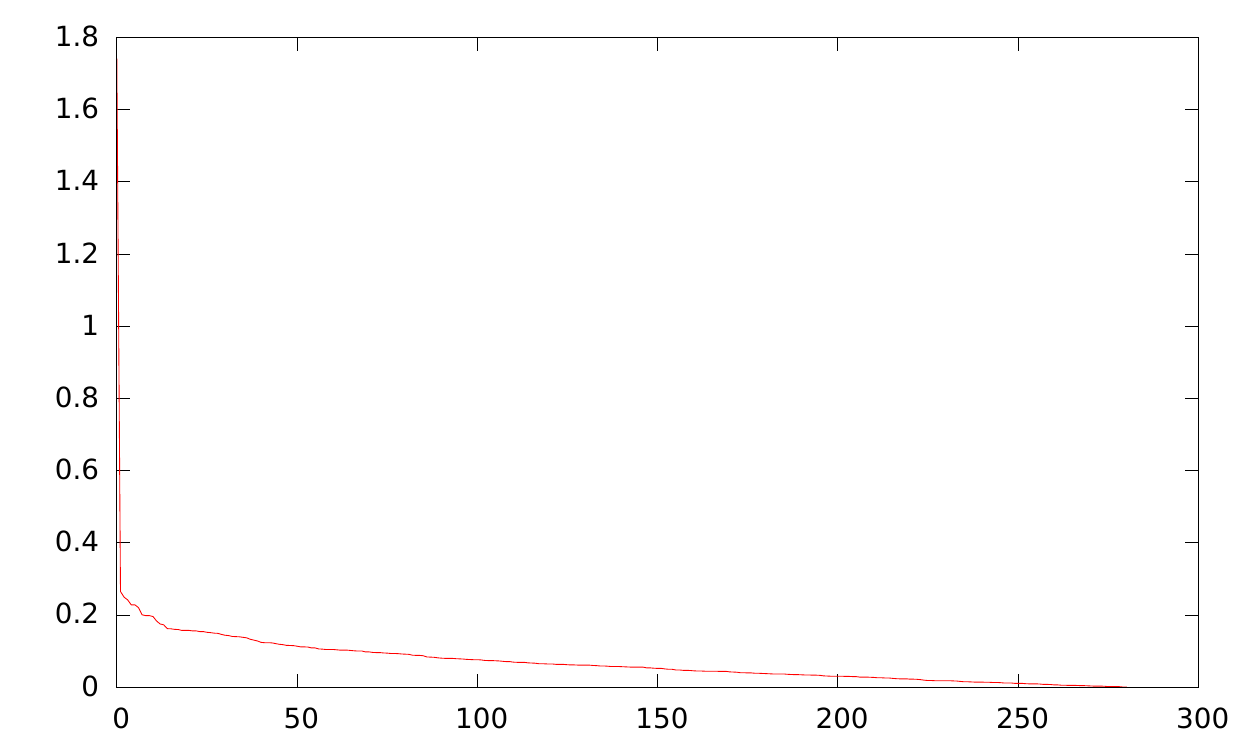}\\
Panel (a): Noise = 1 & Panel (b): Dominant interval length \\
\includegraphics[width=0.45\textwidth, height=30mm]{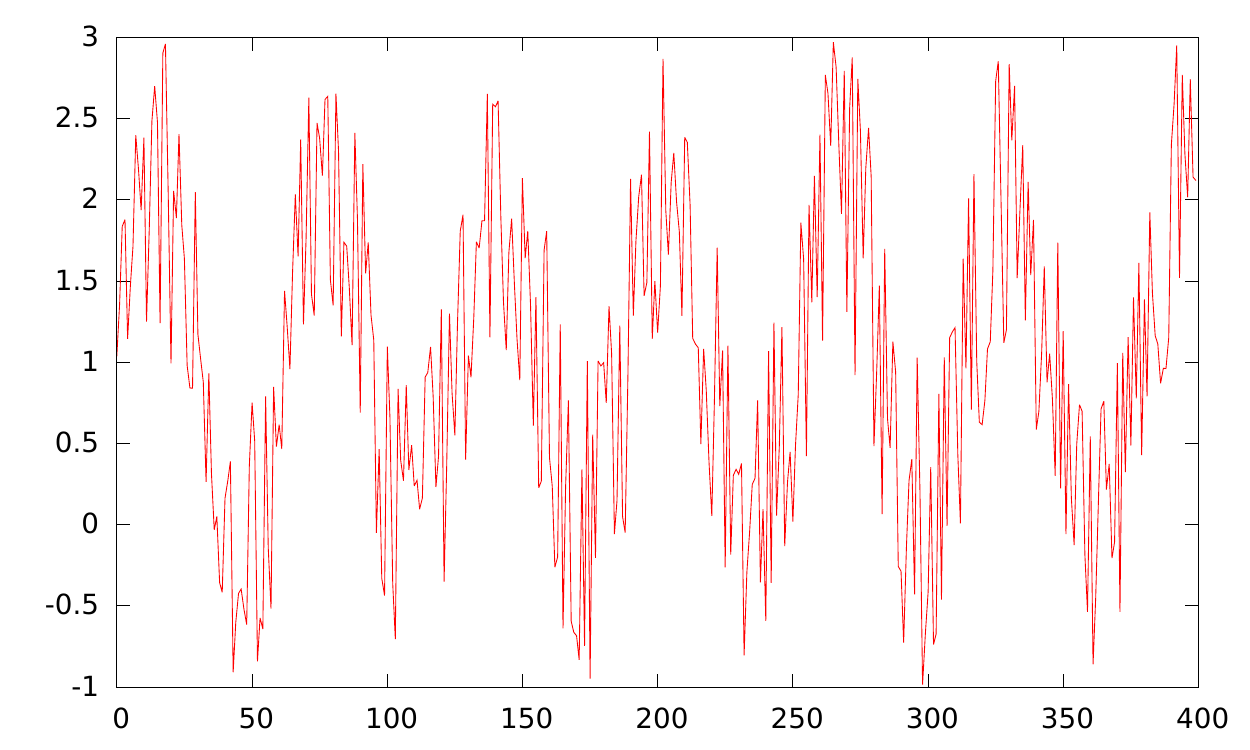}
&
\includegraphics[width=0.45\textwidth, height=30mm]{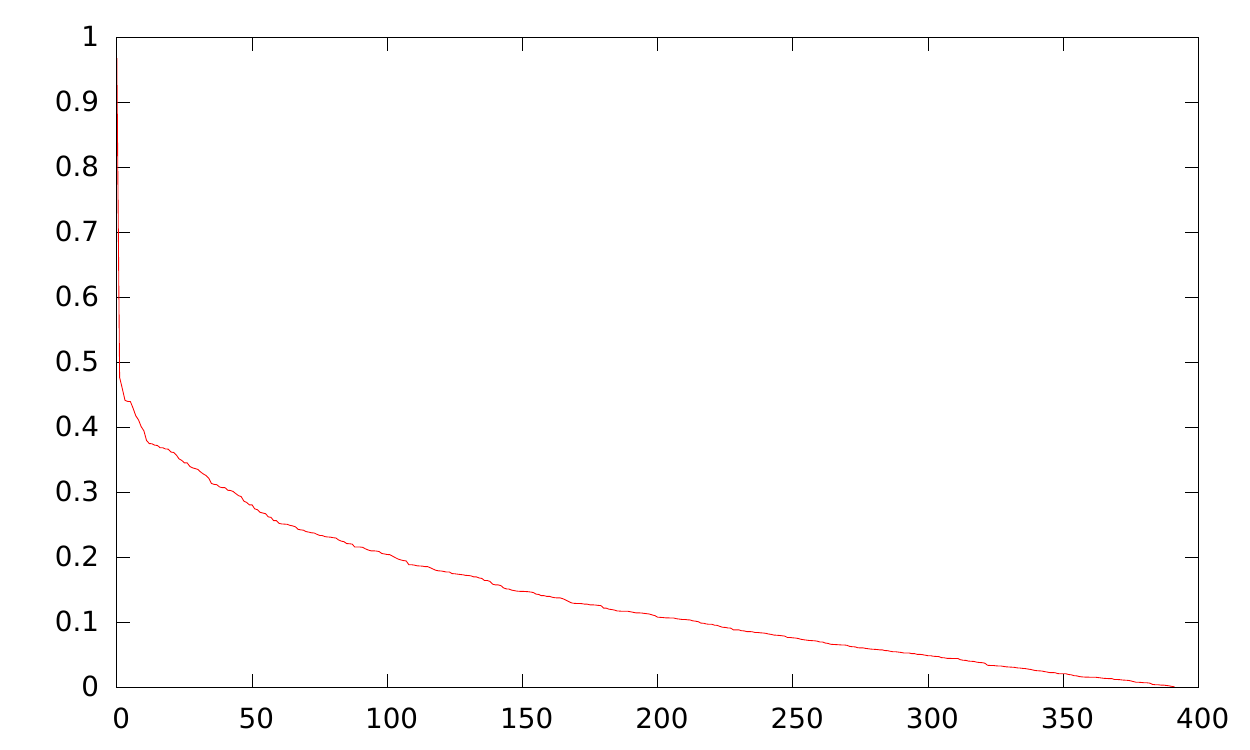}\\
Panel (c): Noise = 2 & Panel (d): Dominant interval length \\
\includegraphics[width=0.45\textwidth, height=30mm]{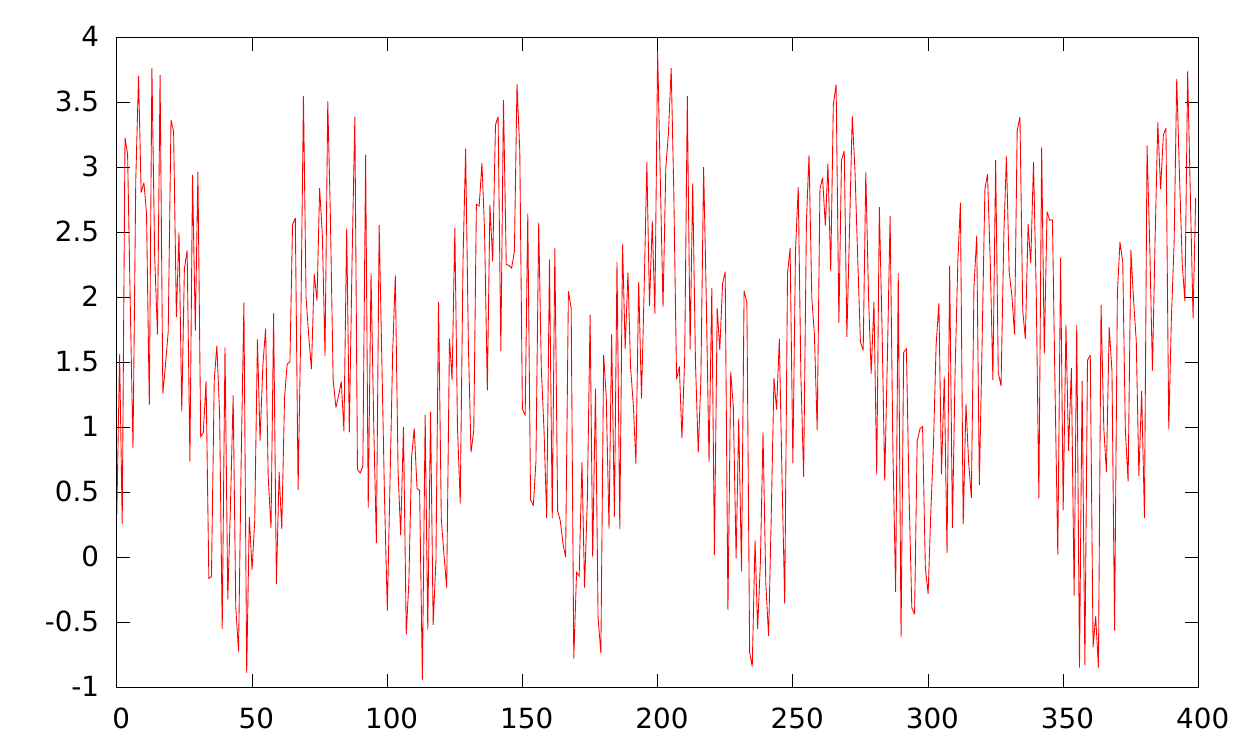}
&
\includegraphics[width=0.45\textwidth, height=30mm]{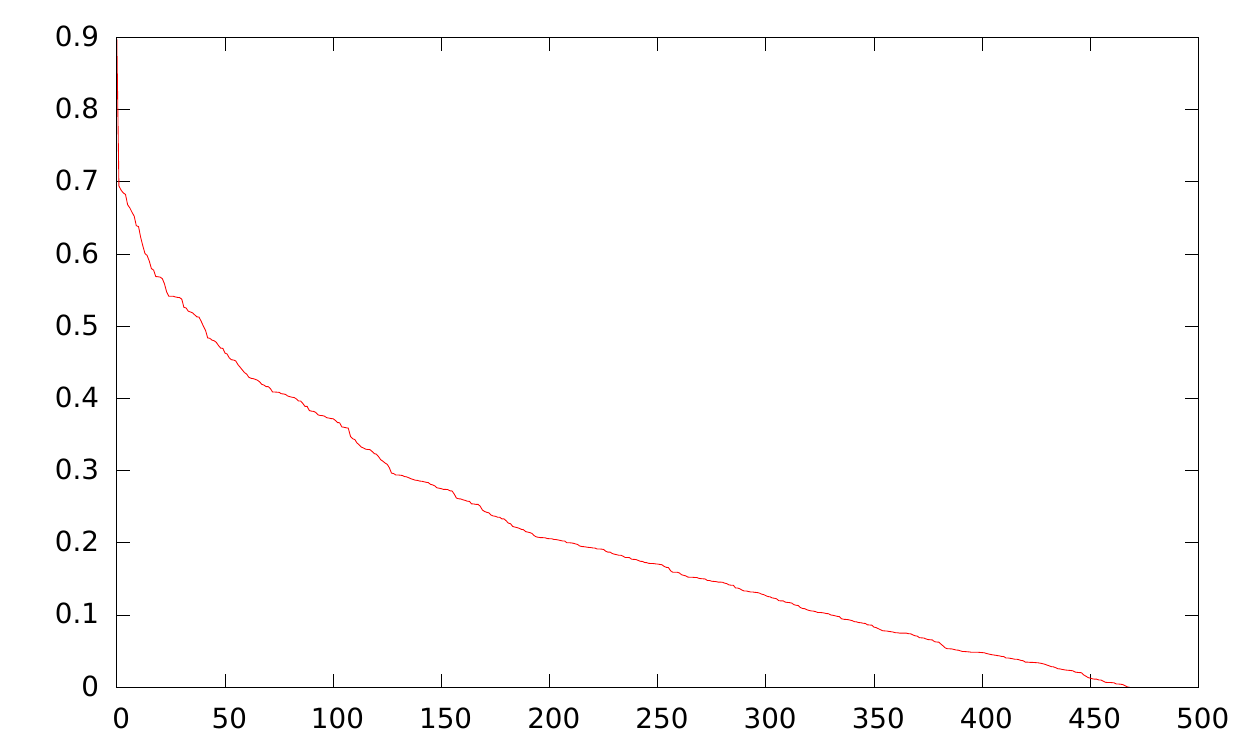}\\
Panel (e): Noise = 3 & Panel (f): Dominant interval length \\
\includegraphics[width=0.45\textwidth, height=30mm]{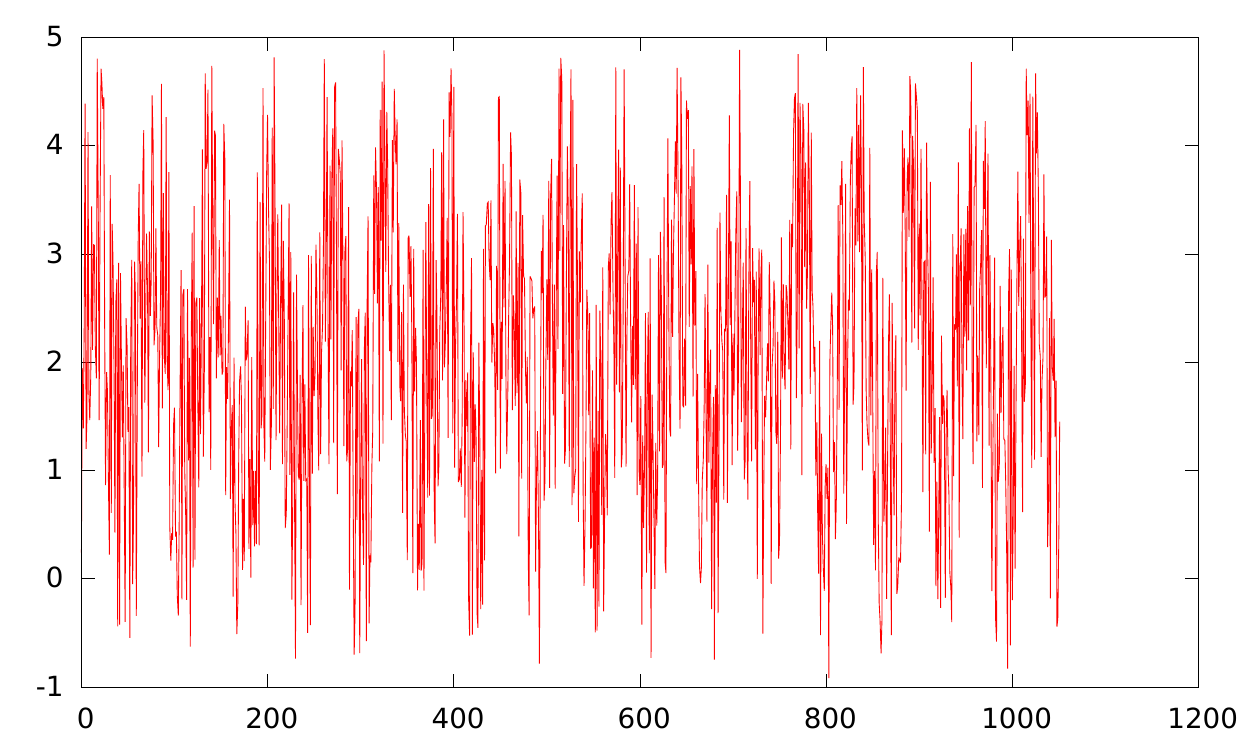}
&
\includegraphics[width=0.45\textwidth, height=30mm]{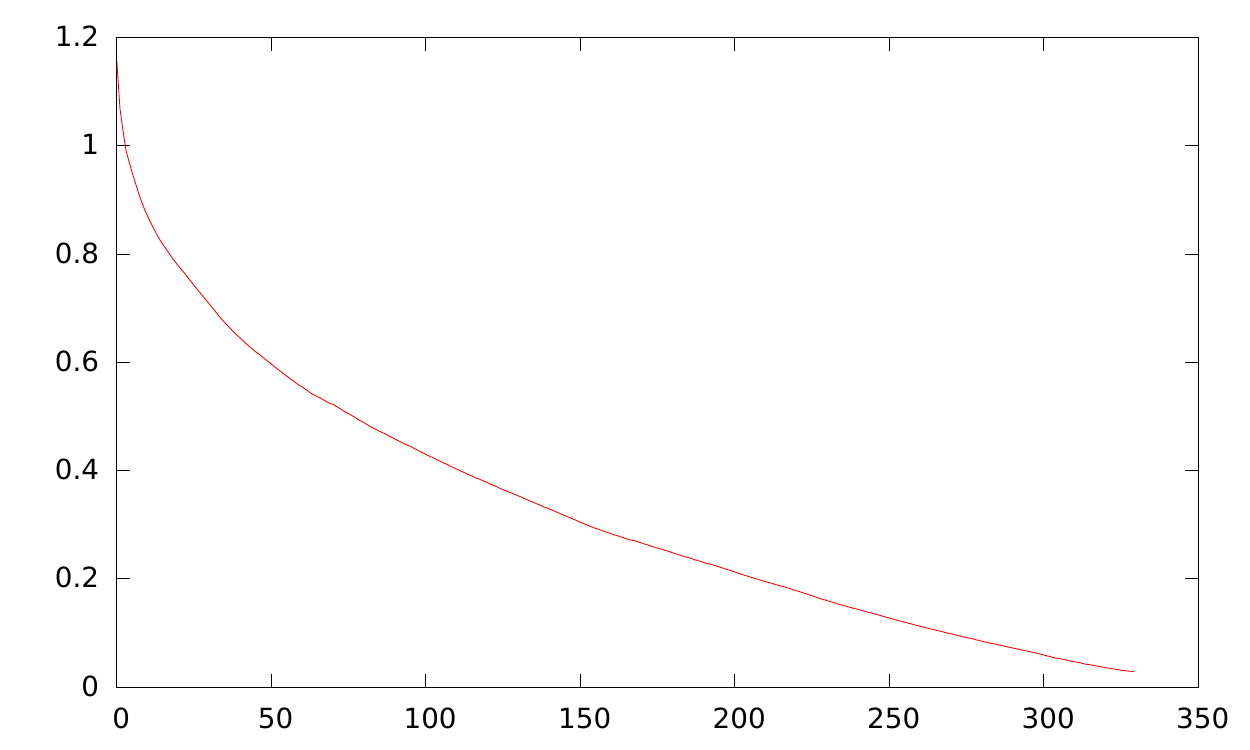}\\
Panel (g): Noise = 4 & Panel (h): Dominant interval length \\
\end{tabular}
\end{center}
\raggedright
\footnotesize{Notes: Panels (a), (c), and (e) provide illustrations of sin functions with increasing levels of noise. In each case we plot the corresponding dominant interval length, death minus birth, for the persistence intervals derived from the application of persistent homology to the sliding window embedding of the time series. The plots of the right present the histograms of the lengths of persistence intervals (from the longest to the shortest) in dimension one. Note that for the small level of noise (top figure) there is clearly dominant interval. With the increasing level of noise (second and third row) the shape of the histogram significantly changes.}
\end{figure}

\begin{figure}
    \begin{center}
    \caption{Noisy $sin(2x)+sin(\frac{x}{2})$ Functions \label{fig:noisy_sin3}}
    \begin{tabular}{ c  c  }
\includegraphics[width=0.45\textwidth, height=30mm]{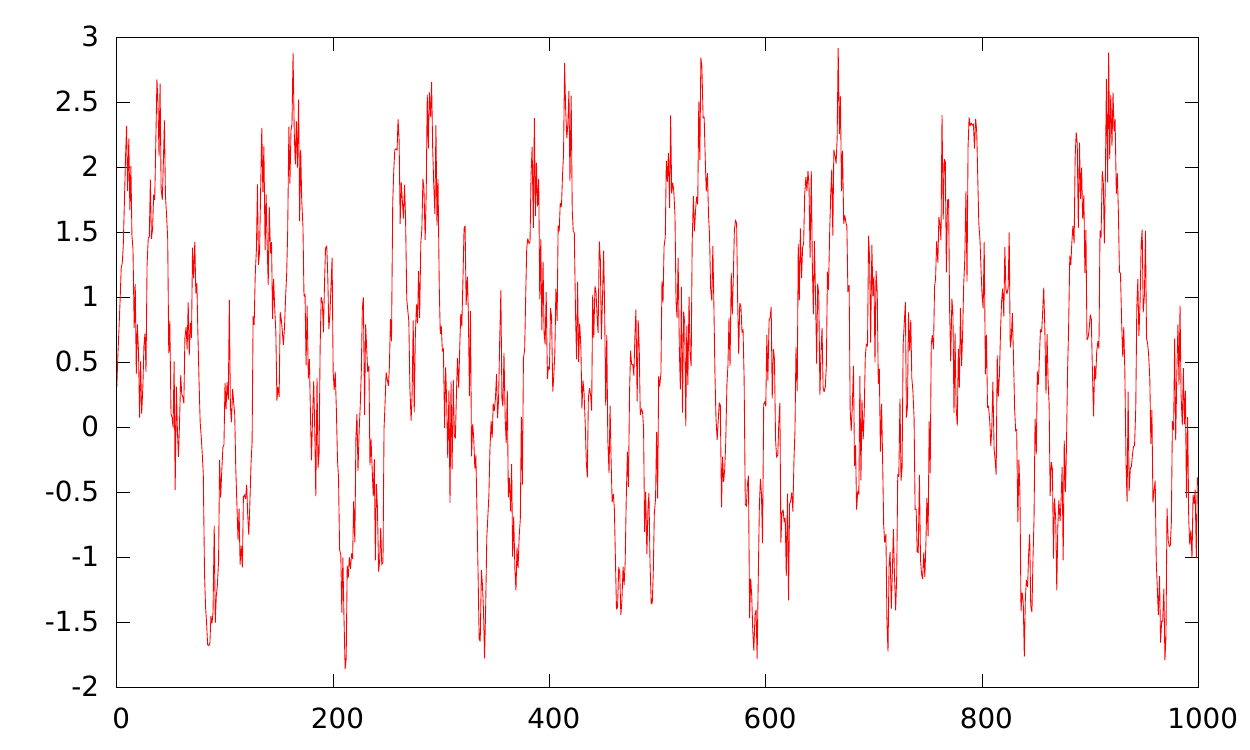}
&
\includegraphics[width=0.45\textwidth, height=30mm]{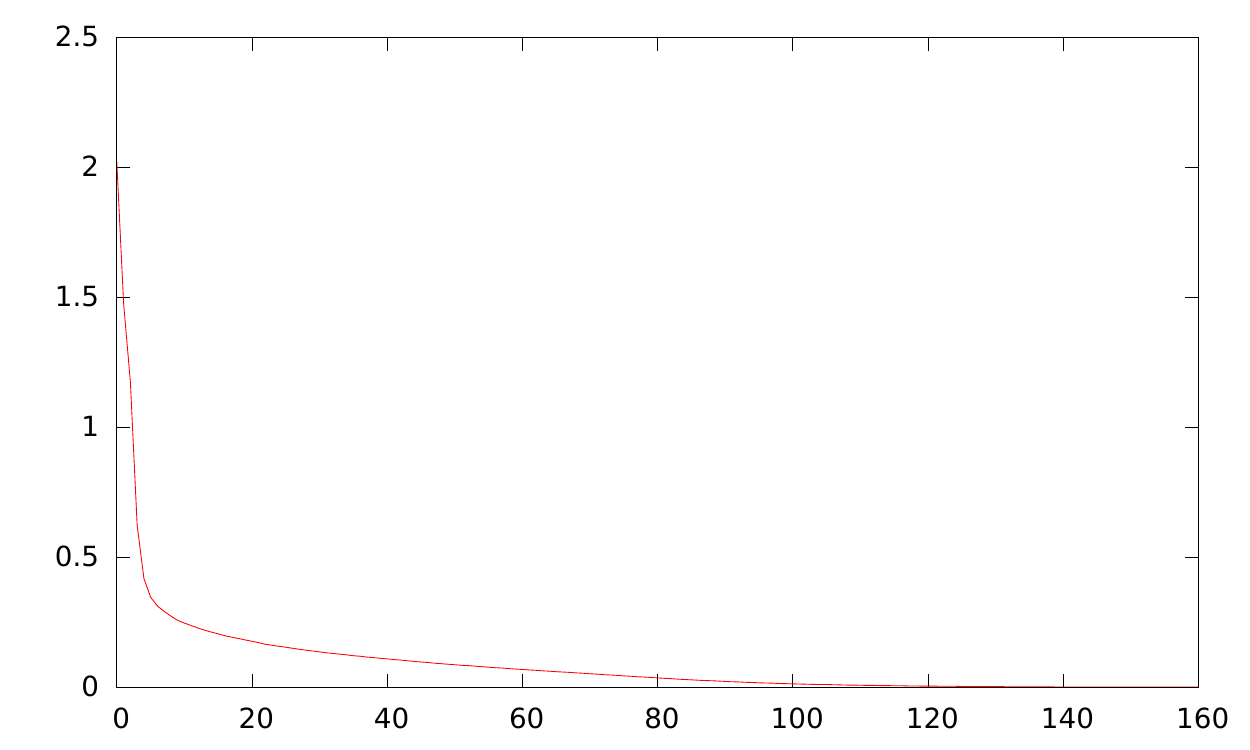}\\
Panel (a): Noise = 1 & Panel (b): Dominant interval length \\
\includegraphics[width=0.45\textwidth, height=30mm]{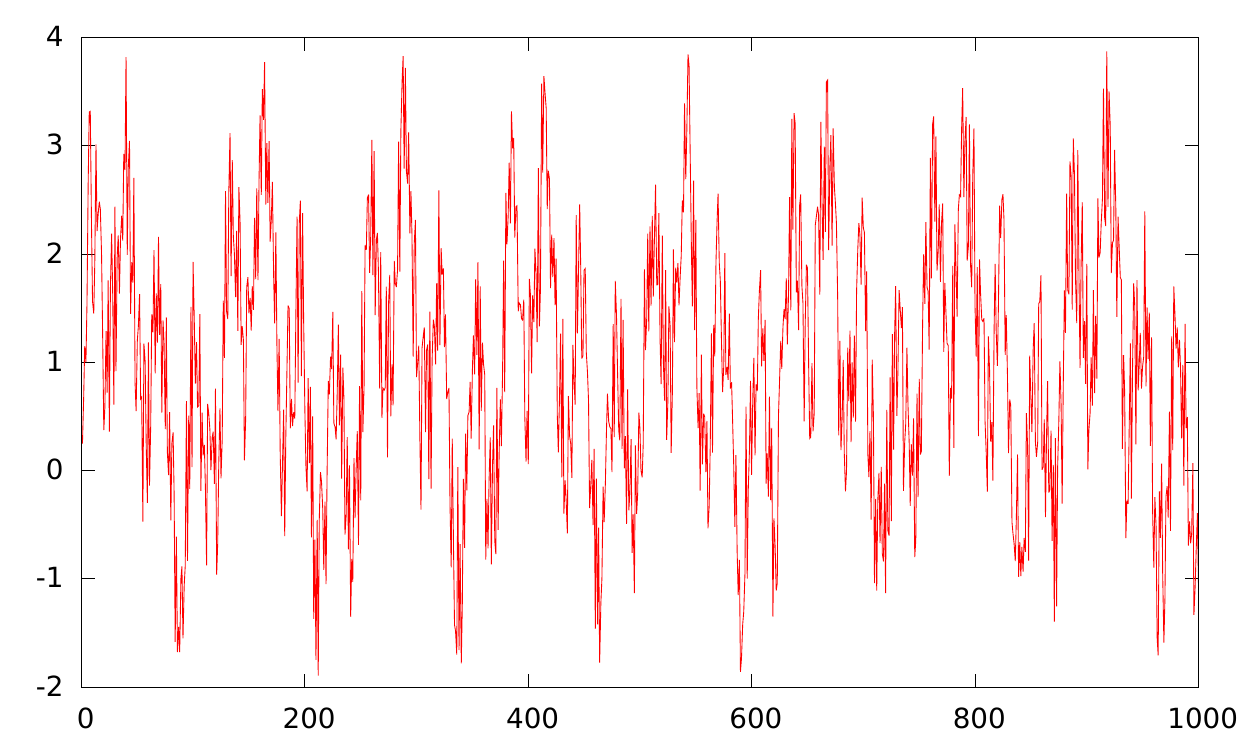}
&
\includegraphics[width=0.45\textwidth, height=30mm]{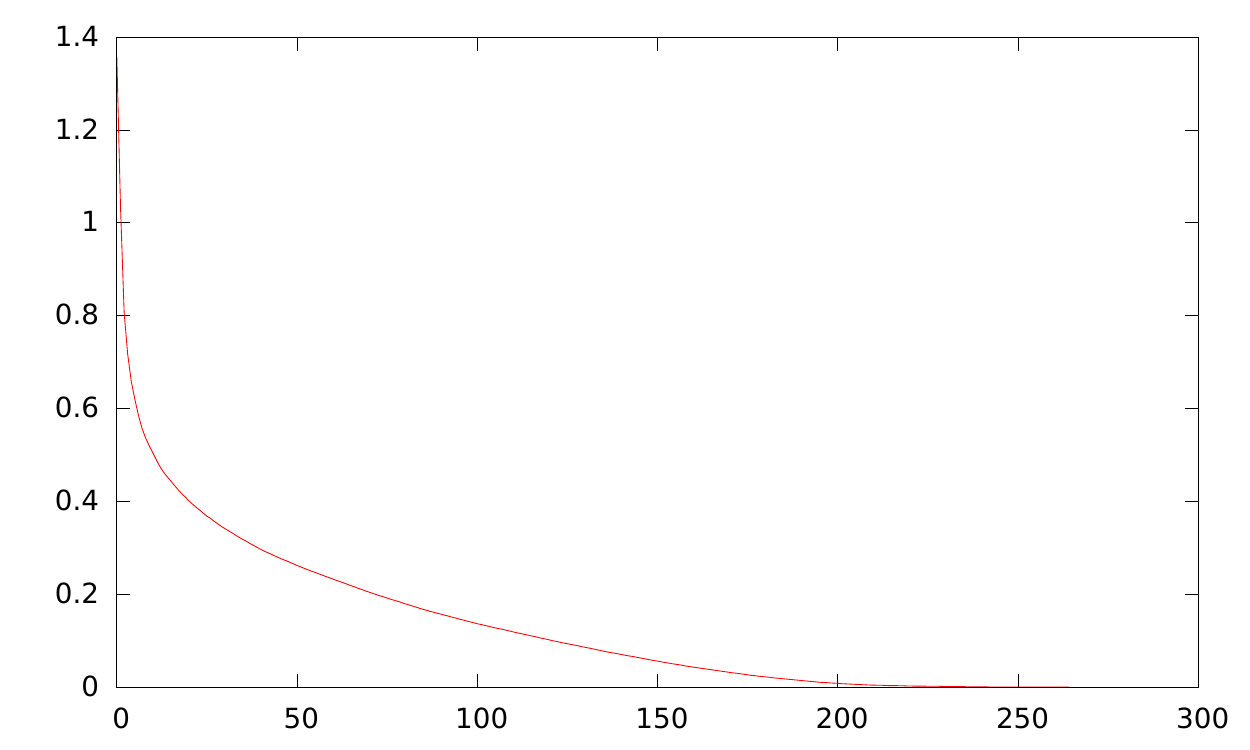}\\
Panel (c): Noise = 2 & Panel (d): Dominant interval length \\
\includegraphics[width=0.45\textwidth, height=30mm]{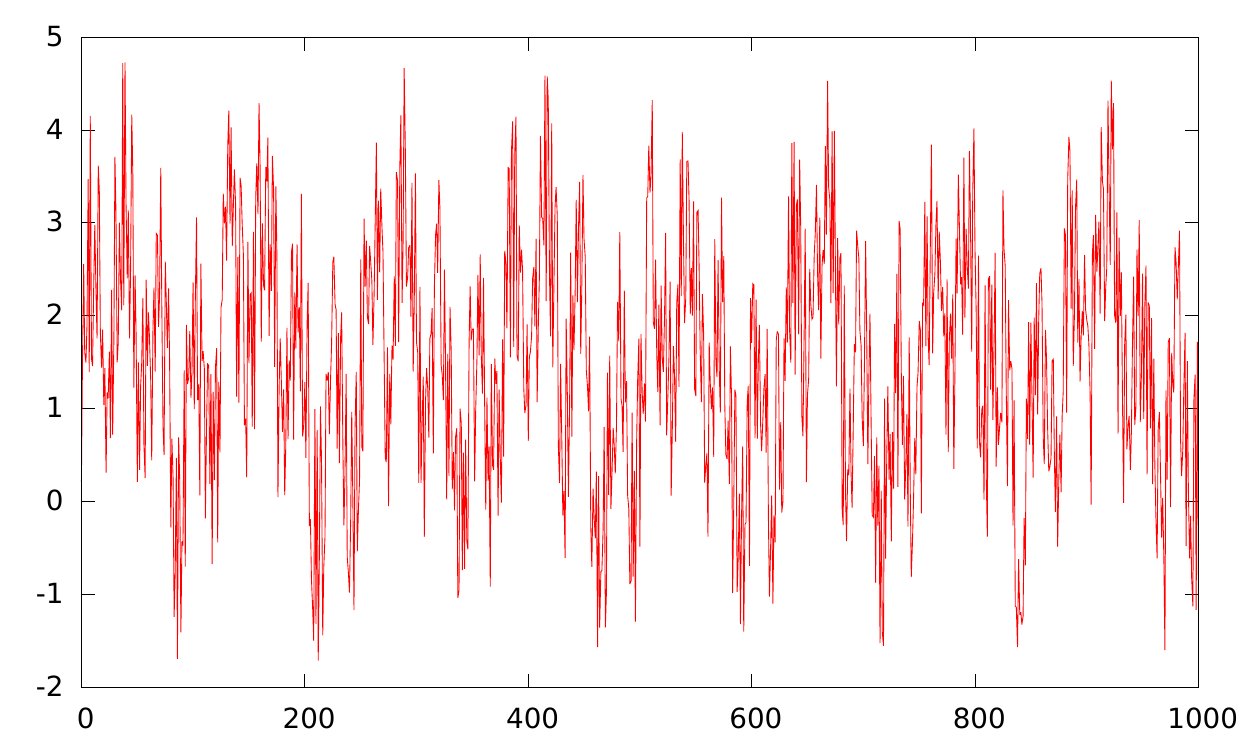}
&
\includegraphics[width=0.45\textwidth, height=30mm]{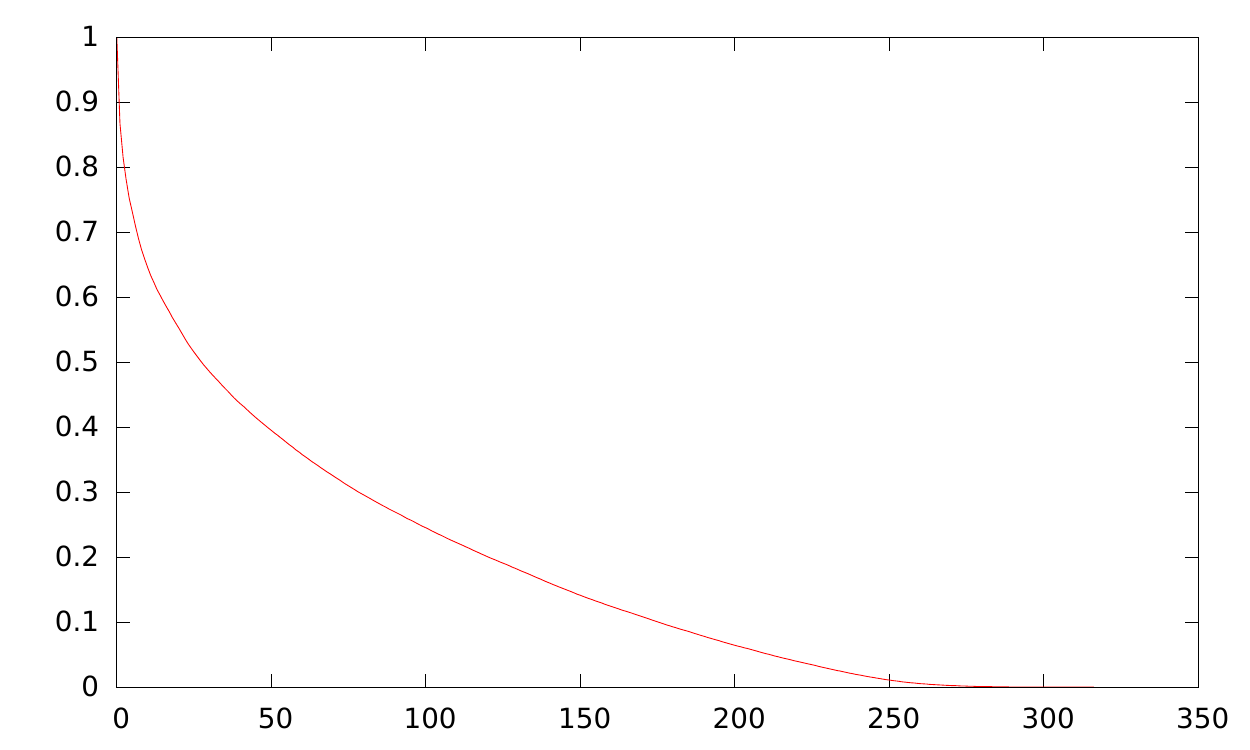}\\
Panel (e): Noise = 3 & Panel (f): Dominant interval length \\
\includegraphics[width=0.45\textwidth, height=30mm]{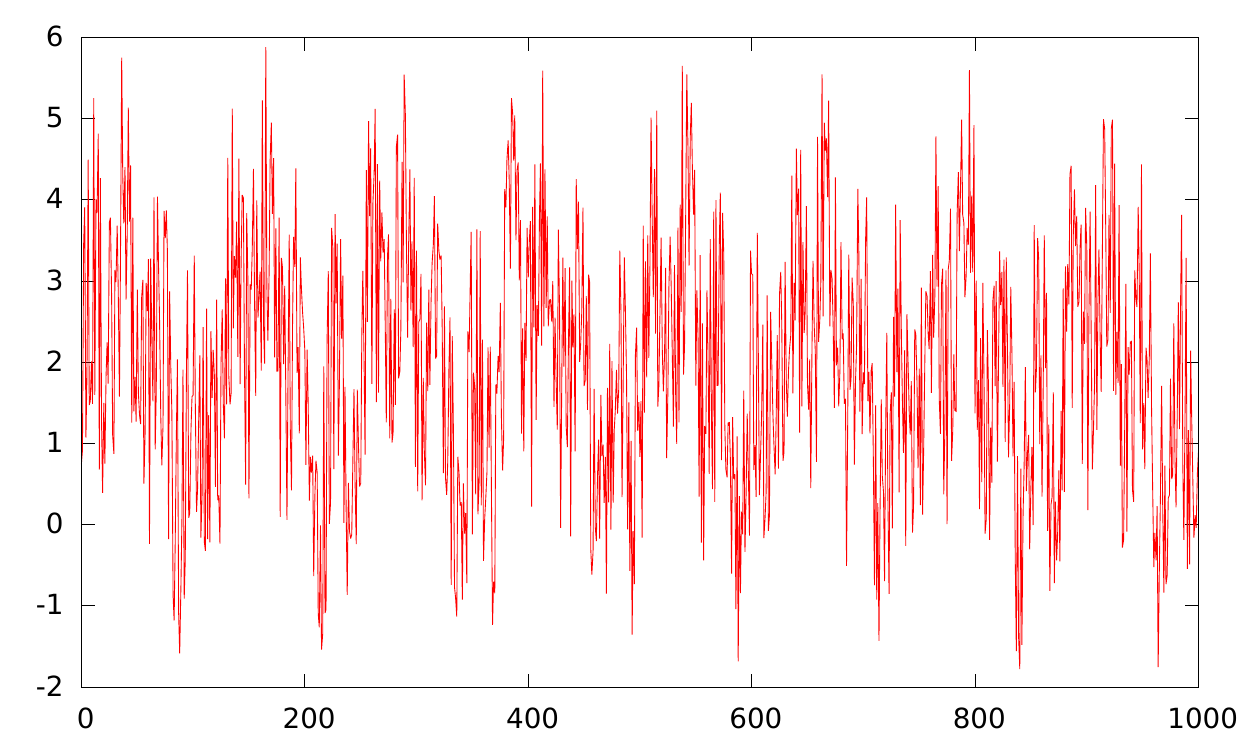}
&
\includegraphics[width=0.45\textwidth, height=30mm]{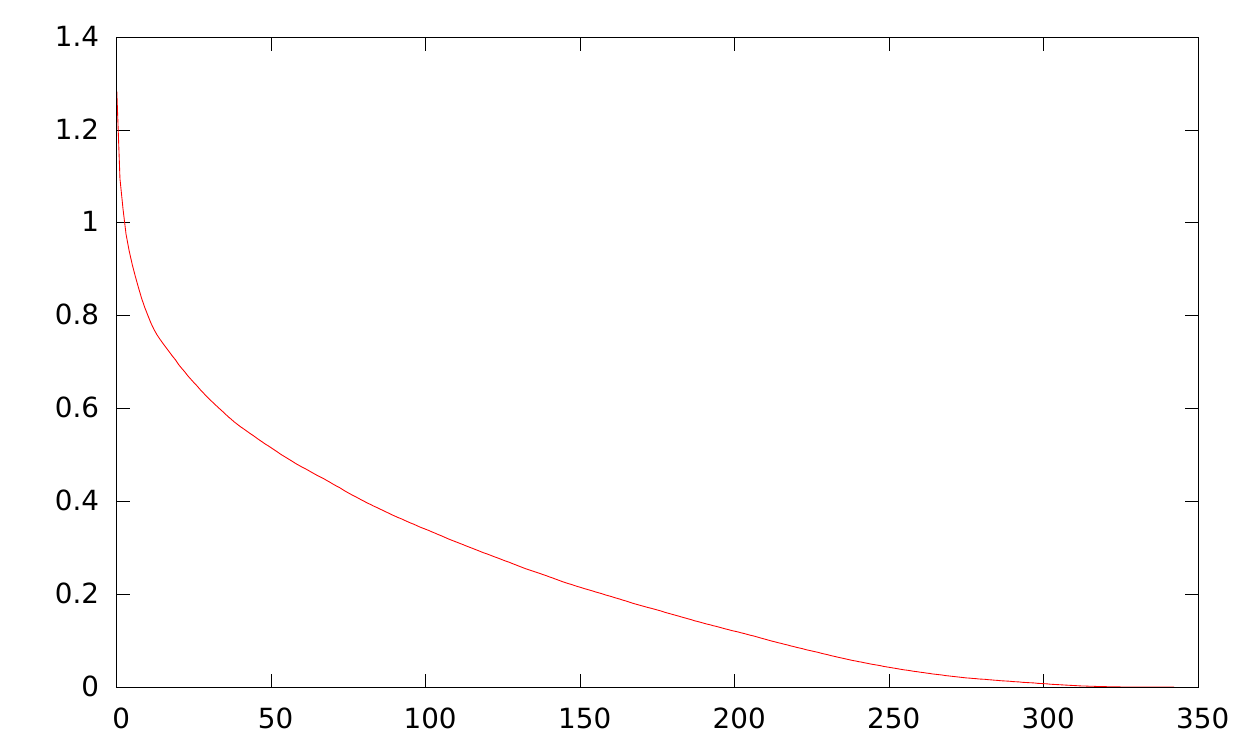}\\
Panel (g): Noise = 4 & Panel (h): Dominant interval length \\
\end{tabular}
\end{center}
\raggedright
\footnotesize{Notes: Panels (a), (c), (e) and (g) provide illustrations of sin functions with increasing levels of noise. In each case we plot the corresponding dominant interval length, death minus birth, for the persistence intervals derived from the application of persistent homology to the sliding window embedding of the time series. The plots of the right present the histograms of the lengths of persistence intervals (from the longest to the shortest) in dimension one. Note that for the small level of noise (top figure) there is clearly dominant interval. With the increasing level of noise (second to fourth rows) the shape of the histogram significantly changes.}
\end{figure}

In Figure \ref{fig:noisy_sin} we can observe a sinus summed with a uniform noise of an amplitude $1$, $2$, $3$ and $4$; panels (a), (c), (e) and (g) respectively. To the right of each we add the histogram of lengths of lengths of dominant persistence intervals in dimension $1$; panels (b), (d), (f) and (h)  respectively. The histograms were obtained as a average of $100$ random time series with the appropriate amount of noise. In the presence of noise, the fingerprint of a semi--periodic behavior is not the dominant persistence interval anymore, but rather a progressively decreasing sequence of persistence intervals in dimension $1$. This example demonstrates such clearly.

As a second example a composite sin function is presented, $sin(2x)+sin(x/2)$. This function has a part with a period  of just 31 and a second part with a period of 124. Hence, evaluation of the period requires an ability to identify both elements and presents a stronger challenge to the approach developed herein. Figure \ref{fig:noisy_sin3} demonstrates the function when combined with the same random noise draws as were used in the simple sinus case. As before we see the increased noise creates greater interval lengths - the line in panels (b), (d), (f) and (h) shift radially upwards in the way that was seen across the corresponding panels in Figure \ref{fig:noisy_sin}\footnote{Take for example interval length 100. In panel (b) of Figure \ref{fig:noisy_sin} there are almost no intervals of this length. In panel (d) that figure the density of intervals remains just below 0.2. By panel (f), when the noise reaches 3, the density is approximately 0.25. Finally in panel (h) the density is above 0.3.}. Because of the addition of a short period function there are naturally more short dominant intervals in this second example. Results presented are consistent with expectations and demonstrate the ability of the first stage of our approach to function for more complex time series.

\begin{figure}
\begin{center}
    \caption{$L^1$ and $L^2$ Norms of Noisy Sin Functions }
    \label{fig:noisy_sin_norms}
    \includegraphics[scale=0.8]{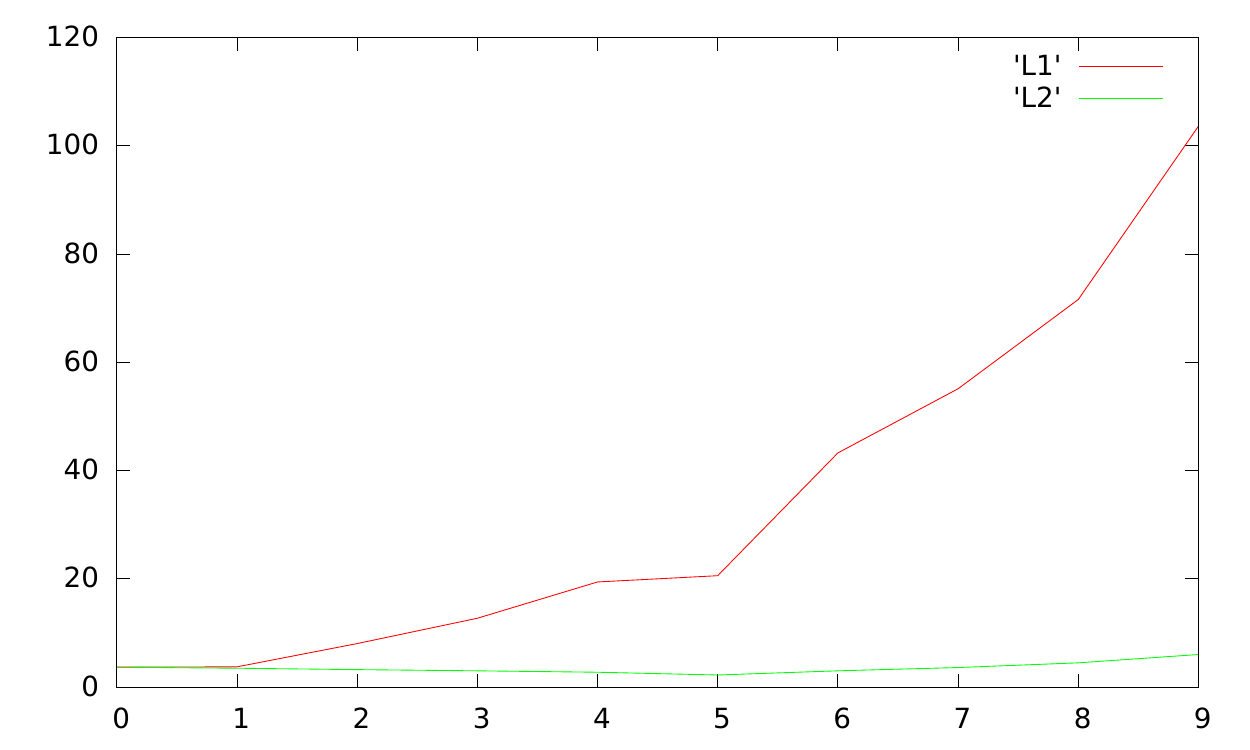}
\end{center}

  \footnotesize{Let us consider a sinusoidal time series with a uniform noise sampled from an interval $[0,i]$ for $i \in \{0,1,\ldots,9\}$. The $L^1$ and $L^2$ norms of its sliding window embedding are presented in the plot: Horizontal axis denotes the noise level and the vertical axis, the magnitude of the norms. Note that $L^1$ norms is significantly growing with $i$, as it accumulate the length of a large number of shot persistence intervals. The $L^2$ norm, which is the square root of sum of the squares of the lengths of the intervals is much less sensitive to this effect and therefore is a much more robust measure in the presence of large noise.}

\end{figure}

To present an aggregative summary of this information we will sum up the lengths (referred to as the $L^1$ norm), or the squares of the lengths (square root of which is the $L^2$ norm), of all persistence intervals in dimension $1$. The results are summarized in Figure~\ref{fig:noisy_sin_norms}. As one can observe the $L^1$ norm is quickly growing with the noise. Meanwhile the $L^2$ norm remains quite stable, and comparable to the length of the unique high persistence point (corresponding to the case with zero noise).

\section{''Chop and Search for Loop'' Procedure}
\label{sec:chop_and_loop}
\begin{figure}
\centering
\caption{Workflow for ``Chop and Search for Loop'' Procedure}
\includegraphics[scale=0.6]{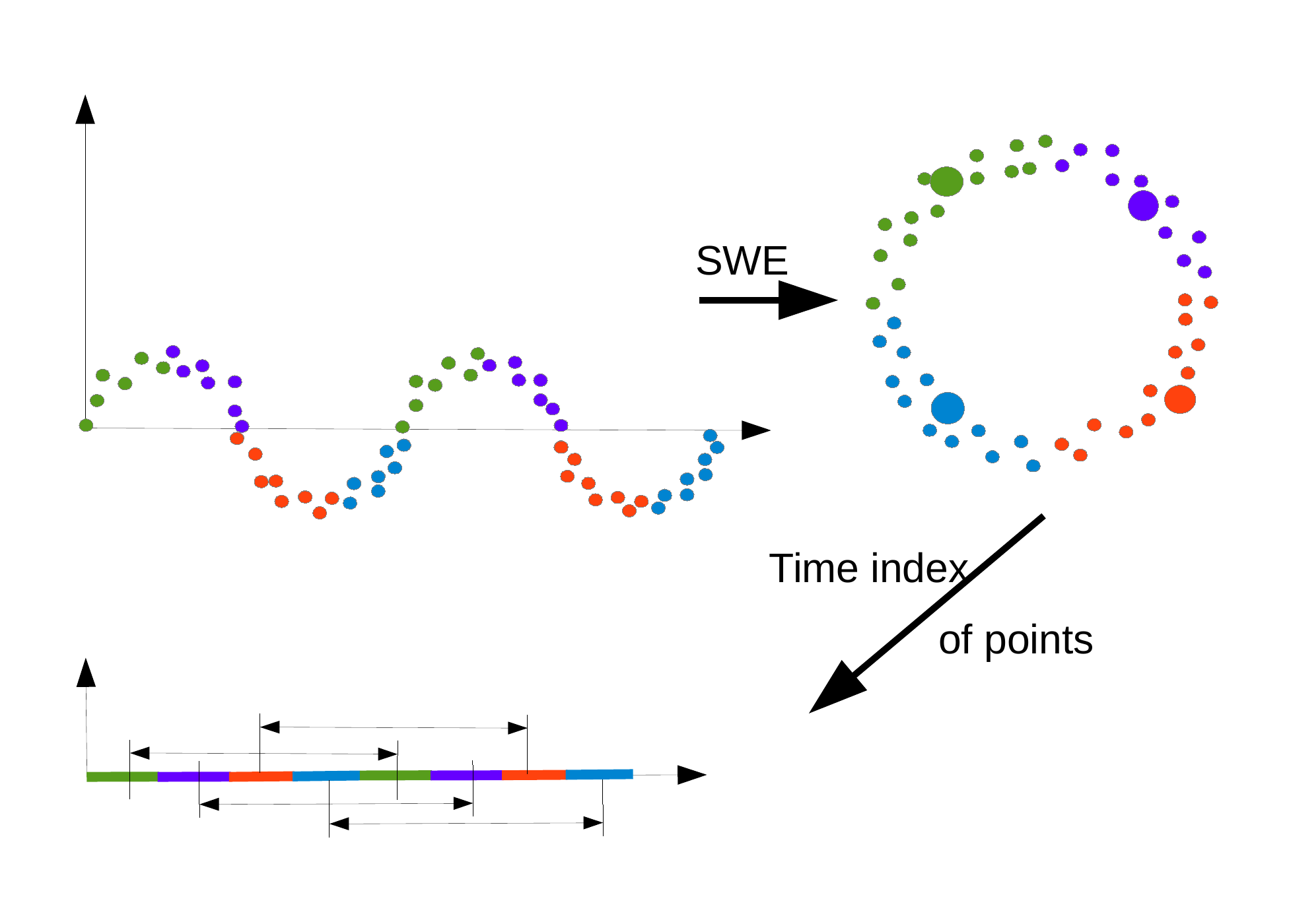}
\raggedright
\footnotesize{Notes:Starting from the time series on the left, we construct it sliding window embedding (uper right). Subsequently, a collection of four equispaced points (marked with larger disks) are selected and their Voronoi cells (colored) on the real line to the lower left of the figure, are located. Note that the colors from the Voroni cells are pulled back to the time series. Bottom, the time indices of the points coming from different Voronoi cells are retrieved. Note that in the bottom the spacing between points falling into the same Voronoi cells is roughly similar, hence we can conclude that the time series is periodic.}
\label{fig:equispaced_oscilations}
\end{figure}

Given the topological signatures of periodic time series and their noisy versions discussed in the Section~\ref{sec:math_approach} and~\ref{sec:noisy_time_series}, we will use the following three step mechanism to verify the semi--periodicity of a time series $f$:
\begin{enumerate}
\item Firstly the initial time series $f$ of a length $N$ gets converted into its SWE (using the parameters $n$ and $d$, see the Section~\ref{sec:math_approach} for details). Given a collection of high dimensional points we consider a sequence of subsets of points indexed between $[  max(0,i-M),i ]$ for $i \in \{ 1,\ldots,N \}$.
%will get chopped into $N-M$ sub-series, each of the length $M$. The first chopped time series, denoted as $f_M$, consist of all data points of $f$ of an index smaller than $M$. Second one, denoted as $f_{M+1}$, consist of all data points of $f$ in the range between $2$ and $M+1$, etc. The sub-index of $f$ will be used as a \emph{time indicator} of the sub series. Note the parameter $M$ in this step.
%
\item Compute the persistent homology of each subset and subsequently the $L^1$ and $L^2$ norm of it.
%For each of the sub series we will construct a sliding window embedding (using the parameters $n$ and $d$, see the Section~\ref{sec:math_approach} for details) and compute the persistent homology and subsequently the $L^1$ and $L^2$ norm of it.
%
\item For each subset of points for which the $L^1$ or $L^2$ norm are sufficiently large\footnote{The precise definition of this parameter is left to the user as its magnitude would be relative to the series being studied.}, we select a $k$ equispaced points $l_1,\ldots,l_k$ in the SWE of $f_i$. We will get this using either $k$-means~\cite{k_means} or max-min algorithm~\cite{de2004topological}. For every point $l_i$, we select all the points in SWE of $f_i$ that are closer to $l_i$ than to $l_j$, for $i \neq j$ and denote them by $V_{l_i}$ and call \emph{Voronoi cells} of $l_i$. Subsequently for every Voronoi cell $V_{l_i}$ we extract the time indices of points in it. If for every Voronoi cell they form an equispaced collection of time--indices, we get evidence of semi periodic oscillations. The idea of the procedure is presented in Figure~\ref{fig:equispaced_oscilations}.
\end{enumerate}

By looking at the ``time of return'' depicted in Figure~\ref{fig:equispaced_oscilations} one can observe this time is a proxy for a period of function. We would like to stress that when making observations on the time of return (as presented at the bottom of the picture) one should not only concentrate on one color (landmark point), but consider all of them and observe equispaced repetitions of the time indices of points.

Let us remark that making such observations by simply looking at the indices of points is possible, but time consuming and not very efficient unless some proper algorithm are to be used for the purpose. Suppose that the landmark points are marked with $0,1,\ldots,n$. For every landmark point $i$ we collect, in the increasing order, the time indices of points for which the landmark $i$ is the closest one into a vector $t_i$. Taking $t_i$ we will consider the difference between an element and the element directly before it, $t_i[l]-t_i[l-1]$ for suitable values of $l$. That will allow us to isolate all the ``jumps'' in the indices. For periodic time series with a constant, $i$ invariant, period $p$, we expect to see a sequence of points centered abound $t_{0}, p+t_{0}, 2p+t_{0}$ etc. The constant $p$ is a proxy for period and can be recovered from the dominant jumps in $t_i$. Note that this jump will approximate the time between \emph{leaving}
 the given landmark to \emph{entering it back again}. Consequently in order to retrieve the period from it, we still need to add the number of steps on which the time series is \emph{staying} in the given landmark. Equivalently for any landmark point we estimate the distance between two constitutive dominant jumps. The described procedure is formalized in Algorithm~\ref{alg:perod_detection}.

\begin{algorithm}[H]
\SetAlgoLined
\caption{Period Estimation}
\KwResult{Period extraction}
 \For{every landmark point i}{
  $t_i$ = sorted time indices of points for which $i$ is the closest landmark\;
  $j_i$ = vector gathering $t_i[k]-t_i[k-1]$ for $k \in \{1,\ldots,len(t_i)\}$\;
  %Estimate the constitutive jumps in the values. The distance between their positions is a proxy for the period.
  Sort $j_i$ from largest to shortest.
  %The longest elements of $j_i$ gives an estimation for a \emph{time of return} to a given Voronoi cell. In the unsorted version of $j_i$ they will be present as a \emph{picks}. An averaged distance between two constitutive picks is an estimation for a period of the time series.
  Locate the dominant picks in $j_i$ and verify if they are approximately equally spread. If so, their average spread is the proxy for period at the landmark $i$ Collect the period's proxies from each landmark point, report their average and standard deviation.
 }
 \label{alg:perod_detection}

\end{algorithm}

\begin{figure}[h!]
\caption{Estimation of a Period for sin with Different Levels of Noise}
\label{fig:j_i_histograms_single_sin}
\begin{tabular}{ c  c  }
\includegraphics[width=0.45\textwidth, height=60mm]{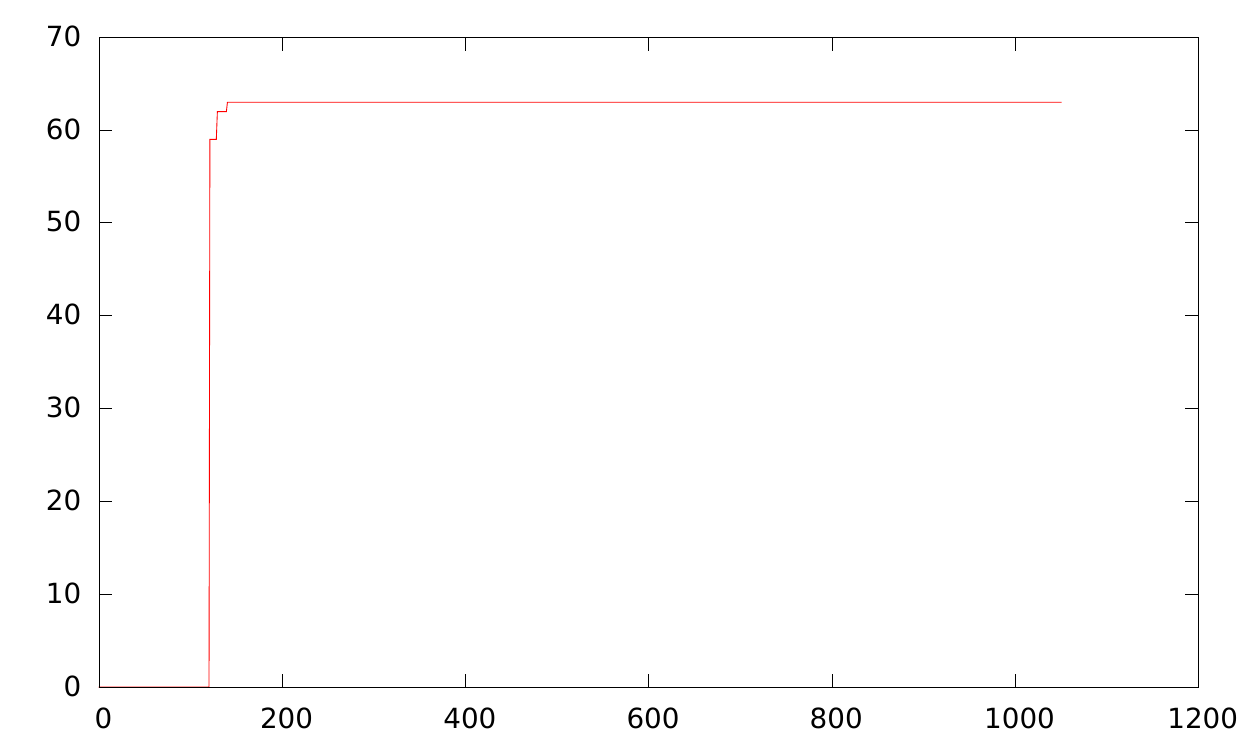}&
\includegraphics[width=0.45\textwidth, height=60mm]{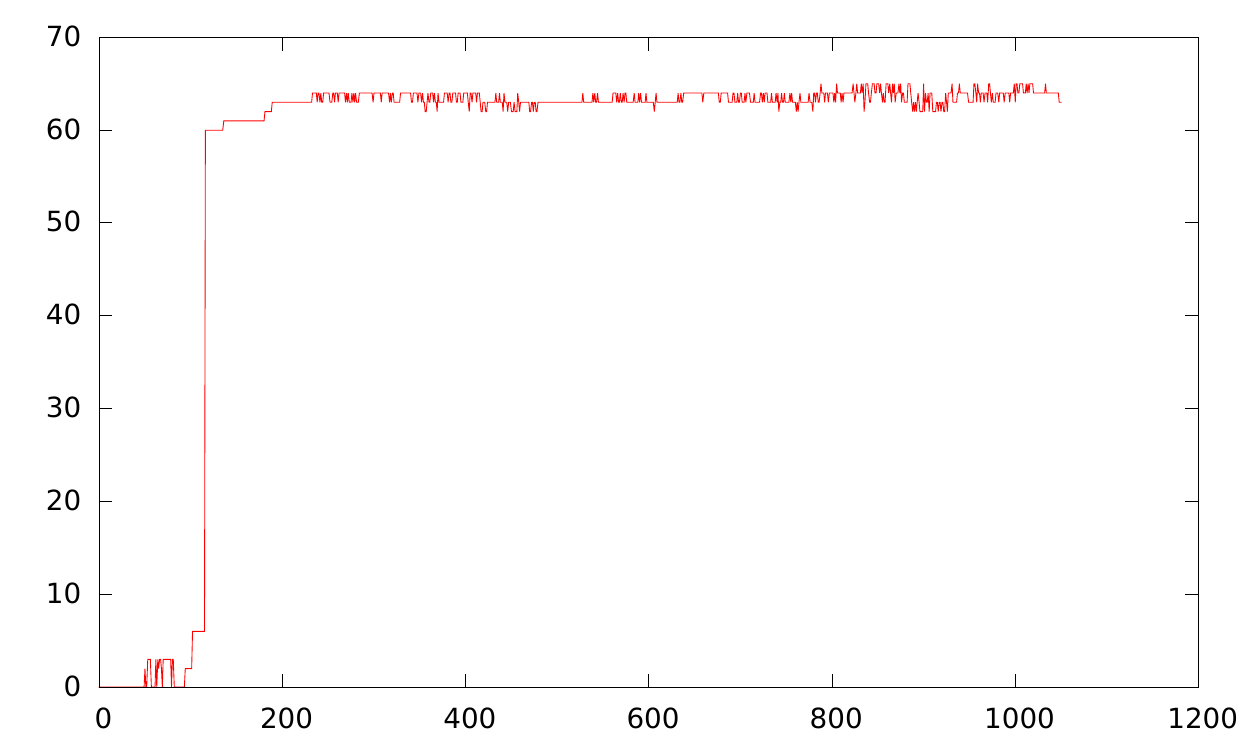}
\\
Panel (a): Noise 1 & Panel (b): Noise 2 \\
\includegraphics[width=0.45\textwidth, height=60mm]{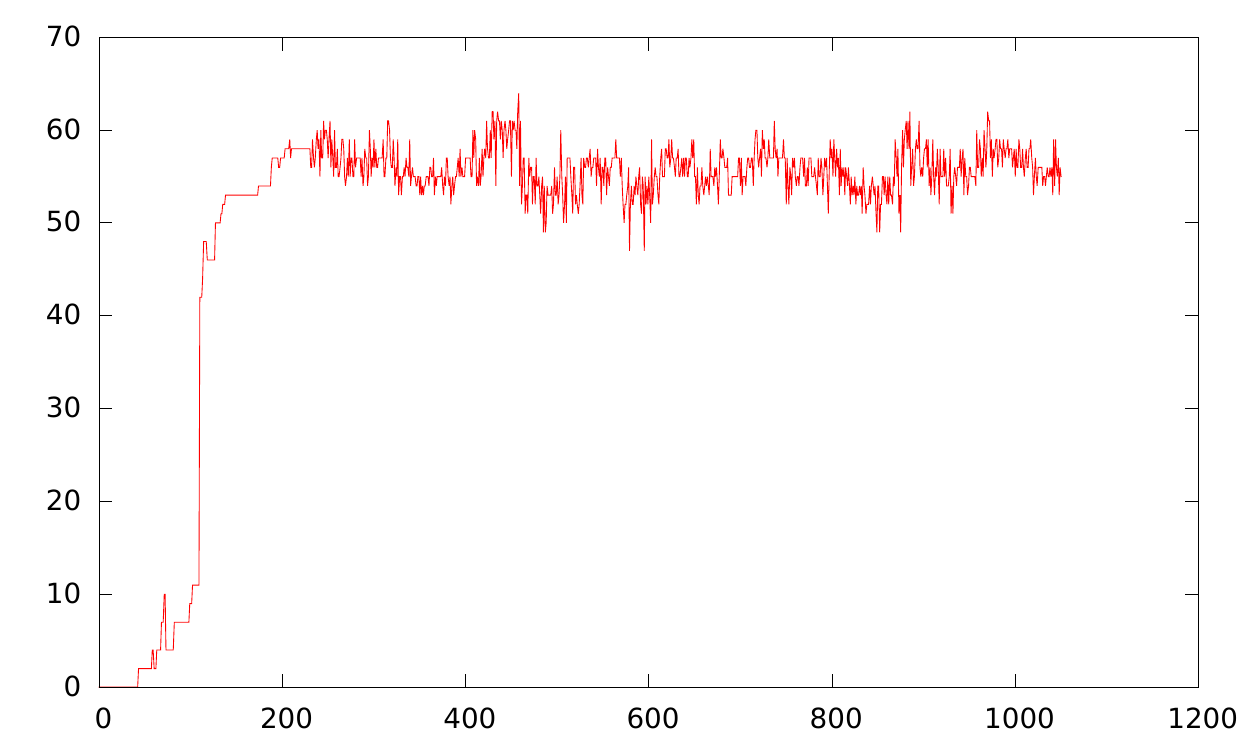}&
\includegraphics[width=0.45\textwidth, height=60mm]{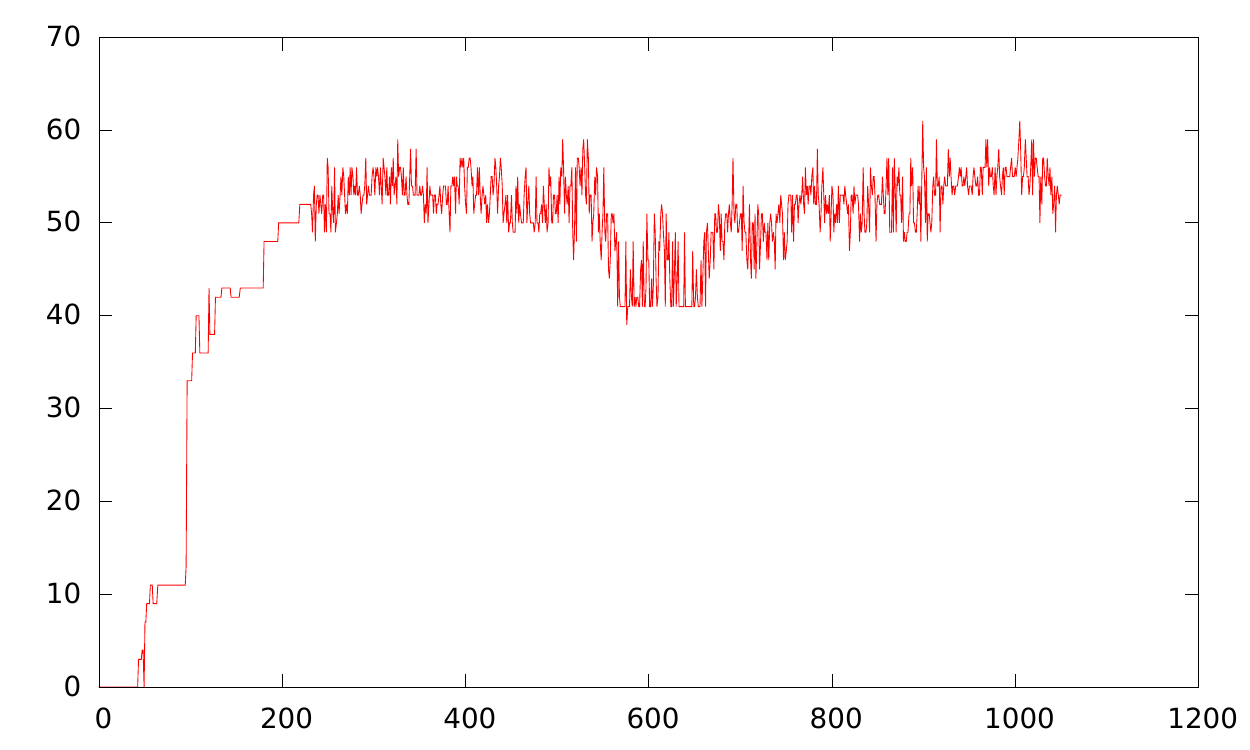}\\
Panel (c): Noise 3 & Panel (d): Noise 4\\
\end{tabular}
\raggedright
\footnotesize{Notes: Horizontal axis denotes observation number. True underlying period is $2\times10\times\pi \approx 63$}
\end{figure}

Let us test the presented methodology of the noisy version of sinus time series presented in the Figure~\ref{fig:noisy_sin}. The time series was obtained by taking a grid of $sin(x_i)$ for $x_i$ between $0$ and $100$ with $0.1$ stepsize. The plots of the $j_i$'s are provided on the Figure~\ref{fig:j_i_histograms_single_sin}.

As one can observe all the plots in the first two panels  of Figure~\ref{fig:j_i_histograms_single_sin} stabilize around $62$. When the noise increases we do notice a reduction in the estimated period but it is still evident that there is periodic behaviour for these series. The initial growth at lower observation numbers is a consequence of an insufficient number of points being available in the initial point clouds to close up the cycle. As soon as there are sufficiently many points to notice periodic behavior, it is captured in the plots. The stabilization at the level of roughly $60$ steps is a correct answer as the function repeats itself every $63$ steps. We can see only one step in the plots as the size of the sub-domains was not large enough to capture more than two periods of the function. When extending the step size and/or the size $N$ of the subseries, we would have obtain another step at the plots at the  level $125$.

\begin{figure}
\caption{Composite sin Functions and Periodicity Detection}
\begin{tabular}{ c  c  }
\includegraphics[width=0.4\textwidth, height=50mm]{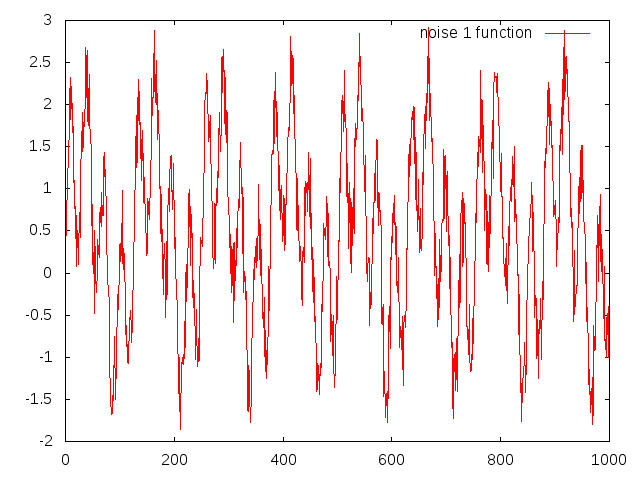}&
\includegraphics[width=0.4\textwidth, height=50mm]{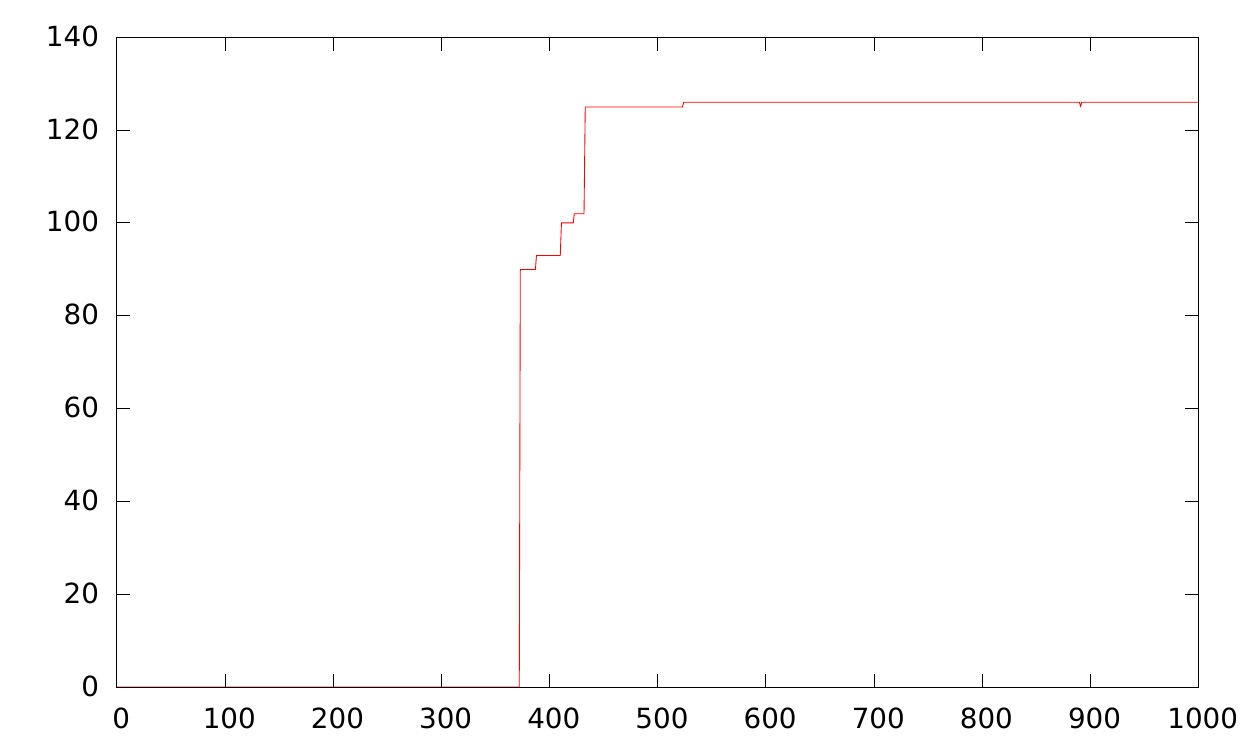}\\
$sin(2x)+sin(\frac{x}{2})$ with uniform noise from$[0,1]$ & Period estimation plot, no noise \\
\\
\includegraphics[width=0.4\textwidth, height=60mm]{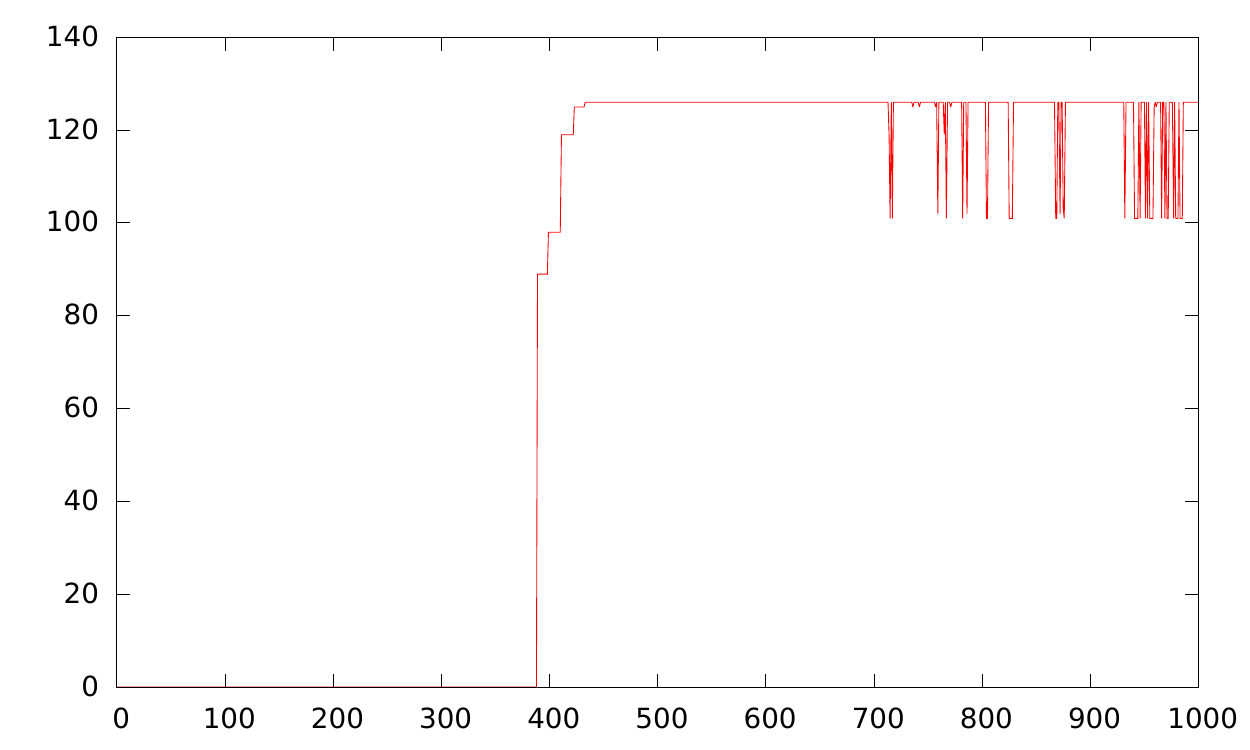}&
\includegraphics[width=0.4\textwidth, height=50mm]{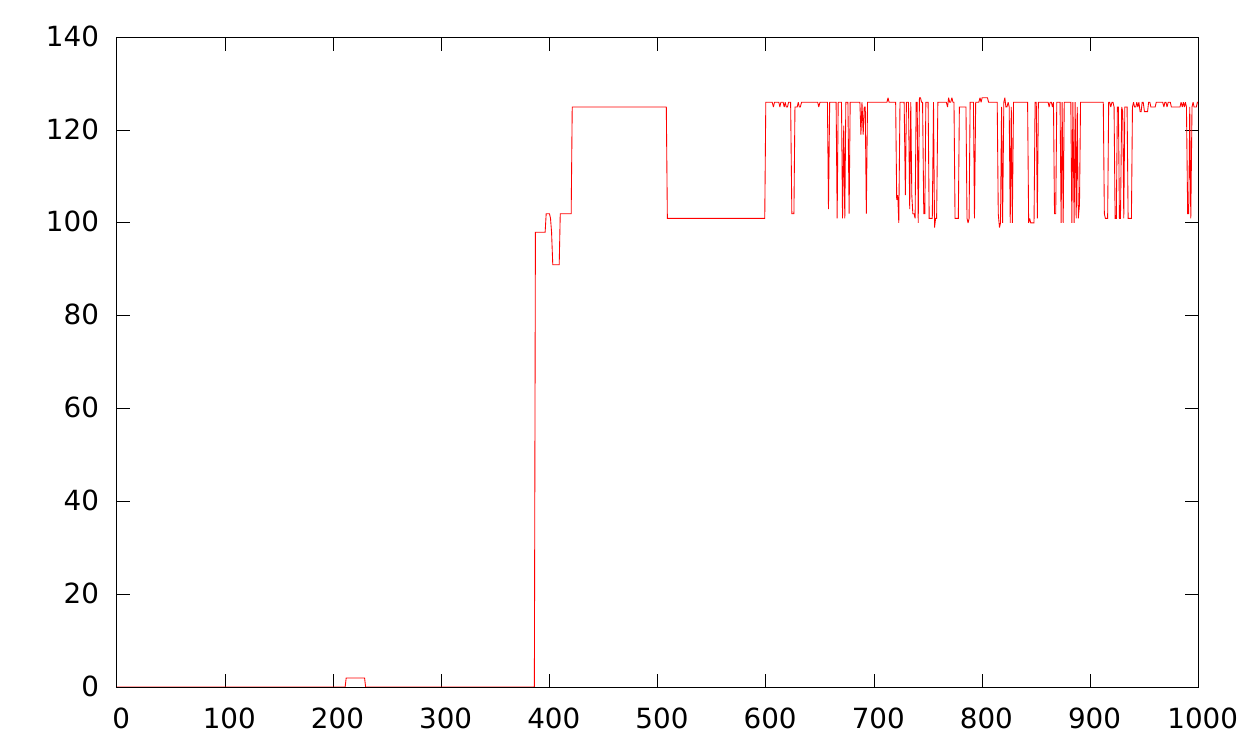}\\
Period estimation plot, $[0,1]$ noise & Period estimation plot $[0,2]$ noise \\
\includegraphics[width=0.4\textwidth, height=50mm]{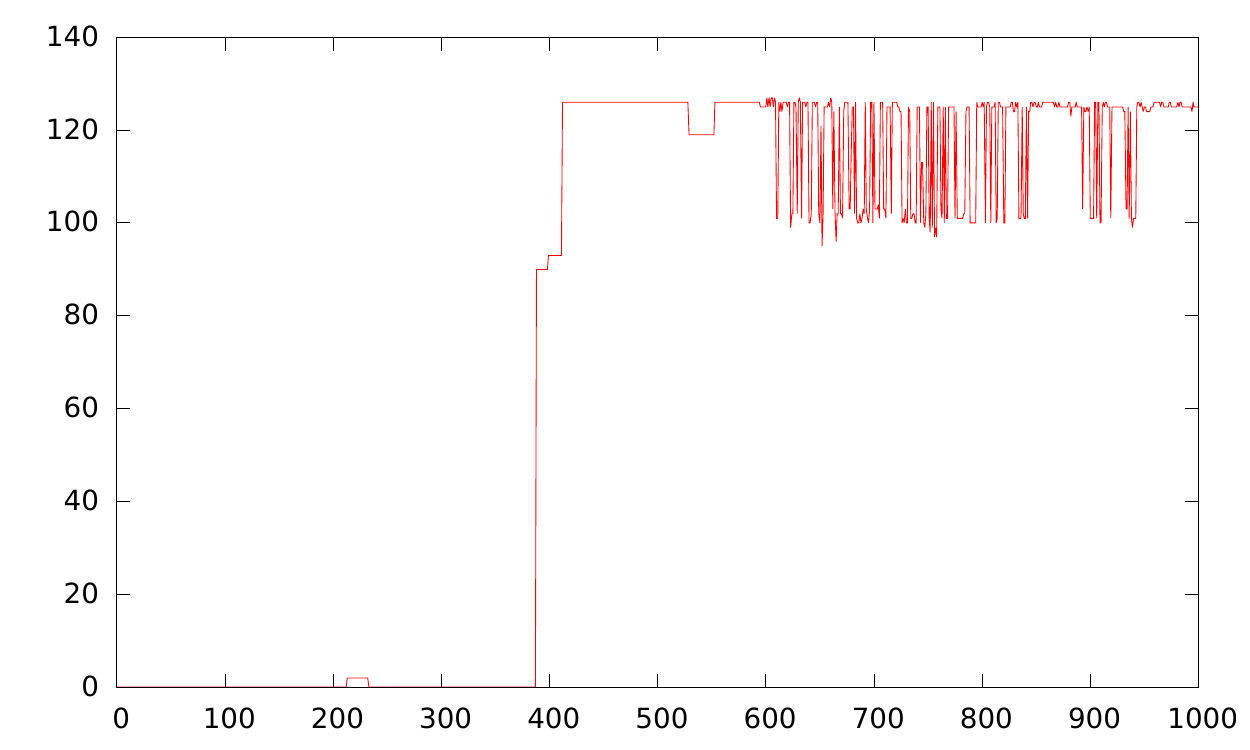}&
\includegraphics[width=0.4\textwidth, height=50mm]{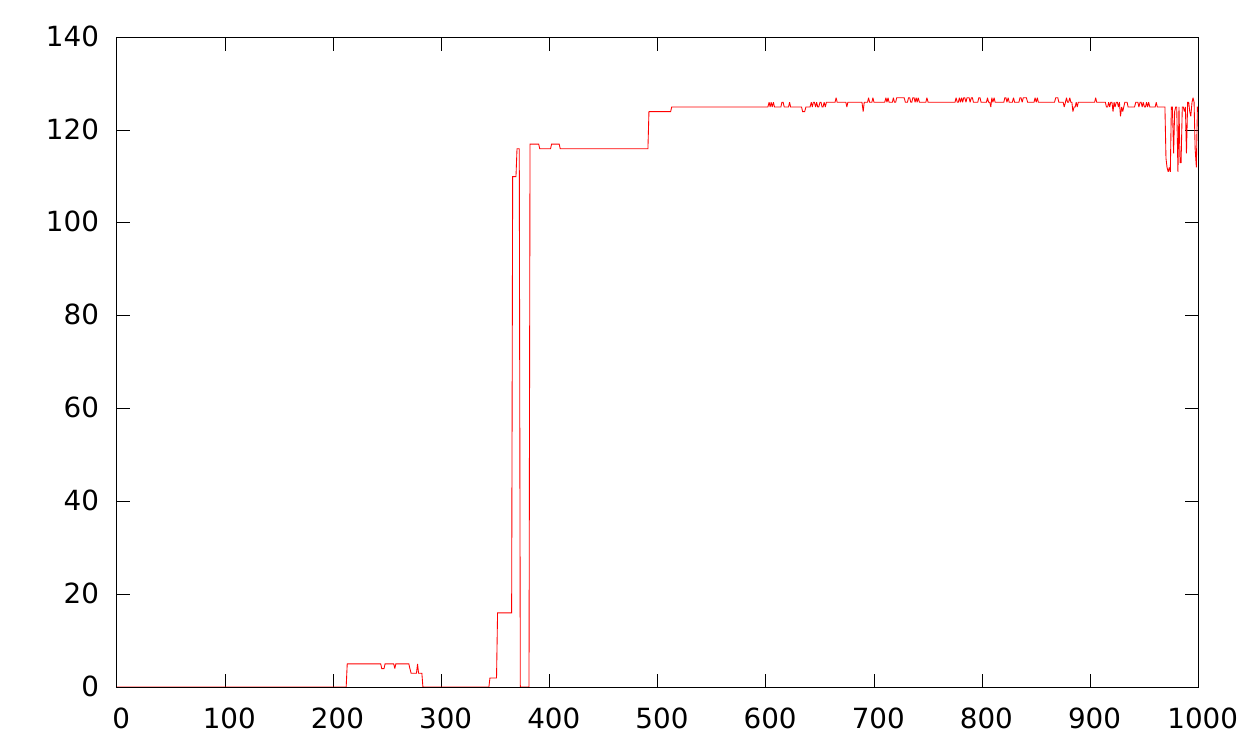}\\
Period estimation plot, $[0,3]$ noise & Period estimation plot, $[0,4]$ noise \\
\end{tabular}
\raggedright
\footnotesize{Notes: We consider the time series of $sin(2x) + sin( \frac{x}{2} )$ with the uniform noise sampled from
$[0,i]$, for $i \in \{ 0,1,2,3,4 \}$. Note that the period of the function is $10\ 4\ \pi \approx 125$ The
top right picture present the function with the nose sampled uniformly from $[0,1]$. The remaining
pictures present period estimation plot. In each of the case the plot stabilizes at the level $120$ which gives
a good estimate for the period.}
\label{fig:j_i_histograms}
\end{figure}

The presented techniques are again tested upon a composition of two periodic functions. We have chosen $sin(2\ x)+sin(0.5\ x)$ sampled on the same grid as described above. The function with uniform noise sampled from $[0,1]$ is presented in the Figure~\ref{fig:j_i_histograms}. On the same picture we consider plots of period estimation of the function with no noise, a random noise sampled from $[0,1]$, $[0,2]$, $[0,3]$ and $[0,4]$ as explained in the previous cases.

We have demonstrated with artificial examples how TDA can recover periodicity quickly from artificial examples and how even in the presence of strong noise there is still identification of the underlying cyclical process. We next consider how a commonly employed alternative compares to this work as a mechanism to reaffirm the value of the TDA approach exposited over the past two sections.

%The input for an automatic analysis would be a collection of time indicators of points falling into each Voronoi cell. Given this, in this paper wpropose two standard options detailed below. Note that each of them should be applied to each Voronoi cell separately and the results of obtained analysis should be compared in order 5 draw a conclusion.
%\begin{enumerate}
%\item By using dendrograms, as presented in the Fig ToDo. As one can see, the time indices that are close by cluster together quickly. At the next level of the dendrogram the clusters that are one period away cluster together. The distance between those clusters, of more or less constant, should be considered the proxy for the period of the time series.
%
%\item By using zero dimensional persistent homology in dimension $0$. In this case the input for PH analysis is a collection of time indices (hence points in $\mathbb{R}$). As one can observe at the bottom of Fig~\ref{fig:equispaced_oscilations}, we would expect to have a number of short persistence points (corresponding to time indices in each connected part of a color) merging together. Later, there should be a number of high persistence points. The persistence of each of them will be a proxy for the period. They will correspond to the connected components of colors (and therefore images of points after each period) merging together. Please see Fig ToDo for further reference.
%\end{enumerate}

\section{Comparison}
\label{sec:compare}

Within the time series analysis community a well understood approach to the estimation of periodicity is the \emph{findfrequency} function developed in R package \emph{forecast} \cite{hyndman2019package} . To evaluate the performance of our proposed TDA approach we return to the artificial functions explored in the previous section and apply the \cite{hyndman2019package} code. Spectral methods like these do not have a natural parameter like the sliding window we used in our TDA evaluation meaning that to gain fair comparison it is instructive to use a rolling window approach. In all that follows we use rolling windows of various sizes, with the variation between the sizes dependent on the length of the series being considered.

Table \ref{tab:sinhynd} provides summary statistics for the sin function with noise 1,2,3 and 4. The overall length of the series is 400 and so we use intervals from a minimum of 40 to a maximum of 250 in increments of 10. For the shortest interval we thus have 361 frequency observations over which to average\footnote{The first estimate from an interval of length 40 uses observations 1 to 40 and provides an estimate of the period at time 40. Consequentially we have estimates from 40 to 400 inclusive, or 361 estimates.}. To summarise across the estimates we report the mean, the standard deviation and the minimum and maximum estimates. In every case we see minimum of 0. Whilst this might be expected on an interval shorter than the period, we would expect to see estimates appearing after more observations. When the noise is low the \textit{forecast} package \cite{hyndman2019package} code is able to recover a period around the true value as we approach the third cycle, 150 observations on a true period of 62. However as the noise increases this closeness quickly goes. An initial conclusion from Table \ref{tab:sinhynd} is that TDA is significantly outperforming the \cite{hyndman2019package} approach on the artificial examples tested so far.

\begin{sidewaystable}
    \caption{Period Estimation from Rolling-Windows and \cite{hyndman2019package}}
    \label{tab:sinhynd}
	\begin{tabular}{l c c c c c c c c c c c c c c c c c }
		\hline
		&&\multicolumn{4}{l}{Sin with Noise 1}&\multicolumn{4}{l}{Sin with Noise 2}&\multicolumn{4}{l}{Sin with Noise 3}&\multicolumn{4}{l}{Sin with Noise 4}\\
		Length & Obs & Mean & S.d. & Min & Max & Mean & S.d. & Min & Max & Mean & S.d. & Min & Max & Mean & S.d. & Min & Max\\
		\hline
		40    & 360   & 0.975 & 0.203 & 0     & 3     & 1.047 & 0.429 & 0     & 3     & 1.091 & 0.699 & 0     & 5     & 1.091 & 0.699 & 0     & 5 \\
    50    & 350   & 1.185 & 0.695 & 0     & 4     & 1.316 & 0.891 & 0     & 5     & 1.490 & 2.516 & 0     & 43    & 1.490 & 2.516 & 0     & 43 \\
    60    & 340   & 2.006 & 4.985 & 0     & 55    & 1.628 & 1.132 & 0     & 6     & 3.173 & 7.896 & 0     & 67    & 3.173 & 7.896 & 0     & 67 \\
    70    & 330   & 12.157 & 34.933 & 0     & 333   & 2.414 & 4.931 & 0     & 48    & 5.483 & 13.466 & 0     & 91    & 5.483 & 13.466 & 0     & 91 \\
    80    & 320   & 21.358 & 66.434 & 0     & 998   & 2.555 & 5.067 & 0     & 67    & 5.006 & 12.685 & 0     & 125   & 5.006 & 12.685 & 0     & 125 \\
    90    & 310   & 24.550 & 33.737 & 0     & 333   & 3.849 & 9.614 & 0     & 59    & 6.473 & 30.554 & 0     & 499   & 6.473 & 30.554 & 0     & 499 \\
    100   & 300   & 29.711 & 30.957 & 0     & 91    & 5.176 & 12.341 & 0     & 62    & 4.369 & 8.214 & 0     & 143   & 4.369 & 8.214 & 0     & 143 \\
    110   & 290   & 24.667 & 29.423 & 0     & 100   & 10.835 & 20.273 & 0     & 62    & 6.278 & 24.067 & 0     & 333   & 6.278 & 24.067 & 0     & 333 \\
    120   & 280   & 24.562 & 28.484 & 0     & 77    & 12.171 & 20.833 & 0     & 59    & 4.676 & 3.293 & 0     & 55    & 4.676 & 3.293 & 0     & 55 \\
    130   & 270   & 38.657 & 30.287 & 0     & 83    & 15.100 & 23.262 & 0     & 62    & 5.694 & 11.164 & 0     & 143   & 5.694 & 11.164 & 0     & 143 \\
    140   & 260   & 50.015 & 24.664 & 0     & 83    & 15.866 & 23.652 & 0     & 62    & 10.134 & 25.171 & 0     & 250   & 10.134 & 25.171 & 0     & 250 \\
    150   & 250   & 57.586 & 15.915 & 0     & 77    & 15.355 & 23.968 & 0     & 62    & 13.016 & 30.332 & 0     & 200   & 13.016 & 30.332 & 0     & 200 \\
    160   & 240   & 55.834 & 18.192 & 0     & 71    & 15.564 & 24.105 & 0     & 62    & 6.925 & 17.230 & 0     & 200   & 6.925 & 17.230 & 0     & 200 \\
    170   & 230   & 52.186 & 21.216 & 0     & 77    & 14.160 & 23.318 & 0     & 62    & 15.978 & 31.959 & 0     & 166   & 15.978 & 31.959 & 0     & 166 \\
    180   & 220   & 55.946 & 15.178 & 0     & 67    & 12.321 & 21.957 & 0     & 62    & 18.380 & 47.832 & 0     & 499   & 18.380 & 47.832 & 0     & 499 \\
    190   & 210   & 55.028 & 16.811 & 0     & 67    & 17.180 & 25.186 & 0     & 62    & 19.919 & 47.404 & 0     & 499   & 19.919 & 47.404 & 0     & 499 \\
    200   & 200   & 55.463 & 15.841 & 0     & 62    & 21.975 & 27.209 & 0     & 62    & 21.721 & 27.773 & 0     & 111   & 21.721 & 27.773 & 0     & 111 \\
    210   & 190   & 56.822 & 14.163 & 0     & 62    & 27.267 & 28.637 & 0     & 67    & 23.068 & 25.867 & 0     & 62    & 23.068 & 25.867 & 0     & 62 \\
    220   & 180   & 57.320 & 14.685 & 0     & 62    & 36.702 & 28.274 & 0     & 62    & 29.099 & 27.201 & 0     & 62    & 29.099 & 27.201 & 0     & 62 \\
    230   & 170   & 56.801 & 15.023 & 0     & 62    & 41.842 & 27.177 & 0     & 62    & 30.006 & 27.436 & 0     & 62    & 30.006 & 27.436 & 0     & 62 \\
    240   & 160   & 56.298 & 15.373 & 0     & 62    & 43.273 & 26.500 & 0     & 62    & 30.863 & 27.171 & 0     & 62    & 30.863 & 27.171 & 0     & 62 \\
    250   & 150   & 56.152 & 15.858 & 0     & 62    & 45.146 & 25.907 & 0     & 62    & 38.093 & 26.324 & 0     & 62    & 38.093 & 26.324 & 0     & 62 \\

		\hline

	\end{tabular}
\raggedright
\footnotesize{Notes: Findfrequency function from R package \textit{forecast} \cite{hyndman2019package} is used on a rolling window basis with the size of the window given in column 1. This then leaves the remaining observations over which the period is calculated. Statistics report the average, standard deviation, minimum and maximum, of the periods that are calculated for the remaining observations.}
\end{sidewaystable}

\begin{figure}
    \begin{center}
        \caption{Frequency Estimates on Noisy sin Functions Using \cite{hyndman2019package}}
        \label{fig:sinhynd}
        \begin{tabular}{c c}
             \includegraphics[width=5cm]{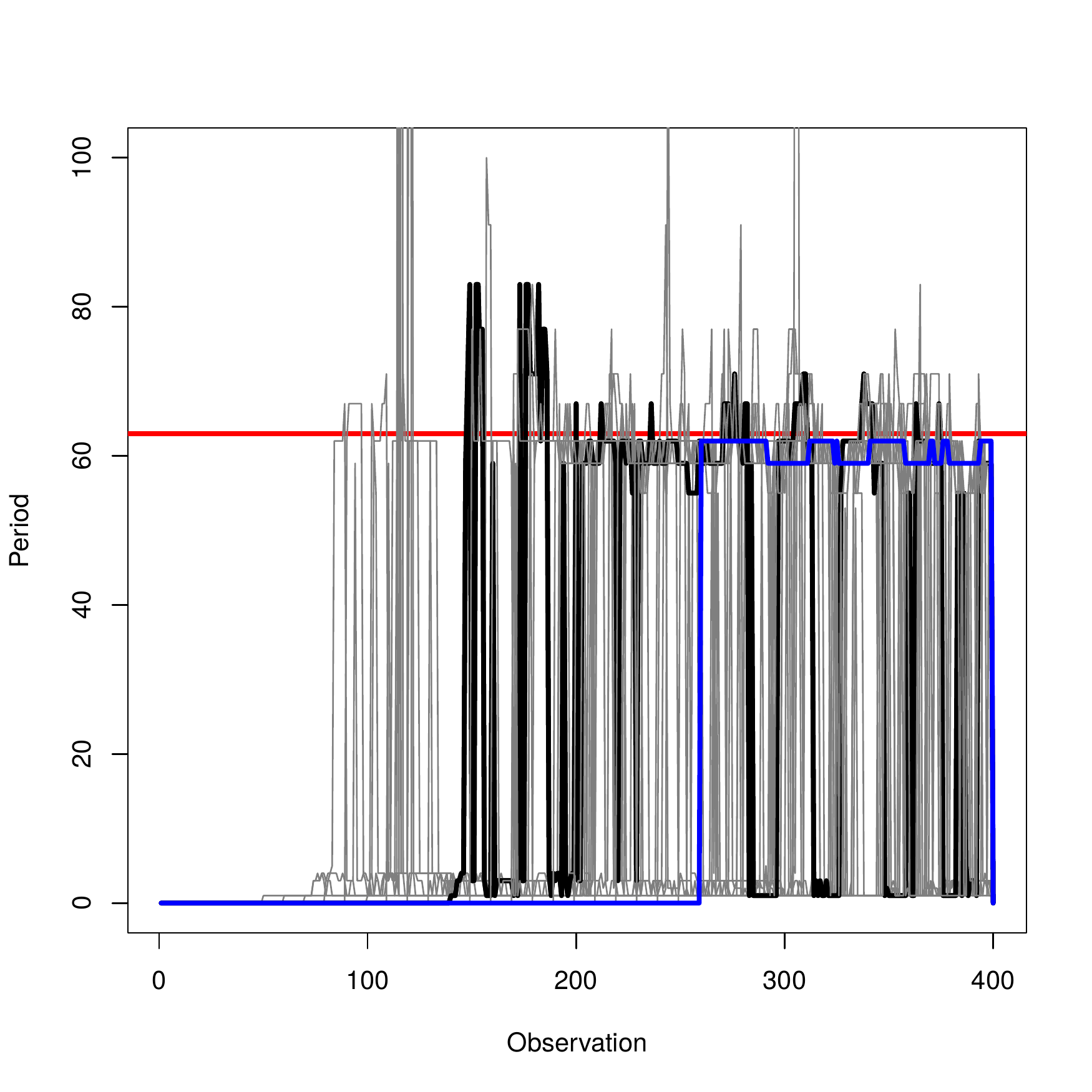}&
             \includegraphics[width=5cm]{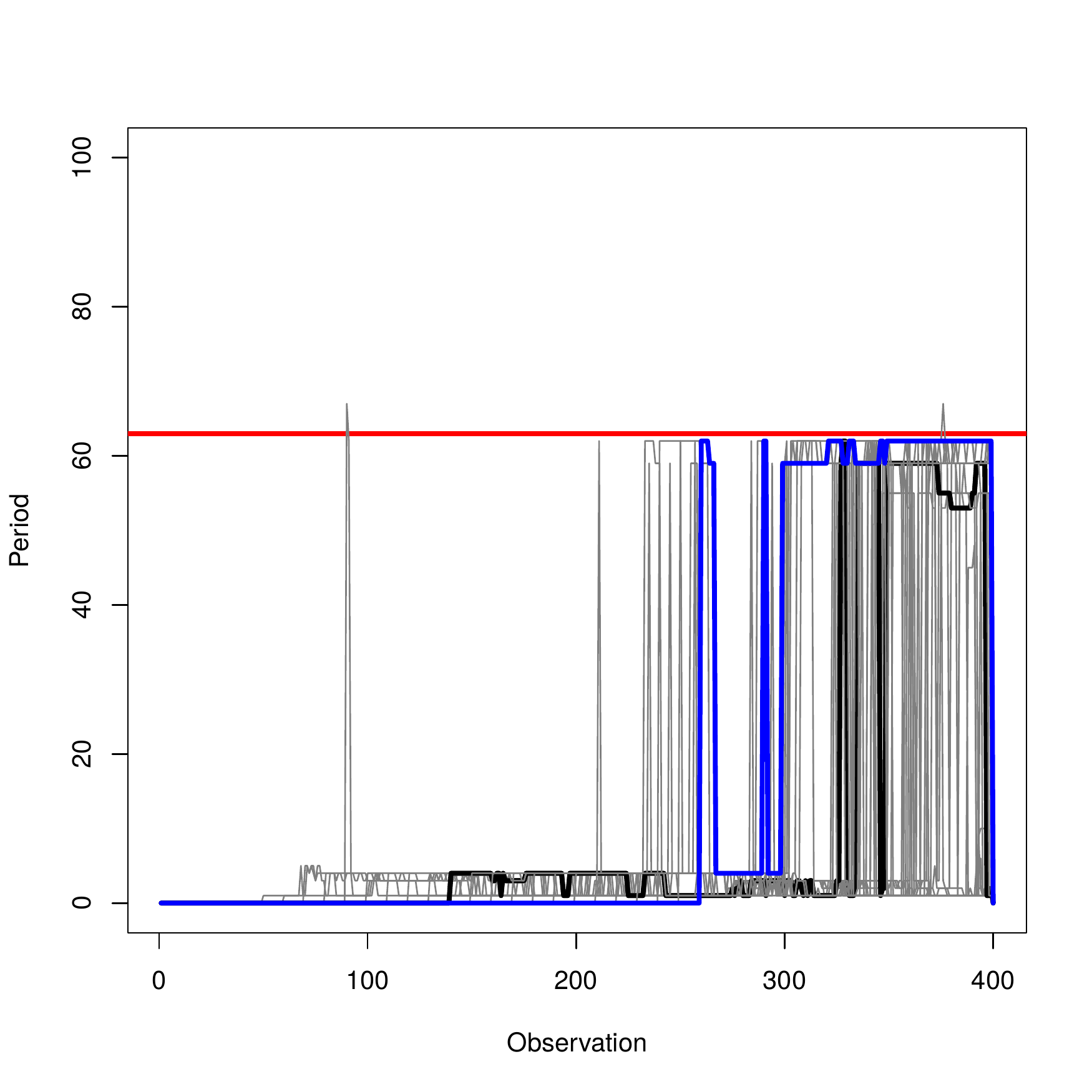}\\
             Panel (a): Noise 1 & Panel (b): Noise 2 \\
             \includegraphics[width=5cm]{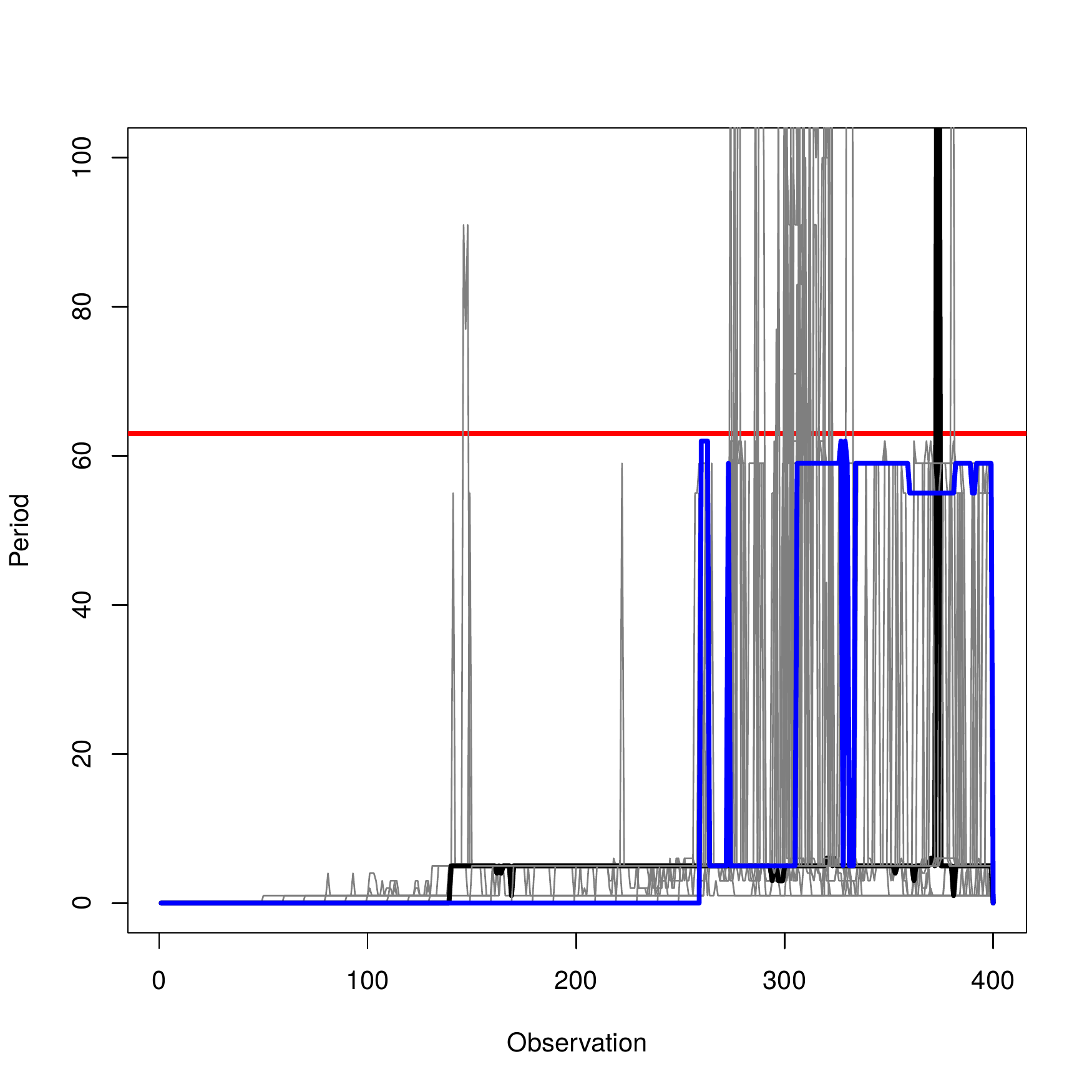}&
             \includegraphics[width=5cm]{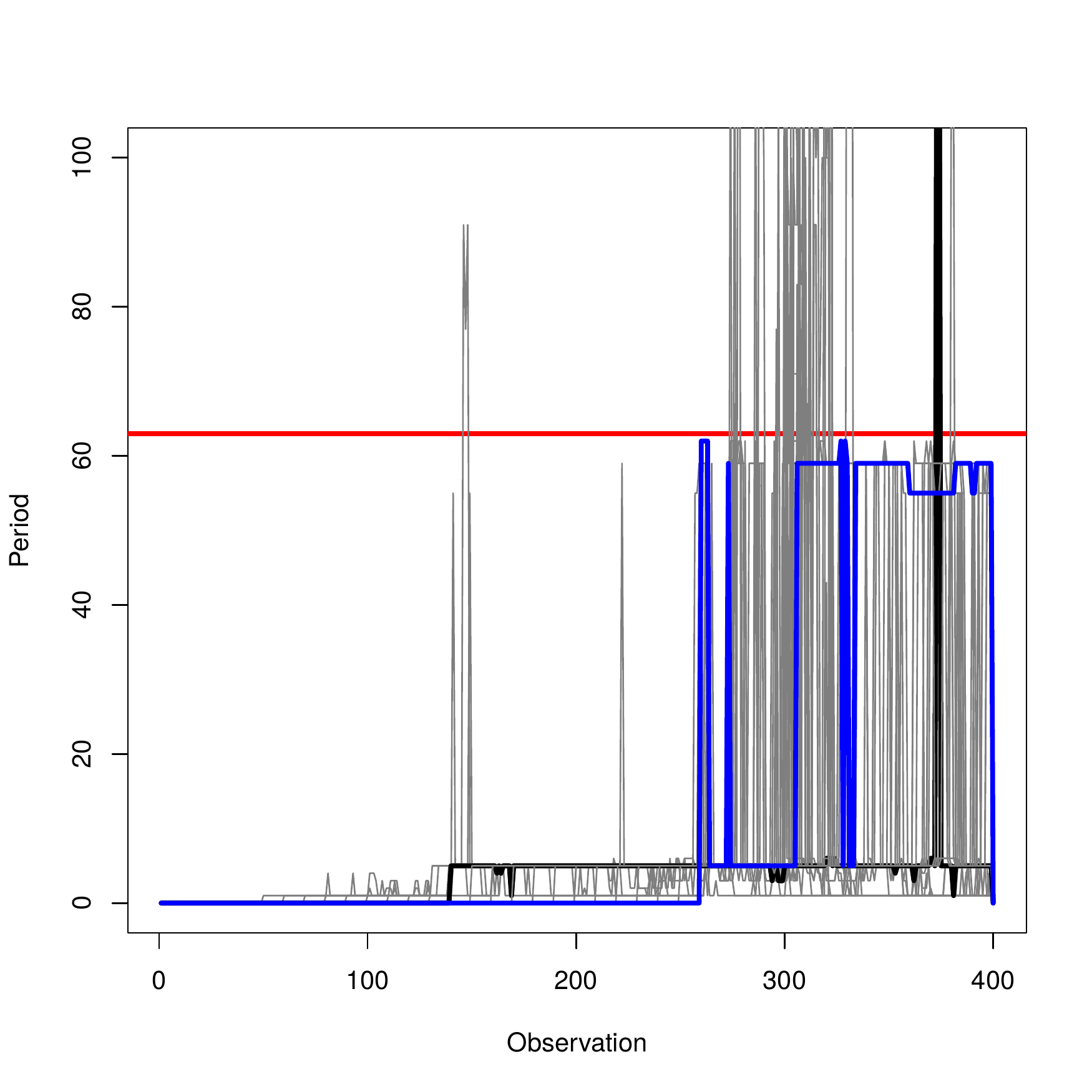}\\
             Panel (c): Noise 3 & Panel (d): Noise 4\\
        \end{tabular}
    \end{center}
    \raggedright
    \footnotesize{Estimation performed for basic sin function, $sin x$, with varying levels of noise. Horizontal axis reports the time period from 0 to 400. Red line provides the true period estimate, whilst a thick black line is plotted for the estimate with a window size of 130. A blue line is added with the maximum window size.}
\end{figure}

However, to really understand what is happening within these estimated periods it is helpful to plot the series of estimates. We present the plots in Figure \ref{fig:sinhynd}. There are ten series for each level of noise. For brevity when plotting we place all series on the same plot showing the series in grey. A thicker black line is plotted using the interval 130 because it is the one that represents the length of two periods and is therefore the one which the TDA approach needs to identify any periodicity. In panel (a) we see that the estimates are quite volatile but there there are often values close to the true $2\pi$ recovered. As the noise level increases the number of such observations falls dramatically. At noise level 3 in panel (c) we see some large overestimates of the period, whilst at noise level 4 in panel (d) there are a few estimates coming close to identifying the true value.

From this first example it is clear that TDA is able to more consistently estimate the period at all levels of noise. At noise levels 1 and 2 TDA almost always reported the period to be 62, whilst at the higher levels TDA is still getting much closer to that true value. An advocacy for the approach developed in this paper is thus found.

\begin{figure}
    \begin{center}
        \caption{Frequency Estimates on Noisy sin Functions with Trend Using \textit{forecast} \cite{hyndman2019package} Package in R}
        \label{fig:sinthynd}
        \begin{tabular}{c c}
             \includegraphics[width=5cm]{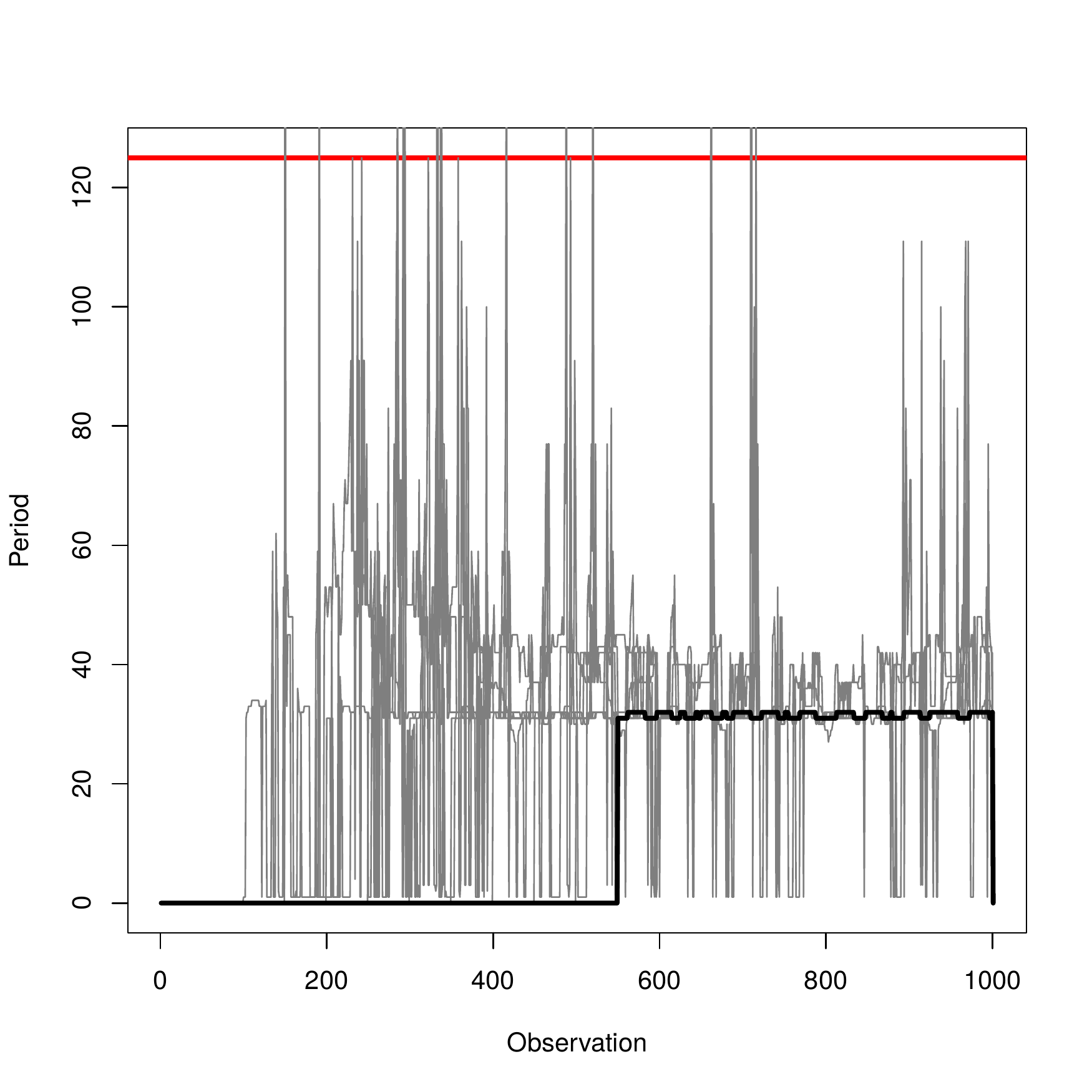}&
             \includegraphics[width=5cm]{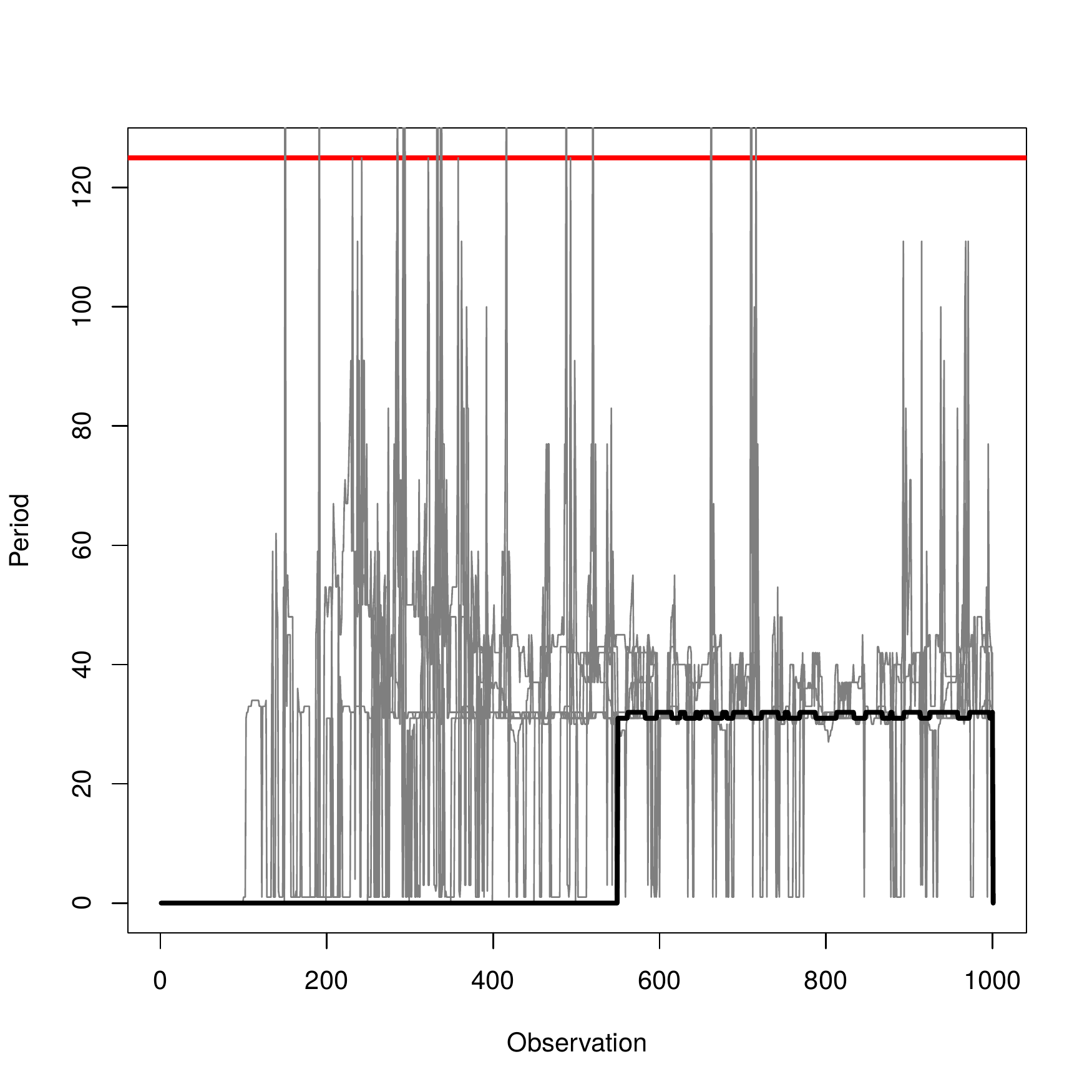}\\
             Panel (a): Noise 1 & Panel (b): Noise 2 \\
             \includegraphics[width=5cm]{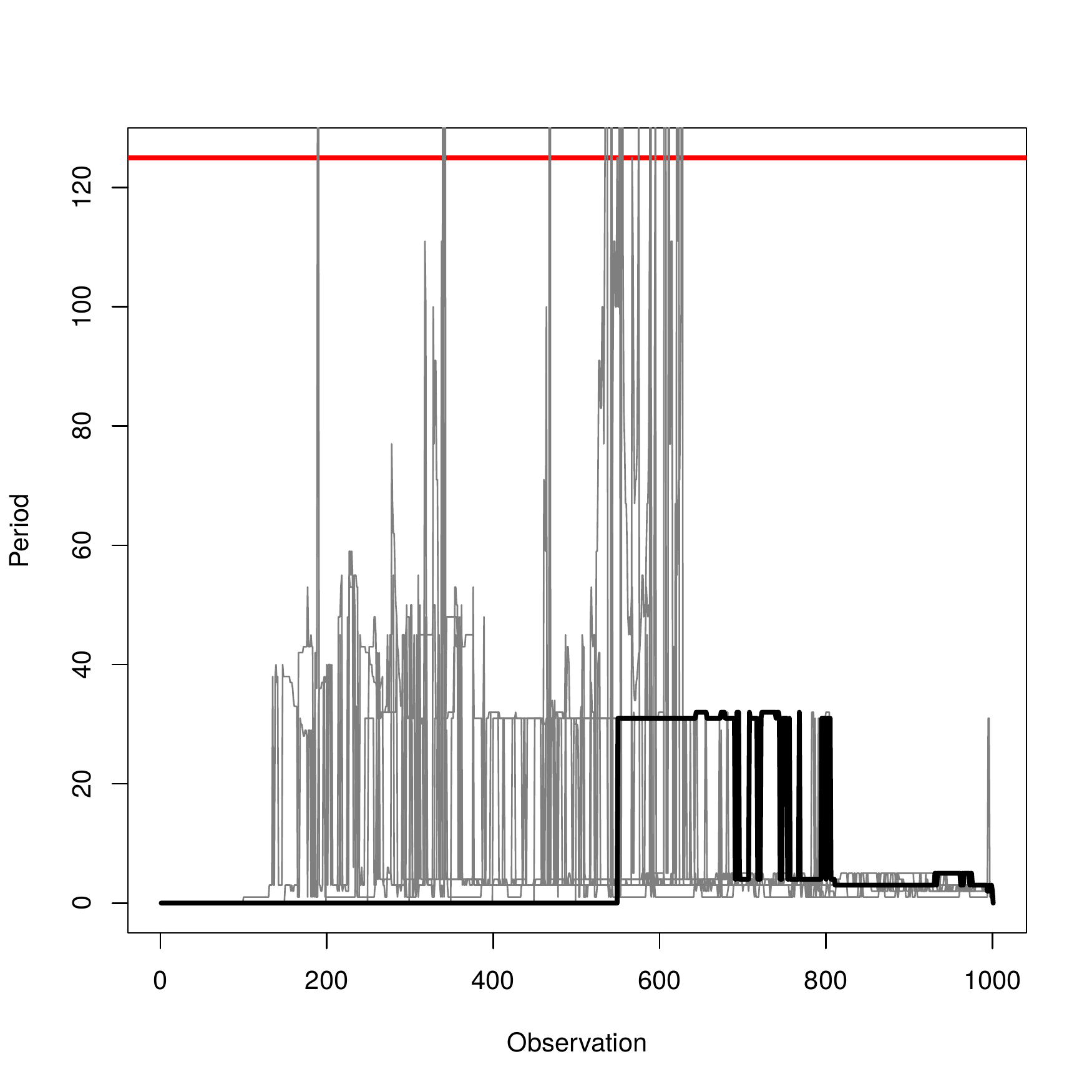}&
             \includegraphics[width=5cm]{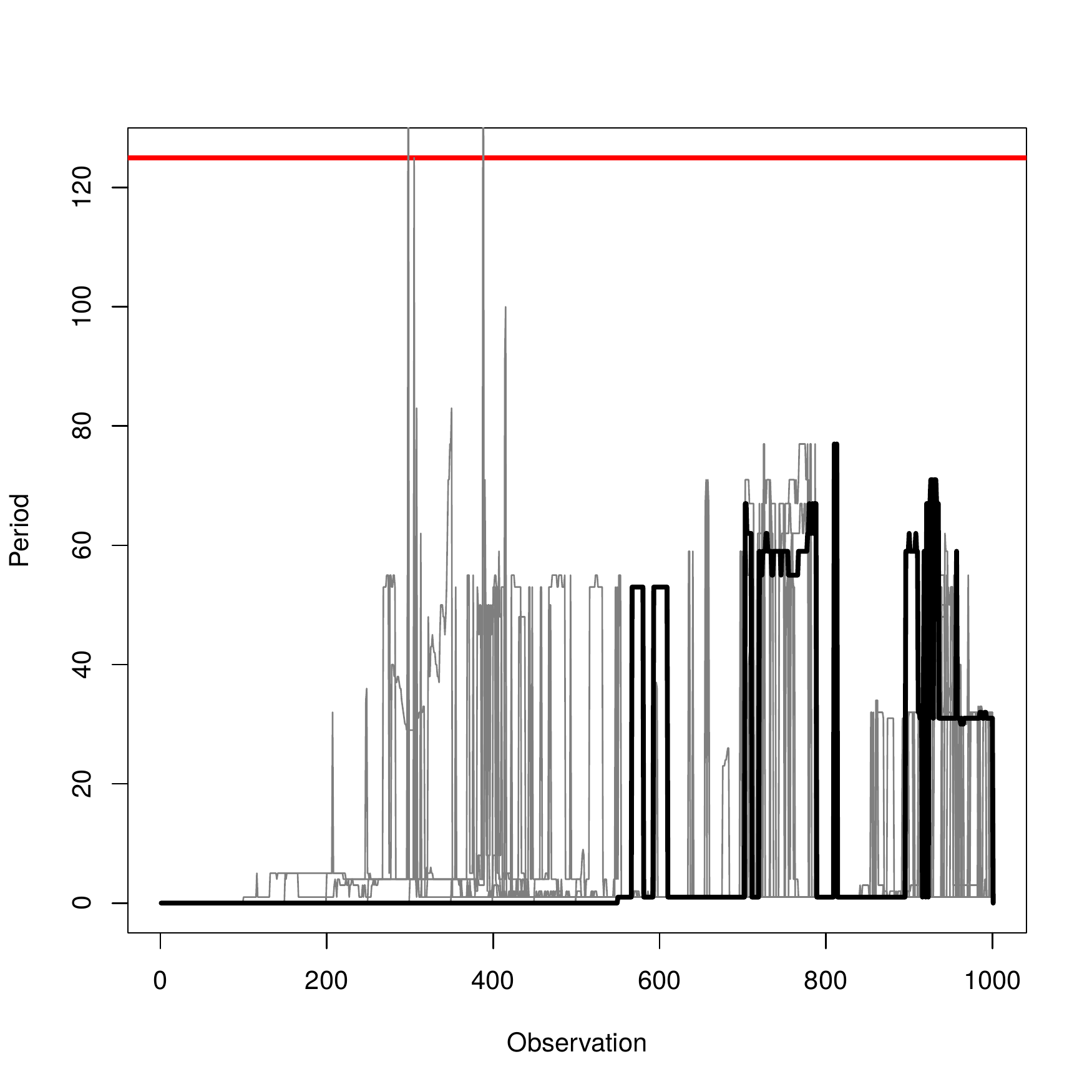}\\
             Panel (c): Noise 3 & Panel (d): Noise 4\\
        \end{tabular}
    \end{center}
    \raggedright
    \footnotesize{Notes: Estimation performed for sin function with varying levels of noise. Horizontal axis reports the time period from 0 to 400. Red line provides the true period estimate, whilst a thick black line is plotted for the estimate with a window size of 550. All other window sizes, 100, 150, 200... etc. are included as thin grey lines.}
\end{figure}

Our second artificial example considers the function $sin(2x)+0.5x$. This is a function with a linear trend, something which the \cite{hyndman2019package} code is designed to deal effortlessly with. Table \ref{tab:sinthynd} gives the estimates and Figure \ref{fig:sinthynd} shows these graphically for an easier visual comparison. This is a longer series so we increase the increment to 50 and provide intervals from 100 to 550.

\begin{table}
    \caption{Period Estimation from Rolling-Windows on Sin Function with Slope}
    \label{tab:sinthynd}
	\begin{tabular}{l c c c c c c c c c  }
		\hline
		&&\multicolumn{4}{l}{Noise 1}&\multicolumn{4}{l}{Noise 2}\\
		Length & Obs & Mean & S.d. & Min & Max & Mean & S.d. & Min & Max \\
		\hline

            100   & 901   & 28.016 & 25.423 & 0     & 499   & 28.016 & 25.423 & 0     & 499 \\
    150   & 851   & 41.971 & 21.571 & 0     & 200   & 41.971 & 21.571 & 0     & 200 \\
    200   & 801   & 37.456 & 12.494 & 0     & 166   & 37.456 & 12.494 & 0     & 166 \\
    250   & 751   & 34.895 & 12.876 & 0     & 250   & 34.895 & 12.876 & 0     & 250 \\
    300   & 701   & 32.036 & 3.317 & 0     & 67    & 32.036 & 3.317 & 0     & 67 \\
    350   & 651   & 31.342 & 1.323 & 0     & 32    & 31.342 & 1.323 & 0     & 32 \\
    400   & 601   & 31.394 & 1.375 & 0     & 32    & 31.394 & 1.375 & 0     & 32 \\
    450   & 551   & 31.380 & 1.427 & 0     & 32    & 31.380 & 1.427 & 0     & 32 \\
    500   & 501   & 31.474 & 1.493 & 0     & 32    & 31.474 & 1.493 & 0     & 32 \\
    550   & 451   & 31.489 & 1.565 & 0     & 32    & 31.489 & 1.565 & 0     & 32 \\
    \hline
		&&\multicolumn{4}{l}{Noise 3}&\multicolumn{4}{l}{Noise 4}\\
		Length & Obs & Mean & S.d. & Min & Max & Mean & S.d. & Min & Max \\
	\hline
    100   & 901   & 6.724 & 20.213 & 0     & 499   & 4.093 & 8.860 & 0     & 40 \\
    150   & 851   & 14.246 & 30.158 & 0     & 333   & 4.513 & 13.178 & 0     & 200 \\
    200   & 801   & 10.741 & 34.891 & 0     & 499   & 4.107 & 12.606 & 0     & 200 \\
    250   & 751   & 8.355 & 25.517 & 0     & 499   & 5.110 & 12.284 & 0     & 59 \\
    300   & 701   & 7.192 & 13.354 & 0     & 200   & 4.843 & 11.023 & 0     & 71 \\
    350   & 651   & 9.147 & 12.510 & 0     & 48    & 6.739 & 13.006 & 0     & 62 \\
    400   & 601   & 9.699 & 11.677 & 0     & 32    & 10.718 & 18.702 & 0     & 77 \\
    450   & 551   & 10.601 & 12.445 & 0     & 53    & 16.313 & 24.832 & 0     & 77 \\
    500   & 501   & 12.721 & 13.053 & 0     & 45    & 18.988 & 24.222 & 0     & 77 \\
    550   & 451   & 15.480 & 13.742 & 0     & 32    & 24.144 & 26.272 & 0     & 77 \\
    \hline
    \end{tabular}
\raggedright
\footnotesize{Notes: Function underlying time series is $sin x$ with noise. Estimates using the \textit{findfrequency} function from R package \textit{forecast} \cite{hyndman2019package} is used on a rolling window basis with the size of the window given in column 1. This then leaves the remaining observations over which the period is calculated. Statistics report the average, standard deviation, minimum and maximum, of the periods that are calculated for the remaining observations.}
\end{table}

Once again Table \ref{tab:sinthynd} reveals the average period estimate to be a long way from the true values. There are a number of zeros within the range that we would be expecting the code to find a period. There is also large variability among estimates evidenced in the standard deviation. As with the sin example it is instructive to see this graphically and so we plot the results in Figure \ref{fig:sinthynd}.

Confirming the findings shown by Figure \ref{fig:sinhynd}, we see that the estimates are consistently below the true value. Indeed it appears that the noise does not make as much difference in this case compared to the simple sin function. However, it is also true that the \cite{hyndman2019package} code does come much closer than it did for the simple sin function.

Across the two artificial examples a very strong case is made for the use of TDA in the identification of periodicity of noisy time series. However, thus far the examples have had simple repetitive periods underlying them. We next move to consider alternatives where that no longer holds.

\section{Discontinuous change of period}

\begin{figure}
    \centering
    \caption{Sequence of Three Periodic Functions} \label{fig:sin_changing_period}
    \includegraphics[scale=0.8]{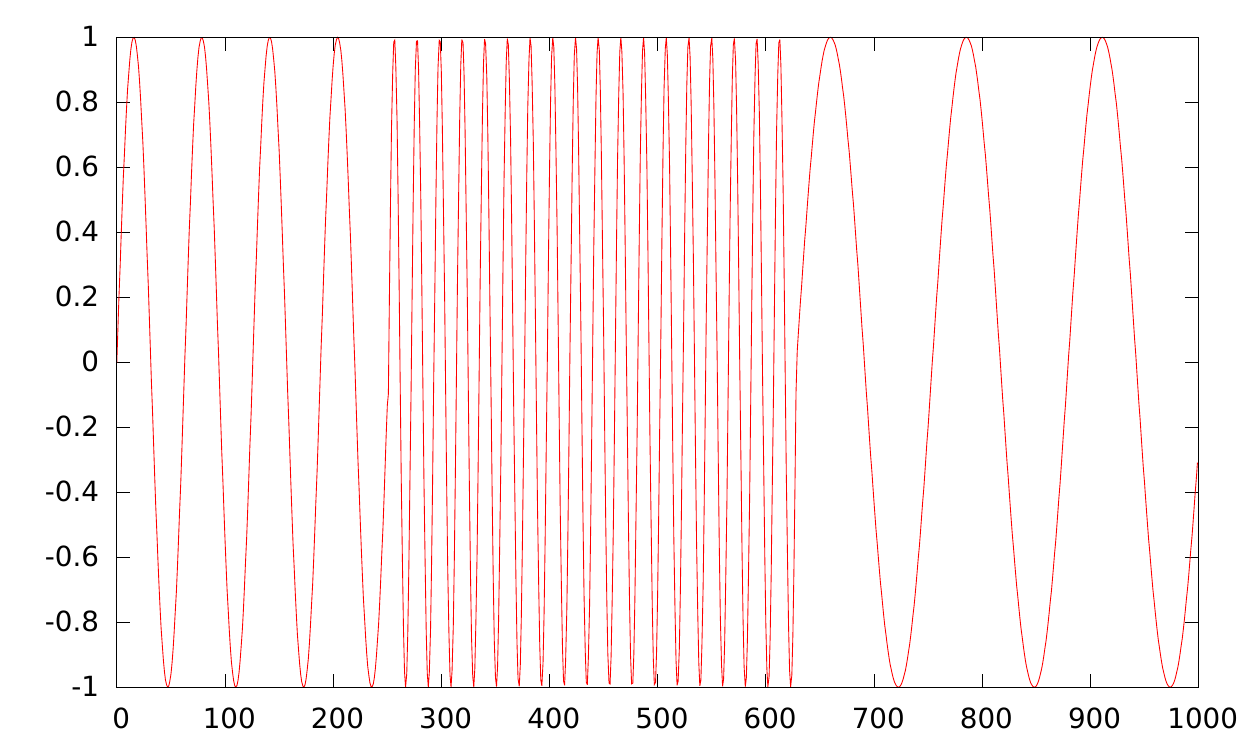}
    \raggedright
    \footnotesize{Notes: From left to right we have $\sin(x)$ followed by $\sin(2x)$ followed by $\sin(\frac{x}{2})$. The transitions are continuous, but not differentiable in the transition points.}
\end{figure}

In many practical cases we may witness a function which is a successive composition of a number of periodic functions. In that case the period is changing in a sharp way between different periodic phases. As an example consider the time series presented in the Figure~\ref{fig:sin_changing_period}. Here there are three clear periodic functions each emerging from the previous. Up to observation 250 the period is 63 $(2\pi)$, with it then accelerating to just 21 $(2\pi/3)$ to observation 625. For the final section it has a period of 125 $(4\pi)$. It should be stressed that the transitions are continuous but not differentiable and that as such any observation numbers are indicative of the transformation points. Note that the presented time series is periodic only when restricted to certain intervals and therefore one cannot assign a \emph{global period} to that function. Such a time series can be found when analyzing motions of motors that can operate on a fixed number of frequencies.

\begin{figure}
    \centering
    \caption{Estimation of Period of the Three-Part Sin Function Presented in Figure~\ref{fig:sin_changing_period}.}
    \label{fig:sin_discontinous_change_of_period_period_estimation.pdf}
    \includegraphics{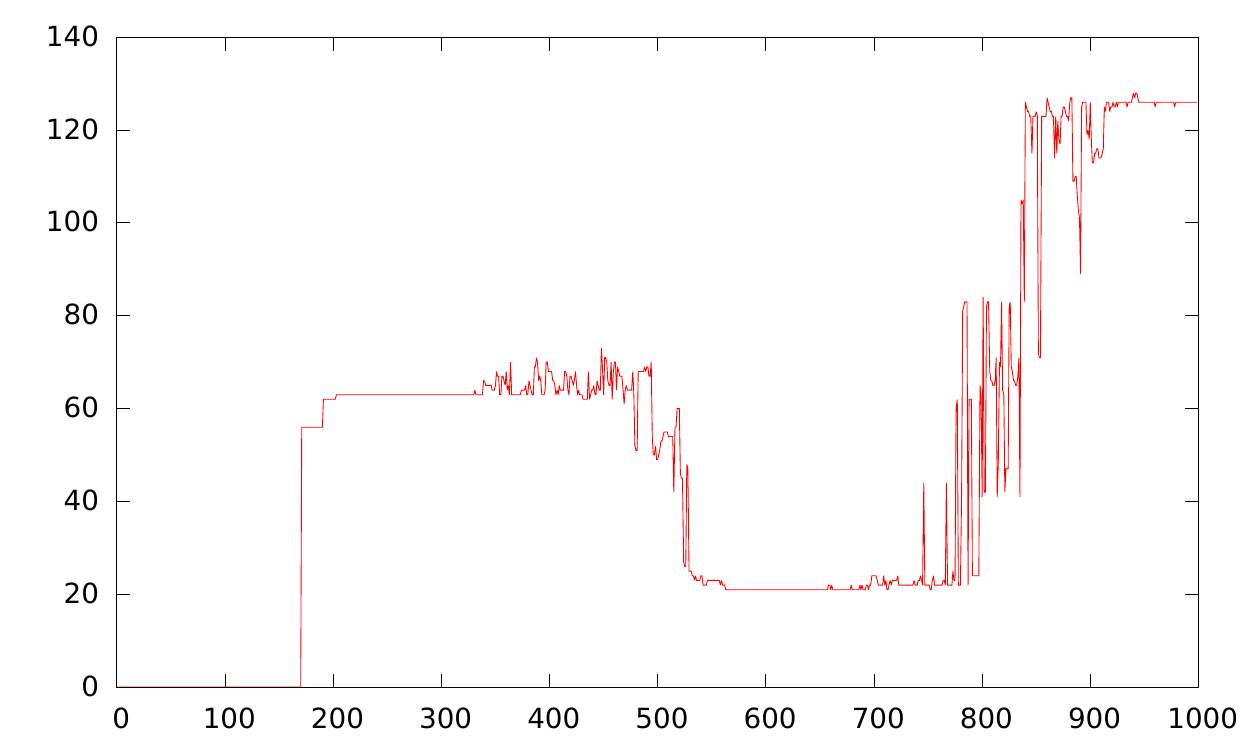}
    \raggedright
    \footnotesize/{Notes: Horizontal axis represents observation number, with vertical axis plotting the period identified by the code.}
\end{figure}

Given that the developed TDA period estimation procedure only takes local information into account it is possible to detect the regions of periodicity of functions like the one presented in the Figure~\ref{fig:sin_changing_period}. The period estimation for this case is presented in the Figure~\ref{fig:sin_discontinous_change_of_period_period_estimation.pdf}. Two observations emerge immediately, firstly the TDA is able to identify all three subperiods, but it does so with an expected lag meaning that any application of the approach in nowcasting would be restricted by such.

\begin{figure}
    \begin{center}
        \caption{Period Estimates for Three-Part Sin Function Presented in Figure \ref{fig:sin_changing_period}}
        \label{fig:sindhynd}
        \begin{tabular}{c c}
            \includegraphics[width=5cm]{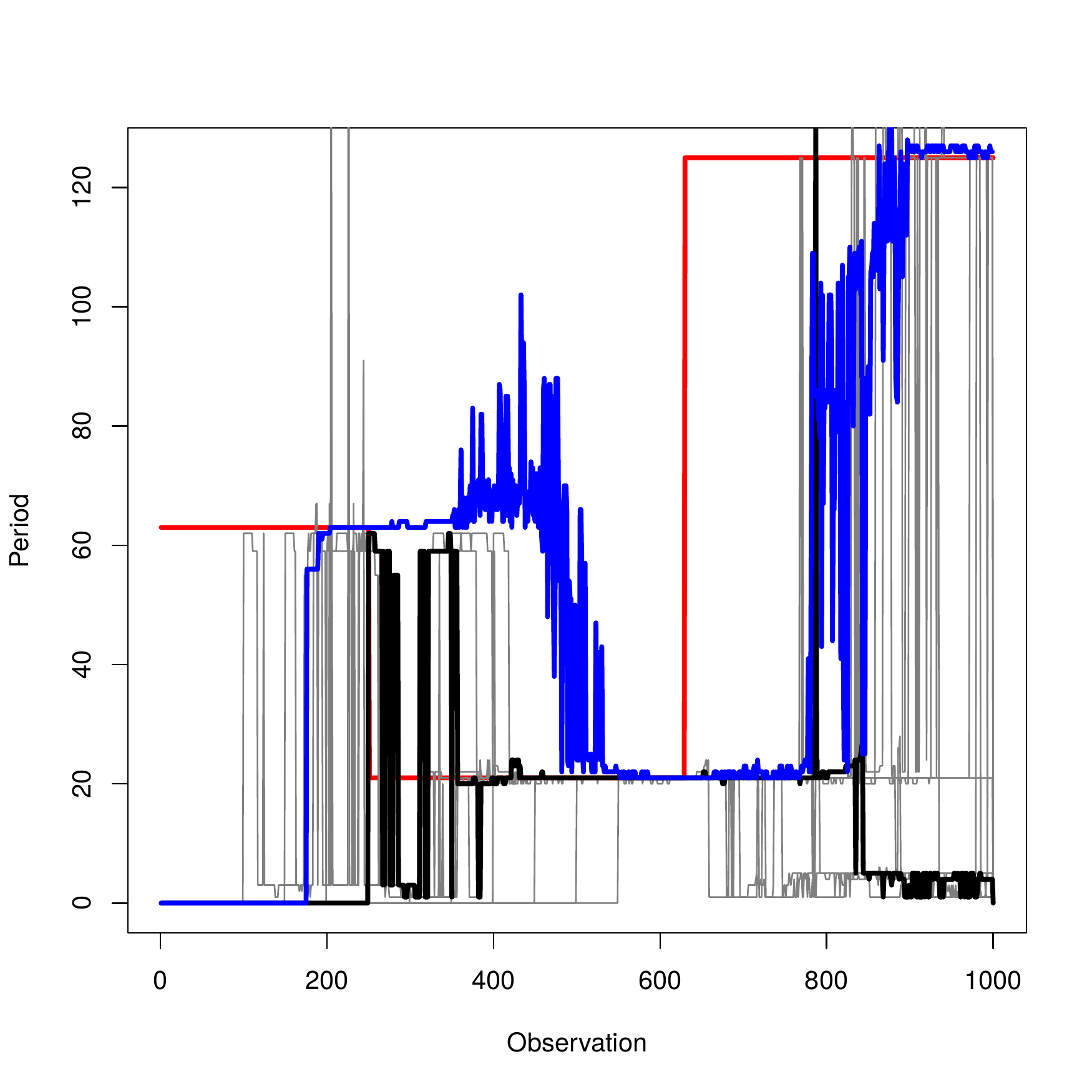} & \includegraphics[width=5cm]{sind1hynd.pdf} \\
            Panel (a): Noise 1 & Panel (b): Noise 2\\
             \includegraphics[width=5cm]{sind1hynd.pdf}& \includegraphics[width=5cm]{sind1hynd.pdf}\\
             Panel (c): Noise 3 & Panel (d): Noise 4 \\
        \end{tabular}
    \end{center}
    \raggedright
    \footnotesize{Notes: Figures show periodicity estimates obtained using our TDA procedure as thick blue lines. Values from the \cite{hyndman2019package} \textit{findfrequency} function are plotted as grey lines. Estimates for the period in the last method are obtained from rolling windows of width 100 to 550 in increments of 50 observations. The result with a window width of  550 are plotted as a thick black line. A red line is added with the true period of the function. All other rolling window sizes are plotted as thin grey lines.}
\end{figure}

As one can observe, the estimation of period is accurate in the regions far away from transitions between periods. In this case it the period is estimated to $63$, $21$ and $134$. As in the previous sections we can introduce noise over the top of the functions and observe that the results are highly robust with respect to such random perturbations. Estimates for the period arising from the TDA are presented in Figure \ref{fig:sindhynd} as blue lines.

Evaluating the benefits of our proposed TDA approach we present comparison results from the \cite{hyndman2019package} \textit{findfrequency} function. In this case a table is not meaningful as summary statistics will not identify the three separate periods so we proceed straight to the graphical analysis. Figure \ref{fig:sindhynd} shows clearly how TDA, blue line, identifies all three periods with all four of the noise levels.

Given the low number of observations at which the series has a particular period, and the number of parameter values we include from the \cite{hyndman2019package} code, the plots in Figure \ref{fig:sindhynd} are quite busy. However, the key results come through clearly from all four panels. TDA finds all three periods, but does so after a lag. In some cases the \textit{findfrequency} function quickly notices the change in period, moving to the new value within a few observations. For the spectral approach the challenge is identifying all three periods, especially the final long period. Across all of the work we find the spectral method needs more than two cycles to detect periodicity. In this section we have thus demonstrated more support for TDA approaches, but with a call to obtain a quicker adjustment when periods change.

\section{Applications}
\label{sec:apply}

Motivation for developing a TDA toolkit for periodicity analysis lies in the multitude of applications in which there are noisy time series with dynamics of interest to researchers and practitioners. In this section we consider three real life examples of time series: temperature, monthly number of sunspots data, and monthly shipment volumes from the U.S.

\subsection{Melbourne Minimum Temperature Data}

Away from the equator it is well understood that the temperature cycles through cold winters and hot summers, producing a clear annual periodicity. Further there will always be noise within the temperature cycle caused by unusual weather events, variations in global weather patterns, human activity, and many other non-periodic influences. As a first demonstration of the power of the proposed TDA methodology a temperature time series makes clear sense. Specifically let us analyze the minimal daily temperature in Melbourne, Australia in the period from 1 Jan. 1981-1 Jan. 1989. This daily data can be found on a number of machine learning example data web-pages\footnote{Our version of the data is taken from \url{https://www.kaggle.com/paulbrabban/daily-minimum-temperatures-in-melbourne}.}.

\begin{figure}
    \begin{center}
    \caption{Daily Minimum Temperatures in Melbourne.}
    \label{fig:tdamelb}
    \begin{tabular}{c}
         \includegraphics[scale=0.8]{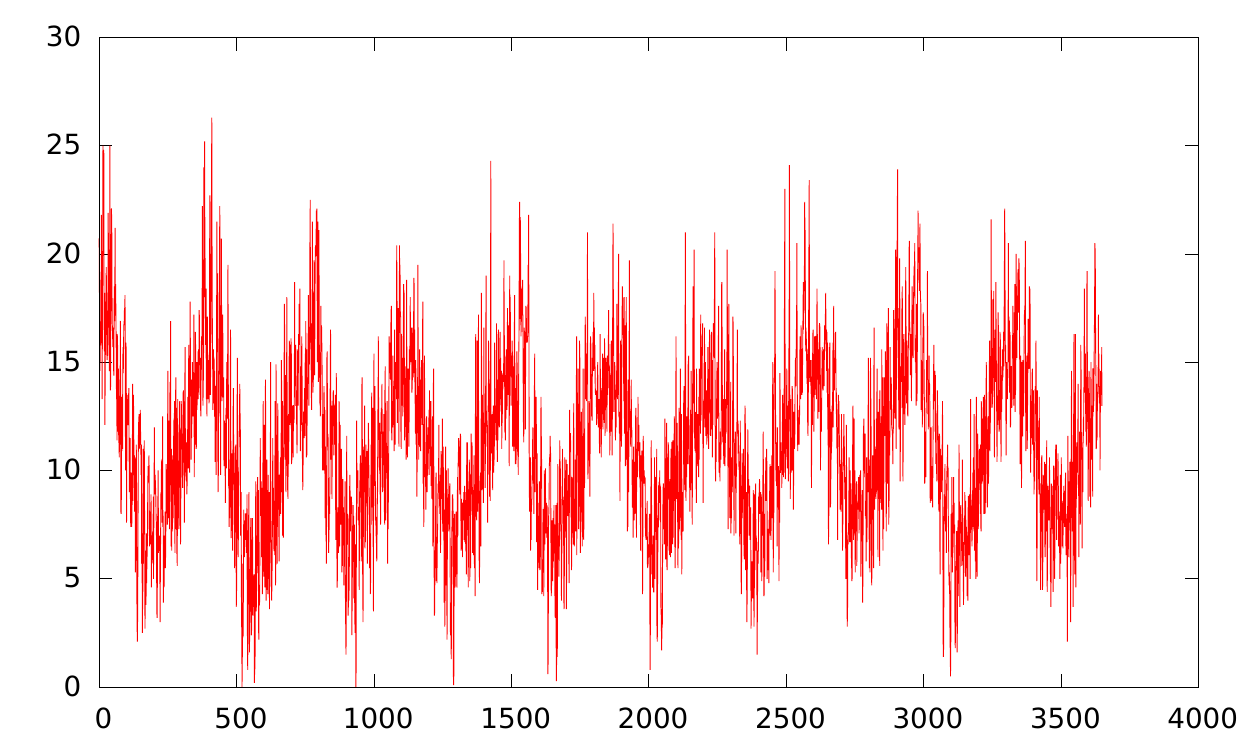}\\
         Panel (a): Time Series of Minimum Temperatures \\
         \includegraphics[scale=0.8]{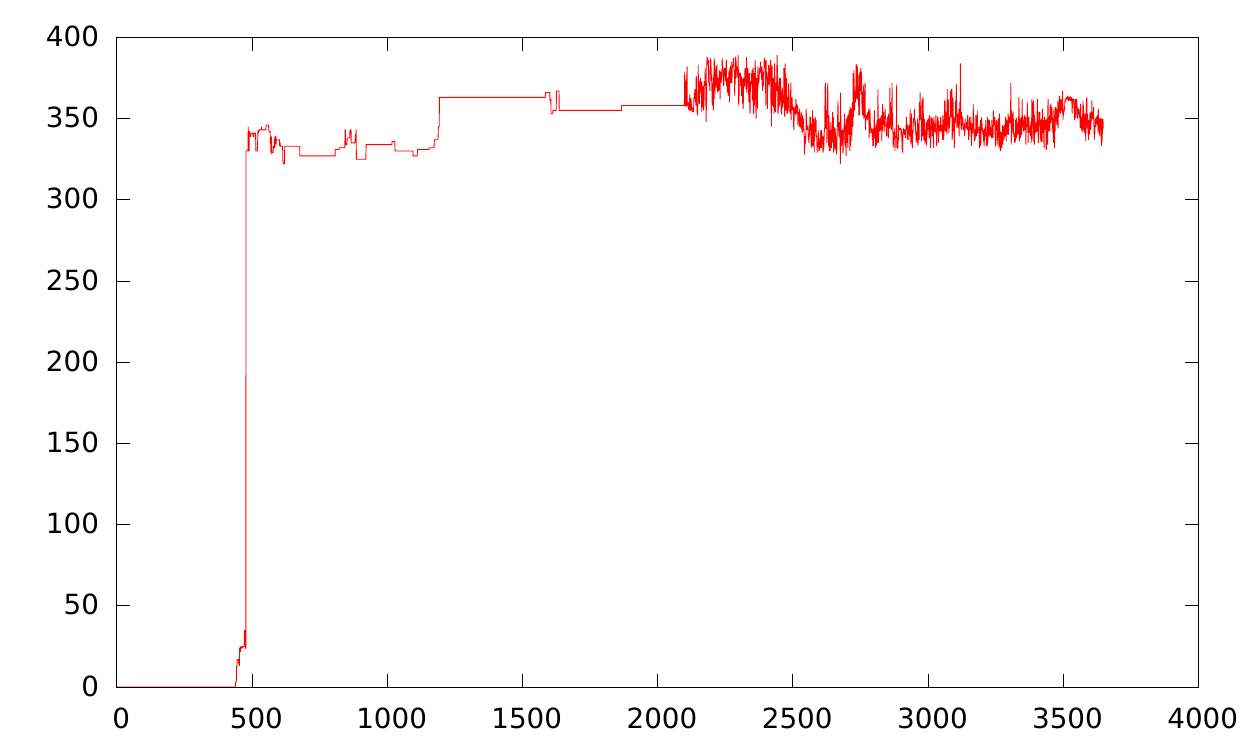}\\
         Panel (b): Period Estimation\\
    \end{tabular}
    \end{center}
    \raggedright
    \footnotesize{Notes: Panel (a) provides a time series plot of minimum daily temperature recorded in Melbourne Australia between 1st January 1981 and 1st January 1989. Panel (b) plots the corresponding period estimate obtained through the TDA procedure developed in this paper.}
\end{figure}

Figure \ref{fig:tdamelb} plots the time series clearly demonstrating the stylised facts of temperature data. We can roughly spot the yearly
pattern in the data with the series starting, and ending, at a similar level. Further the series is clearly subjected to a large variation around the main period. In this way it shares many characteristics of the noisy sin functions discussed in Section \ref{sec:noisy_time_series}. The lower panel of Figure \ref{fig:tdamelb} shows the estimated period from our TDA approach. Although not returning 365 as the number of days in the year, the period estimates are consistently around 350 and often push upwards towards 370. For much of the time we are able to recover an almost constant estimate, but this is replaced by a noiser estimate in the final years of the sample.

\begin{table}[]
    \begin{center}
        \caption{Period Estimation for Melbourne Daily Temperature Data Using \cite{hyndman2019package}}
        \label{tab:melbhynd}
        \begin{tabular}{l l l l l l}
            \hline
            Length & Obs & Mean & S.d. & Min & Max \\
            \hline
             100   & 3550  & 4.241 & 7.055 & 0     & 166 \\
    200   & 3450  & 3.427 & 10.210 & 0     & 499 \\
    300   & 3350  & 6.312 & 3.889 & 0     & 25 \\
    400   & 3250  & 10.093 & 4.444 & 0     & 24 \\
    500   & 3150  & 12.901 & 4.703 & 0     & 23 \\
    600   & 3050  & 13.911 & 4.669 & 0     & 30 \\
    700   & 2950  & 14.745 & 5.205 & 0     & 34 \\
    800   & 2850  & 15.781 & 5.995 & 0     & 32 \\
    900   & 2750  & 16.156 & 5.950 & 0     & 31 \\
    1000  & 2650  & 16.767 & 5.926 & 0     & 31 \\
    1100  & 2550  & 16.840 & 5.859 & 0     & 32 \\
    1200  & 2450  & 16.954 & 5.730 & 0     & 34 \\
    1300  & 2350  & 17.777 & 6.166 & 0     & 31 \\
    1400  & 2250  & 18.471 & 6.188 & 0     & 31 \\
    1500  & 2150  & 18.918 & 6.187 & 0     & 32 \\

             \hline
        \end{tabular}
    \end{center}
    \raggedright
    \footnotesize{Notes: Periods estimated on a rolling window basis using the \textit{findfrequency} from the R package \textit{forecast} \cite{hyndman2019package}. Length reports the size of the rolling window used to estimate the periodicity, Obs. then details the number of repetitions of the window used to calculate the summary statistics in the final four columns. Mean is the average period estimate, with S.d. reporting the standard deviation thereof. Min and Max provide the minimum and maximum estimated period.}
\end{table}
\begin{figure}
    \begin{center}
        \caption{Period Estimates for Melbourne Daily Temperature Data using \cite{hyndman2019package}}
        \label{fig:melbhynd}
        \includegraphics[width=8cm]{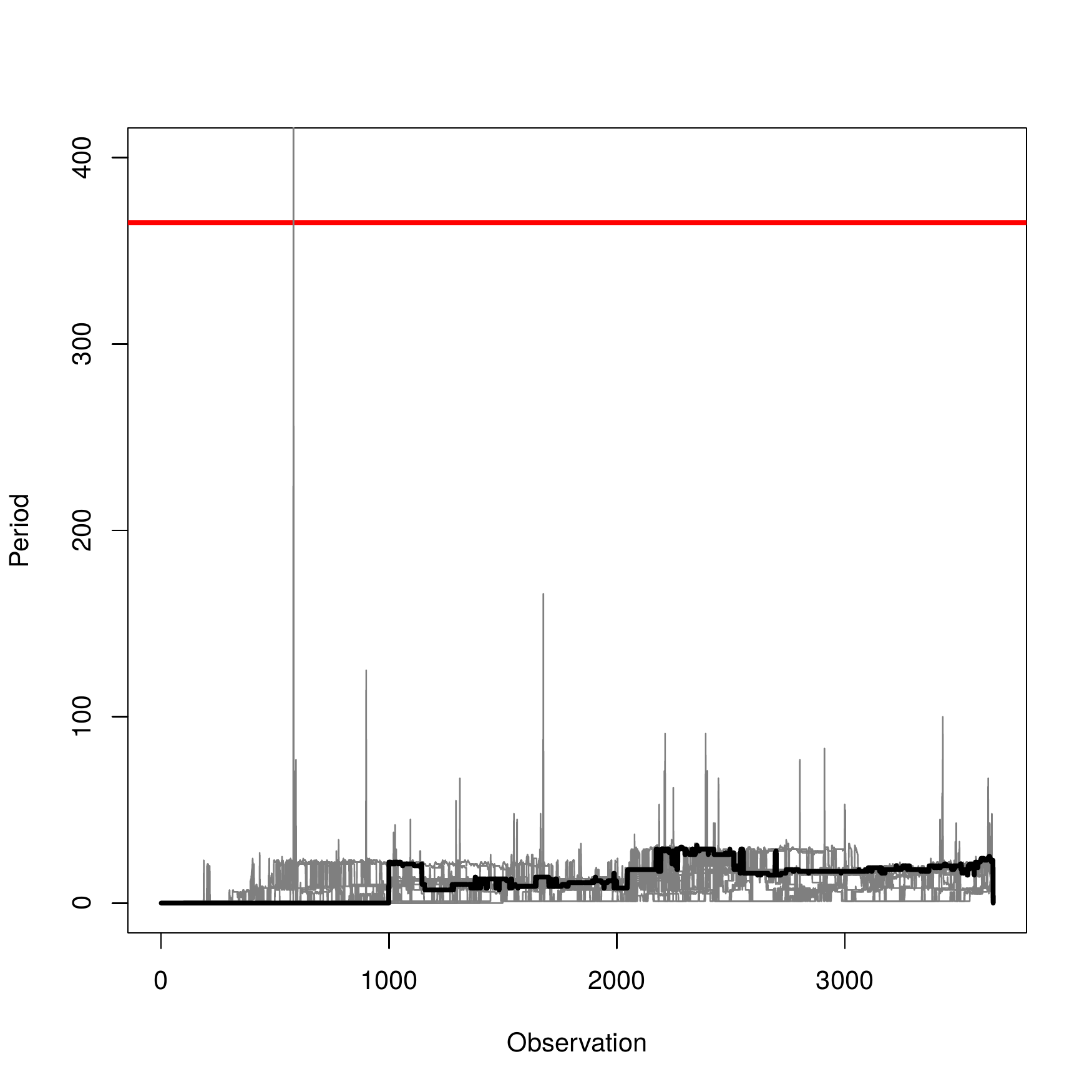}
    \end{center}
    \raggedright
    \footnotesize{Notes: Periods estimated on a rolling window basis using the \textit{findfrequency} from the R package \textit{forecast} \cite{hyndman2019package}. The period estimates from a rolling window of 1000 are plotted as a solid line. All other estimates are plotted using thinner grey lines.}
\end{figure}

Table \ref{tab:melbhynd} presents the estimates from the \textit{findfrequency} function estimations of the period of the Melbourne minimum daily temperature series. We ran the code with windows between 100 and 1500 days at increments of 100 days. Extending beyond 700 days brought the average period estimate higher, but in every case the values suggested were less than one month, a long way short of the true annual periodicity. Likewise, the difference between our proposed TDA procedure and the approach developed in \cite{hyndman2019package} is also large.

For this dataset there is a large amount of noise and, like with the artificial example earlier, this makes uncovering the true period difficult. However, we have demonstrated the TDA approach that is proposed within this paper produces much more accurate predictions about periodicity than do spectral methods of the type calculated by \cite{hyndman2019package}'s \textit{forecast} package in R. There are a number of other temperature datasets available, but each shares similar characteristics to the Melbourne series and develops similar results. We therefore posit that for climatic analysis TDA offers a strong option for recovering dynamic behaviour.

\subsection{Zurich Sun Activity Data}

In the second example we consider a number of Sun spots reported in Zurich between 1749 and 1983. This is an example of a series with a known periodicity and a sufficiently long history to enable meaningful comparison possible. With more than 2500 months of data there is sufficient range to make meaningful use of a monthly frequency series. Another classic example for machine learning, the data set can be found online at many sources\footnote{Our data was downloaded from~\url{https://github.com/PacktPublishing/Practical-Time-Series-Analysis/blob/master/Data\%20Files/monthly-sunspot-number-zurich-17.csv}.}.

\begin{figure}
    \begin{center}
    \caption{Zurich Sunspot Data}
    \label{fig:zurich}
    \begin{tabular}{c}
    \includegraphics[scale=0.5]{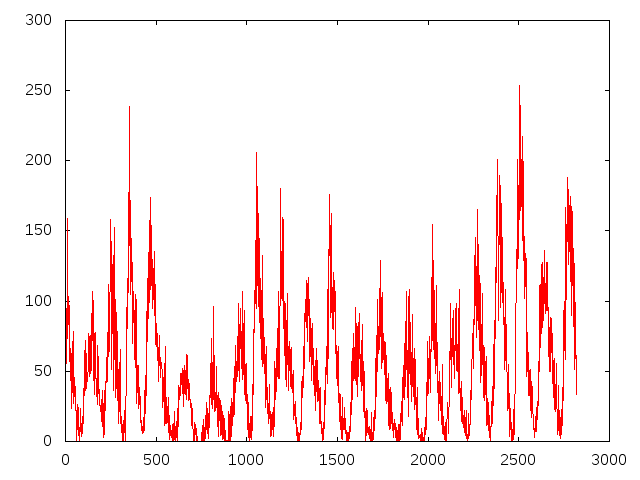}\\
    Panel (a): Number of sunspots measured in month\\
    \includegraphics[scale=0.5]{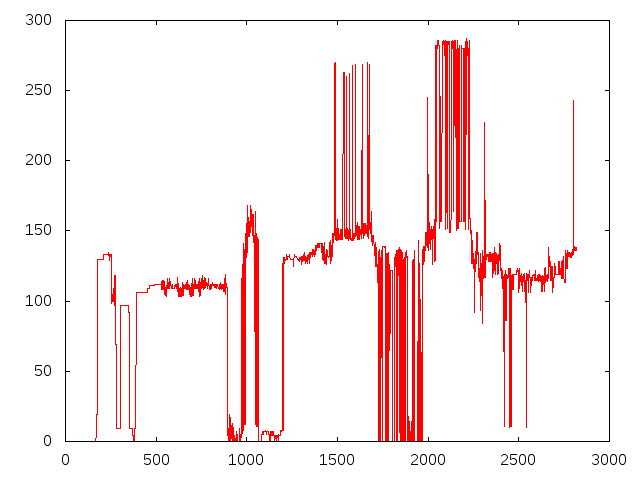}\\
    Panel (b): Estimate of period using TDA\\
    \end{tabular}
    \end{center}
    \raggedright
    \footnotesize{Notes: Panel (a) presents time series of Zurich sunspot activity. Panel (b) provides the estimated period from the TDA approach developed in this paper. Observation numbers on the horizontal axis facilitate approximate identification of the period by eye, and alignment of the two plots.}
\end{figure}

Again aligning the series with the TDA approach periodicity estimate we can see that for large sections of the time series the identification steps presented herein are able to recover a period close to the expected 125. Relative to previous examples we do see more variability in these estimates, with occasion where the value increases towards two actual periods, and other occasions where the code is unable to indentify any periodicity at all.

\begin{table}
    \begin{center}
        \caption{Period Estimation for Zurich Sun Activity Data Using \cite{hyndman2019package}}
        \label{tab:zurhynd}
        \begin{tabular}{l l l l l l}
            \hline
            Length & Obs & Mean & S.d. & Min & Max \\
            \hline
             100   & 2720  & 4.510 & 16.010 & 0     & 499 \\
    200   & 2620  & 13.595 & 36.695 & 0     & 499 \\
    300   & 2520  & 38.354 & 53.631 & 0     & 250 \\
    400   & 2420  & 80.478 & 83.899 & 0     & 998 \\
    500   & 2320  & 110.276 & 91.388 & 0     & 998 \\
    600   & 2220  & 113.236 & 63.319 & 0     & 998 \\
    700   & 2120  & 113.978 & 59.899 & 0     & 998 \\
    800   & 2020  & 115.388 & 62.943 & 0     & 998 \\
    900   & 1920  & 111.225 & 55.084 & 0     & 499 \\
    1000  & 1820  & 107.266 & 56.214 & 0     & 333 \\
             \hline
        \end{tabular}
    \end{center}
    \raggedright
    \footnotesize{Notes: Periods estimated on a rolling window basis using the \textit{findfrequency} from the R package \textit{forecast} \cite{hyndman2019package}. Length reports the size of the rolling window used to estimate the periodicity, Obs. then details the number of repetitions of the window used to calculate the summary statistics in the final four columns. Mean is the average period estimate, with S.d. reporting the standard deviation thereof. Min and Max provide the minimum and maximum estimated period.}
\end{table}
\begin{figure}
    \begin{center}
        \caption{Period Estimates for Zurich Sun Activity Data using \cite{hyndman2019package}}
        \label{fig:zurhynd}
        \includegraphics[width=8cm]{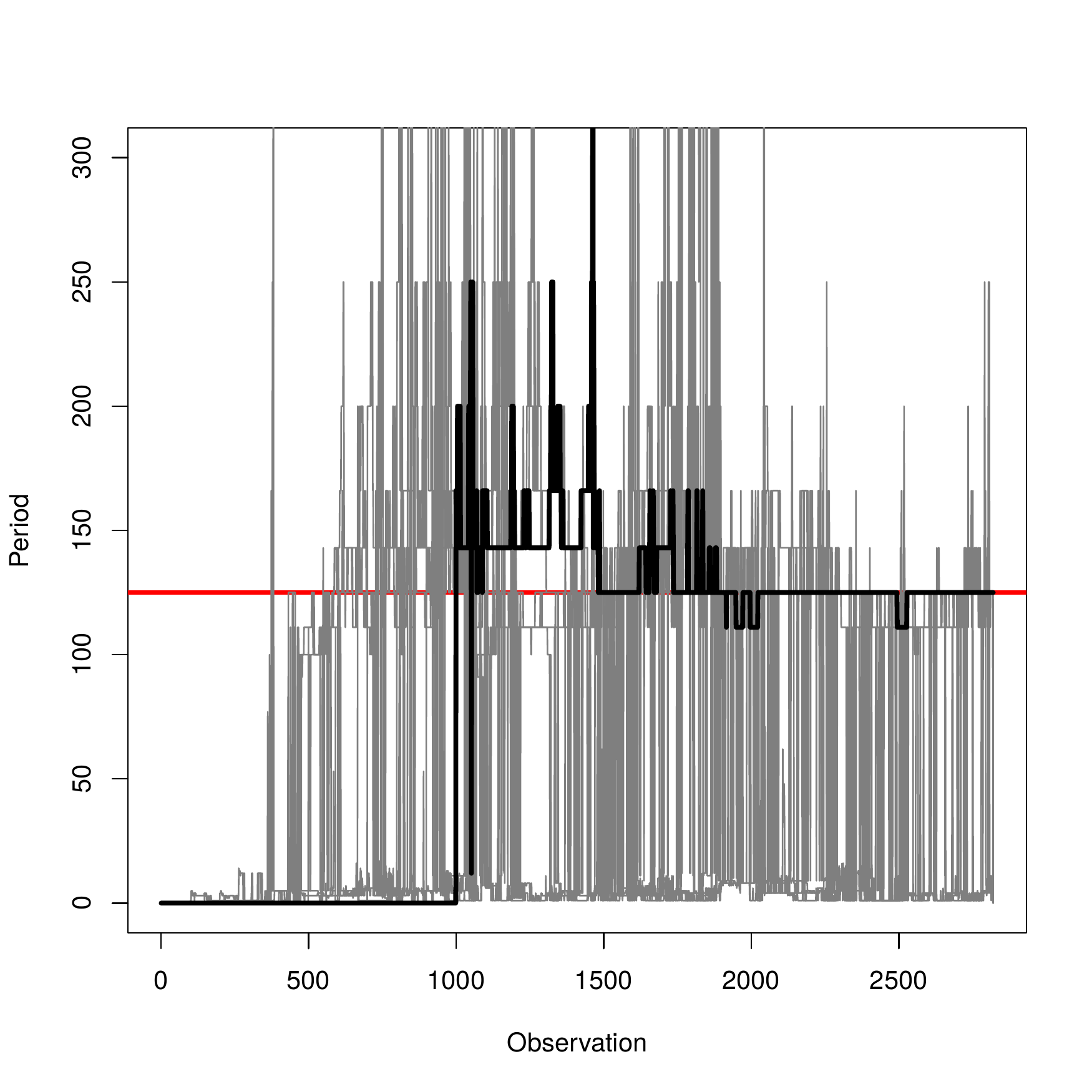}
    \end{center}
    \raggedright
    \footnotesize{Notes: Periods estimated on a rolling window basis using the \textit{findfrequency} from the R package \textit{forecast} \cite{hyndman2019package}. A solid line is plotted for the period estimate on a rolling window of 1000. All other windows are plotted as lighter grey lines.}
\end{figure}

As in past examples Table \ref{tab:zurhynd} shows how it takes a large number of observations before the \cite{hyndman2019package} \textit{forecast} function is able to recover the true period of the function. In this example we see consistency that is absent in the TDA. Shorter windows are not sufficient, but the 1000 observation case shown as the thick line in Figure \ref{fig:zurhynd}, recovers 125 from around observation 2000. Overall the TDA produces more estimates close to the true value but the volatility in this case offers some favour to the established \cite{hyndman2019package} approach.

\subsection{Shipping Volumes}

\begin{figure}
    \begin{center}

    \caption{CASS Shipping Volumes and TDA Period Estimate}
    \label{fig:shipment_time_series}
    \begin{tabular}{c}
    \includegraphics[scale=0.5]{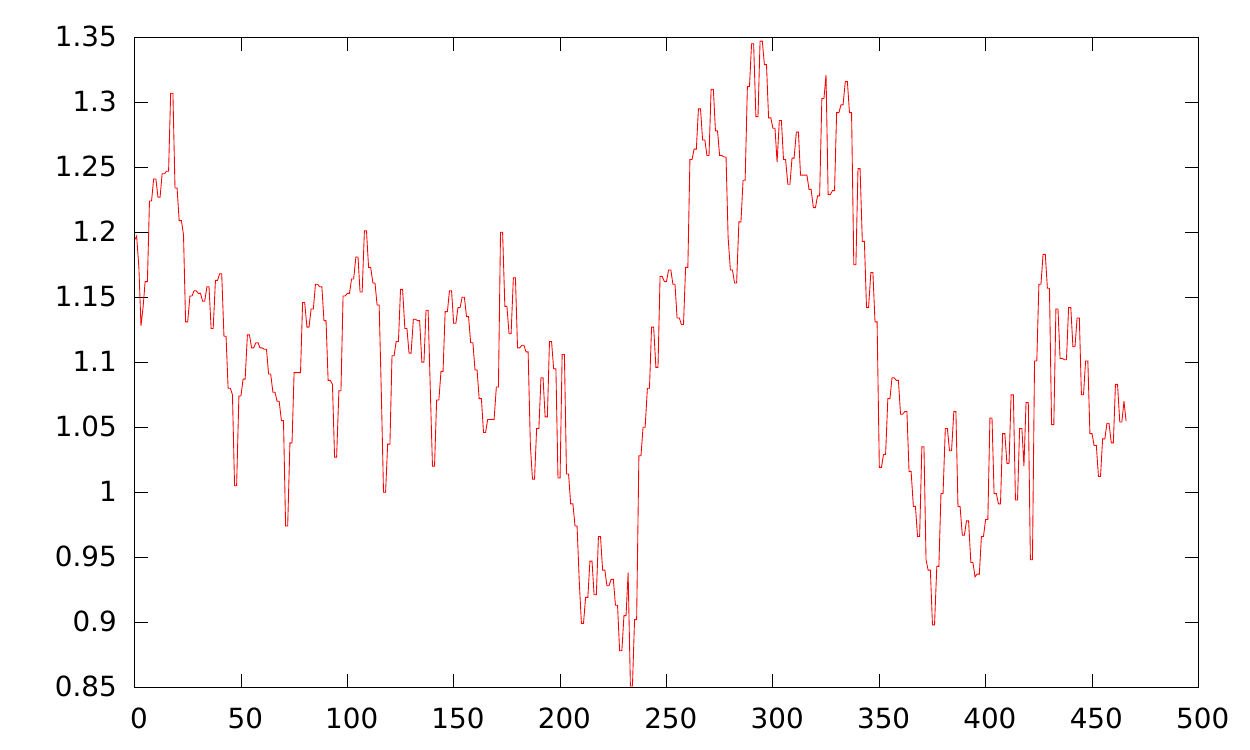}\\
    Panel (a): Shipment Volume Time Series \\
    \includegraphics[scale=0.5]{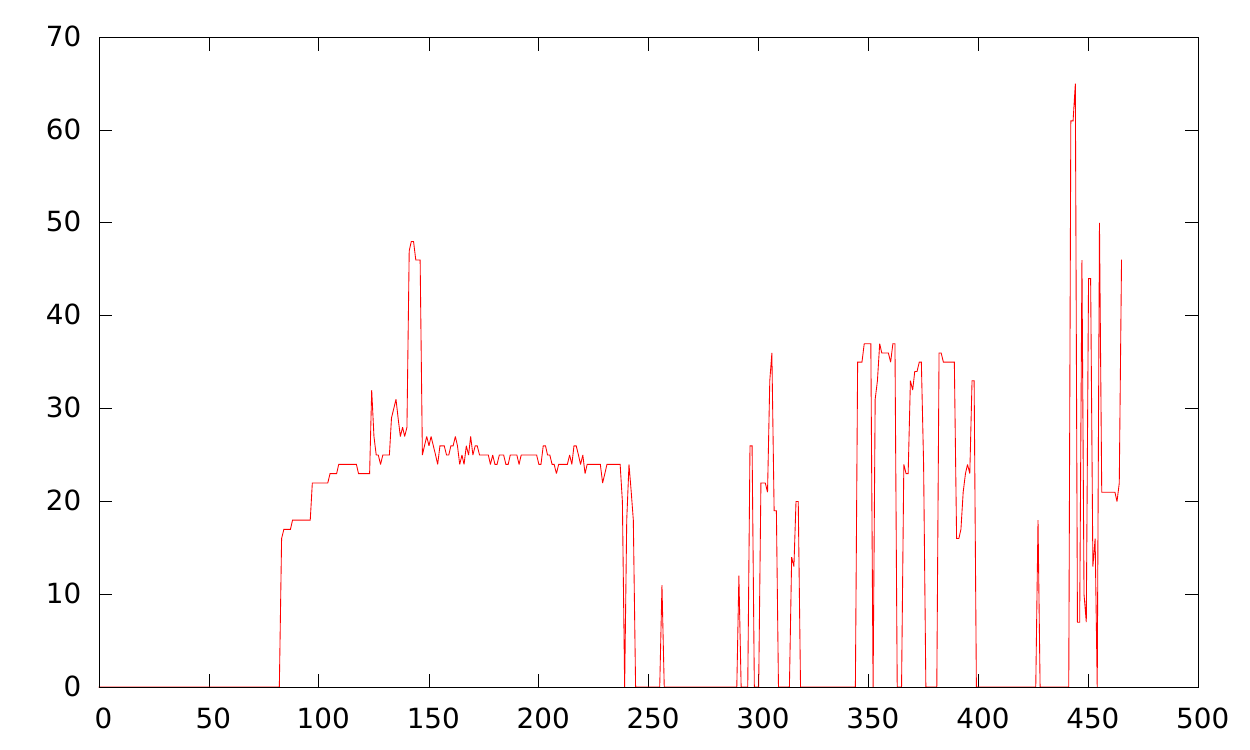}\\
    Panel (b): Estimated Period from TDA Approach\\
    \end{tabular}
    \end{center}
    \raggedright
    \footnotesize{Notes: Panel (a) plots the volume of shipping recorded by the CASS system of freight measurement. Panel (b) provides the estimated period using the TDA approach outlined in this paper. Data for this example sourced from CASS via the Quandl system \url{https://www.quandl.com/data/CASS-Freight-Indices}.}

\end{figure}

The last time series to consider is the shipment dataseries which can be found in
\url{https://www.quandl.com/data/CASS-Freight-Indices}. This index is constructed from company data in the United States to form a measure of the volume, and value, of trade that is passing through American ports. Data is based upon the month that it is processed by CASS and not the time at which the transaction takes place. Nevertheless it is still clear that there are cyclicalities within the data. We plot the time series of the export shipment volumes in Figure~\ref{fig:shipment_time_series}. Using 1.00 as the base for the index means the values of the series do not change by much, but the topological interpretation remains the same as the wider ranging series earlier in this paper.

Our TDA approach reveals a period of around 25 for a long block of the time series after observation 100. This ends around the time that the big drop is seen in the value of the index. Panel (a) of Figure \ref{fig:shipment_time_series} shows a peak in the index emerging from the drop which ends around observation 250. TDA then identifies another short time span over which there is some periodic behaviour towards observation 400. The period estimate in panel (b) for this second block is slightly longer than the first at just over 30. An inspection of panel (a) of Figure \ref{fig:shipment_time_series} suggests that the estimates of the TDA approach are consistent with the observed series. However, there is little immediate motivation for the pattern of almost two year cyclicality appearing in monthly data.

\begin{table}
    \begin{center}
        \caption{Period Estimation for Zurich Sun Activity Data Using \cite{hyndman2019package}}
        \label{tab:shiphynd}
        \begin{tabular}{l l l l l l}
            \hline
            Length & Obs & Mean & S.d. & Min & Max \\
            \hline
            10    & 457   & 0.998 & 0.047 & 0     & 1 \\
            20    & 447   & 0.998 & 0.047 & 0     & 1 \\
            30    & 437   & 0.998 & 0.048 & 0     & 1 \\
    40    & 427   & 0.998 & 0.048 & 0     & 1 \\
    50    & 417   & 0.998 & 0.049 & 0     & 1 \\
    60    & 407   & 0.998 & 0.050 & 0     & 1 \\
    70    & 397   & 0.997 & 0.050 & 0     & 1 \\
    80    & 387   & 0.997 & 0.051 & 0     & 1 \\
    90    & 377   & 0.997 & 0.051 & 0     & 1 \\
    100   & 367   & 0.997 & 0.052 & 0     & 1 \\

             \hline
        \end{tabular}
    \end{center}
    \raggedright
    \footnotesize{Notes: Periods estimated on a rolling window basis using the \textit{findfrequency} from the R package \textit{forecast} \cite{hyndman2019package}. Length reports the size of the rolling window used to estimate the periodicity, Obs. then details the number of repetitions of the window used to calculate the summary statistics in the final four columns. Mean is the average period estimate, with S.d. reporting the standard deviation thereof. Min and Max provide the minimum and maximum estimated period.}
\end{table}
\begin{figure}
    \begin{center}
        \caption{Period Estimates for CASS Shipping Data using \cite{hyndman2019package}}
        \label{fig:shiphynd}
        \includegraphics[width=8cm]{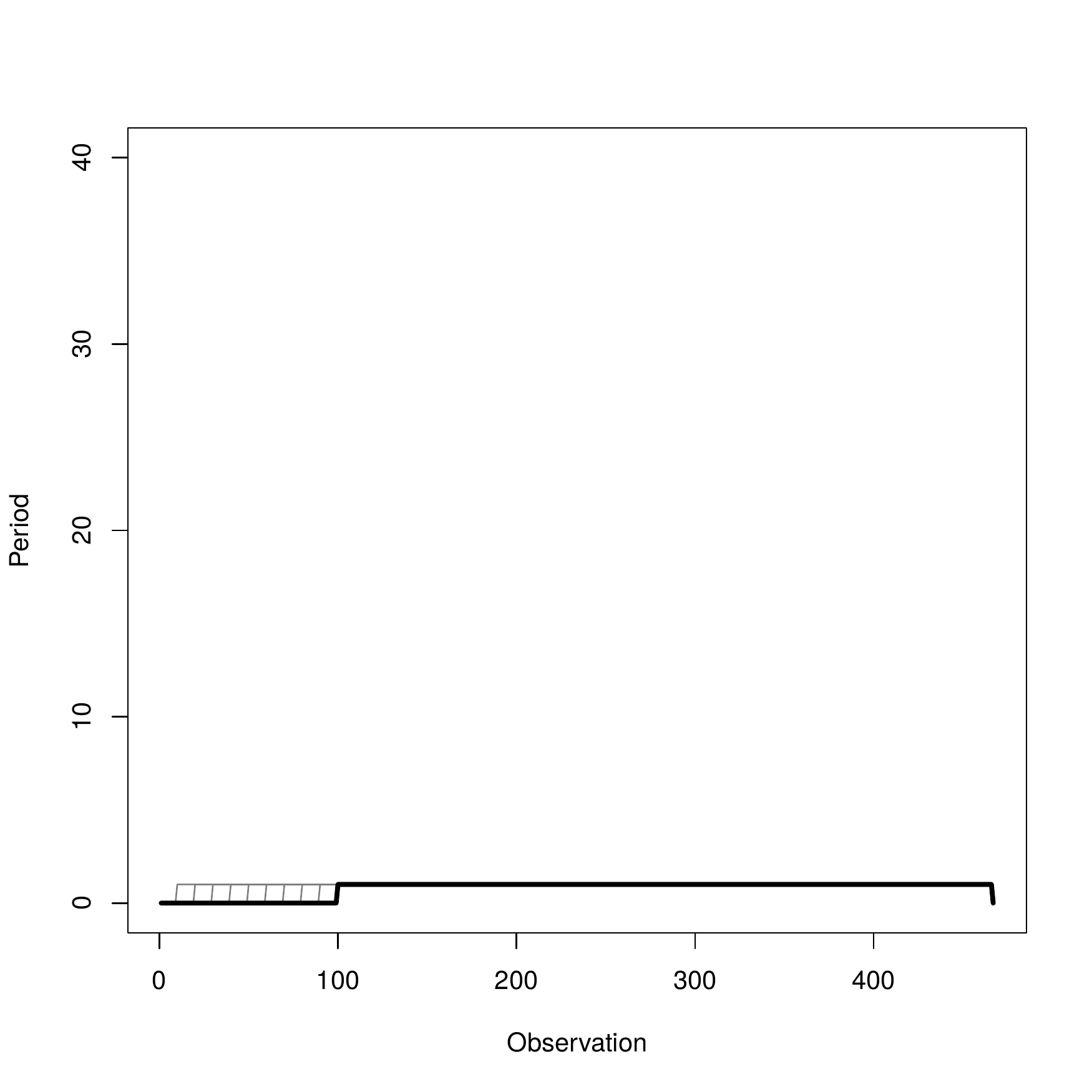}
    \end{center}
    \raggedright
    \footnotesize{Notes: Periods estimated on a rolling window basis using the \textit{findfrequency} from the R package \textit{forecast} \cite{hyndman2019package}. Solid line indicates the estimate from window size 100. All other window sizes added using thin lighter grey lines.}
\end{figure}

Table \ref{tab:shiphynd} and Figure \ref{fig:shiphynd} both illustrate the inability of the \textit{findfrequency} code from \cite{hyndman2019package} to identify any periodicity within the CASS data. In every window size we see a reported average periodicity of less than 1 with many observations recorded as 0. Whilst this is an empirical example and we have no ``true'' period against which to compare, the estimates from our TDA approach are much more in keeping with the feeling of periodicity from observing the time series. Hence another case is found where we are able to make significant improvements over widely implemented alternatives.

\subsection{Summary of Examples}

This section has offered three very different applications for the TDA approach we develop to be tested against the widely adopted spectral based methods. Specifically we have compared the \textit{findfrequency} function from the R package \textit{forecast} \cite{hyndman2019package}, with a TDA approach founded on the topological signatures of the time series. In all cases it was seen that TDA has strong potential to bring deeper appreciation of the periodicity of data. For Melbourne temperature data robustness to noise was key in TDA coming close to recognising the annual pattern; estimates close to 350 being in keeping with the 365 day year. In the Zurich sunspot data both approaches did well, and arguably the TDA was more sensitive to variation giving the spectral code a consistency advantage. However, in the CASS shipping data we saw again how spectral methods can fail to report any period at all even though the eye is telling us that there is a definite pattern there. The conclusion is thus very strongly in favour of undertaking TDA analysis of time series as a primary approach to identifying any cyclicality, or peridoicity, within data.

%\section{Multi varied time series}
%\label{sec:multi_varied_ts}
%In this paper we are concentrating on single valued time series, but all the presented methods works without any change in the multi varied case. As an example let us consider...

\section{Summary}
\label{sec:conclude}

In this paper we present a topological analysis based method to detect and quantify periodicity of noisy time series. The method is shown to be superior to the state of the art methods on a number of synthetic and real world examples. Such results are highly promising for the development of a deeper understanding of time series behaviour across the academic disciplines. For practice too the opportunity to learn more about time series is clear; manufacturers can understand machinery behaviour, business can learn of patterns in demand, and meterologists may learn more of climatic patterns. For ease of exposition this paper concentrates on univariate time series, but all results carry over to multivariate cases. Anywhere that has time series to analyse can benefit from the work presented here.

There are still many research questions left unanswered, not least the best way to operationalise the TDA approach within statistical software. The code used to perform the analysis for this paper is currently under development and it is intended that a public domain implementation of this method will be provided when the full version of the paper is published. There is scope to for applications in many directions and these will also offer avenues for future research. However, the most promising next steps lie in the detrending of data; how might we justifiably remove underlying trends to allow the TDA approach to find periodicity in time series. The \cite{hyndman2019package} \textit{findfrequency} approach offers linear detrending prior to estimation; this may be a helpful first step to add.

Notwithstanding all of these additional questions the tools developed in this paper are of great power and give scope to unlock many interesting avenues of further research in time series analysis. Primary advantages are the requirement of just two repeats of a cycle, one to identify, one to confirm, and the consistently demonstrated robustness to noise. Significant strides in understanding noisy time series are made to the chargrin of researchers and practitioners alike.

%Extracting periodicity from data with noise is lended an additional challenge by the multiple underlying cyclicalities that are present in many observed financial series. A theoretical approach rooted in TDA analysis is offered as a means to recover what the Fourier analysis based techniques cannot; namely the presence of changing periods in the build up to direction changes. In so doing many of the theories of technical analysis, such as Elliot waves, may be captured empirically and attributed neatly to the direction change they proport to capture. Our algorithm as presented has many other useful applications, including the potential to serve as an early warning system for impending crashes. Likewise moving to intraday and tick-by-tick data offers fruitful extension direction that can contribute to a growing literature on analysis of high-frequency data. Developing these applications is a challenge for future empirical work. This paper merely scratches the surface of the potential of TDA in time series analysis, as more is understood of the theory so the practicalities may be driven forwards. However, as a critical advance in the TDA literature, our work provides a practical robust contribution to their understanding.

\end{document}